\documentclass{tac}

\usepackage{mathrsfs,amssymb,stmaryrd,bbm,mathtools, upgreek}
\usepackage{chngcntr}
\usepackage[shortlabels]{enumitem}
\makeatletter
\@addtoreset{equation}{section}
\makeatother
\numberwithin{equation}{subsection}

\newtheorem{theo}{Theorem}[section]
\newtheorem{lem}[theo]{Lemma}
\newtheorem{prop}[theo]{Proposition}
\newtheorem{coro}[theo]{Corollary}
\newtheorem{defi}[theo]{Definition}
\newtheoremrm{rem}[theo]{Remark}

\let\olddefi\defi
\renewcommand{\defi}{\olddefi\normalfont}
\let\oldrem\rem
\renewcommand{\rem}{\oldrem\upshape}

\def\mathrmdef#1{\expandafter\def\csname#1\endcsname{{\rm#1}}}
\def\mathsfdef#1{\expandafter\def\csname#1\endcsname{{\rm\sf#1}}}
\def\mathcaldef#1{\expandafter\def\csname#1\endcsname{{\mathcal#1}}}

\mathrmdef{op}\mathrmdef{co}\mathrmdef{coop}
\mathsfdef{CoAlg}   \mathsfdef{Cat}
\mathsfdef{lim}
\mathsfdef{colim}
\mathsfdef{Alg}
\mathrmdef{id}
\mathsfdef{Top}
\mathsfdef{CmpHaus}
\mathsfdef{Set}
\mathsfdef{ran}
\mathsfdef{lan}
\mathsfdef{ord}

\mathsfdef{cat}

\def\Famf{{\mathsf{Fam}}_{\mathsf{fin}}}

\def\liva{\chi }
\def\livb{\gamma }
\def\livc{\zeta }

\def\!{\overline{!}}
\def\l!{\underline{!}}

\def\AAA{\mathbb{A}}
\def\CCC{\mathbb{C}}

\def\BBB{\mathbb{B}}

\def\TTTT{\mathcal{T}}
\def\RRRR{\mathcal{R}}

\usepackage{ifluatex}
\input diagxy        
\xyoption{curve}     
\xyoption{all}
\xyoption{2cell}

\begin{document}
	
\ifluatex
\directlua{adddednatlualoader = function ()
     require = function (stem)
         local fname = dednat6dir..stem..".lua"
         package.loaded[stem] = package.loaded[stem] or dofile(fname) or fname
       end
   end}
\catcode`\^^J=10
\directlua{dofile "dednat6load.lua"}
\else
%
\def\diagxyto{\ifnextchar/{\toop}{\toop/>/}}
\def\to     {\rightarrow}
\def\defded#1#2{\expandafter\def\csname ded-#1\endcsname{#2}}
\def\ifdedundefined#1{\expandafter\ifx\csname ded-#1\endcsname\relax}
\def\ded#1{\ifdedundefined{#1}
    \errmessage{UNDEFINED DEDUCTION: #1}
  \else
    \csname ded-#1\endcsname
  \fi
}
\def\defdiag#1#2{\expandafter\def\csname diag-#1\endcsname{\bfig#2\efig}}
\def\defdiagprep#1#2#3{\expandafter\def\csname diag-#1\endcsname{{#2\bfig#3\efig}}}
\def\ifdiagundefined#1{\expandafter\ifx\csname diag-#1\endcsname\relax}
\def\diag#1{\ifdiagundefined{#1}
    \errmessage{UNDEFINED DIAGRAM: #1}
  \else
    \csname diag-#1\endcsname
  \fi
}
\newlength{\celllower}
\newlength{\lcelllower}
\def\cellfont{}
\def\lcellfont{}
\def\cell #1{\lower\celllower\hbox to 0pt{\hss\cellfont${#1}$\hss}}
\def\lcell#1{\lower\celllower\hbox to 0pt   {\lcellfont${#1}$\hss}}
\def\expr#1{\directlua{output(tostring(#1))}}
\def\eval#1{\directlua{#1}}
\def\pu{\directlua{pu()}}
%

\defdiag{changeofthebaselaxcommacategorieschangeofbase}{   
  \morphism(0,0)|b|/{@{->}@/_18pt/}/<825,0>[{\AAA{//z}}`{\AAA{/y}};{c^\Leftarrow}]
  \morphism(825,0)|a|/{@{->}@/_18pt/}/<-825,0>[{\AAA{/y}}`{\AAA{//z}};{c\!}]
  \morphism(412,90)/-|/<0,-180>[{\phantom{O}}`{\phantom{O}};]
}
\defdiag{compositionoftwoadjunctionschangeofbasecommacompositionofpullbackcomma}{   
  \morphism(0,0)|b|/{@{->}@/_16pt/}/<675,0>[{\AAA{//z}}`{\AAA{/z}};{\id_z^\Leftarrow}]
  \morphism(675,0)|a|/{@{->}@/_16pt/}/<-675,0>[{\AAA{/z}}`{\AAA{//z}};{\id_z\!}]
  \morphism(675,0)|b|/{@{->}@/_16pt/}/<675,0>[{\AAA{/z}}`{\AAA{/y}};{c^\ast}]
  \morphism(1350,0)|a|/{@{->}@/_16pt/}/<-675,0>[{\AAA{/y}}`{\AAA{/z}};{c!}]
  \morphism(0,0)|b|/{@{->}@/_47pt/}/<1350,0>[{\AAA{//z}}`{\AAA{/y}};{c^\Leftarrow}]
  \morphism(1350,0)|a|/{@{->}@/_47pt/}/<-1350,0>[{\AAA{/y}}`{\AAA{//z}};{c\!}]
  \morphism(338,90)/-|/<0,-180>[{\phantom{O}}`{\phantom{O}};]
  \morphism(1012,90)/-|/<0,-180>[{\phantom{O}}`{\phantom{O}};]
}
\defdiag{whiskering1}{   
  \morphism(0,0)|r|/->/<0,-300>[{w}`{x};{f}]
  \morphism(0,-300)|l|/{@{->}@/_30pt/}/<0,-600>[{x}`{y};{h'}]
  \morphism(0,-300)|r|/{@{->}@/^30pt/}/<0,-600>[{x}`{y};{h}]
  \morphism(0,-900)|r|/->/<0,-300>[{y}`{z};{g}]
  \morphism(-225,-600)|a|/<=/<450,0>[{\phantom{O}}`{\phantom{O}};{\xi}]
}
\defdiag{whiskering2}{   
  \morphism(0,0)|r|/{@{->}@/_30pt/}/<0,-450>[{w}`{x};{f}]
  \morphism(0,0)|r|/{@{->}@/^30pt/}/<0,-450>[{w}`{x};{f}]
  \morphism(0,-450)|l|/{@{->}@/_30pt/}/<0,-600>[{x}`{y};{h}]
  \morphism(0,-450)|r|/{@{->}@/^30pt/}/<0,-600>[{x}`{y};{h'}]
  \morphism(0,-1050)|r|/{@{->}@/^30pt/}/<0,-450>[{y}`{z};{g}]
  \morphism(0,-1050)|r|/{@{->}@/_30pt/}/<0,-450>[{y}`{z};{g}]
  \morphism(-225,-750)|a|/=>/<450,0>[{\phantom{O}}`{\phantom{O}};{\xi}]
  \morphism(-90,-225)/=/<180,0>[{\phantom{O}}`{\phantom{O}};]
  \morphism(-90,-1275)/=/<180,0>[{\phantom{O}}`{\phantom{O}};]
}
\defdiag{twofoldtwofunctors}{   
  \morphism(0,0)|b|/{@{->}@/_16pt/}/<825,0>[{\AAA}`{\BBB};{G}]
  \morphism(825,0)|a|/{@{->}@/_16pt/}/<-825,0>[{\BBB}`{\AAA};{F}]
}
\defdiag{triangleidentityadjunctiondiagram1}{   
  \morphism(0,0)|a|/->/<600,0>[{\AAA}`{\BBB};{G}]
  \morphism(0,0)/=/<0,-600>[{\AAA}`{\AAA};]
  \morphism(600,0)/=/<0,-600>[{\BBB}`{\BBB};]
  \morphism(0,-600)|b|/->/<600,0>[{\AAA}`{\BBB};{G}]
  \morphism(600,0)|m|/->/<-600,-600>[{\BBB}`{\AAA};{F}]
  \morphism(225,-450)|a|/<=/<375,0>[{\phantom{O}}`{\phantom{O}};{\eta}]
  \morphism(0,-150)|a|/<=/<375,0>[{\phantom{O}}`{\phantom{O}};{\varepsilon}]
}
\defdiag{triangleidentityadjunctiondiagram2}{   
  \morphism(600,0)|a|/->/<-600,0>[{\BBB}`{\AAA};{F}]
  \morphism(0,0)/=/<0,-600>[{\AAA}`{\AAA};]
  \morphism(600,0)/=/<0,-600>[{\BBB}`{\BBB};]
  \morphism(600,-600)|b|/->/<-600,0>[{\BBB}`{\AAA};{F}]
  \morphism(0,-600)|m|/->/<600,600>[{\AAA}`{\BBB};{G}]
  \morphism(225,-450)|a|/<=/<375,0>[{\phantom{O}}`{\phantom{O}};{\eta}]
  \morphism(0,-150)|a|/<=/<375,0>[{\phantom{O}}`{\phantom{O}};{\varepsilon}]
}
\defdiag{basicadjunction}{   
  \morphism(0,0)|b|/{@{->}@/_16pt/}/<825,0>[{\AAA}`{\BBB};{G}]
  \morphism(825,0)|a|/{@{->}@/_16pt/}/<-825,0>[{\BBB}`{\AAA};{F}]
  \morphism(412,90)|r|/{@{-|}@<-7pt>}/<0,-180>[{\phantom{O}}`{\phantom{O}};{\left(\varepsilon{,}\eta\right)}]
}
\defdiag{leftsideoftheequationassociativityofmonad}{   
  \morphism(0,0)|a|/<-/<600,0>[{\BBB}`{\BBB};{T}]
  \morphism(600,-600)|r|/->/<0,600>[{\BBB}`{\BBB};{T}]
  \morphism(0,-600)|m|/->/<600,600>[{\BBB}`{\BBB};{T}]
  \morphism(600,-600)|b|/<-/<-600,0>[{\BBB}`{\BBB};{T}]
  \morphism(0,-600)|l|/->/<0,600>[{\BBB}`{\BBB};{T}]
  \morphism(0,-150)|a|/{@{<=}@<-5pt>}/<375,0>[{\phantom{O}}`{\phantom{O}};{\mu}]
  \morphism(225,-450)|a|/{@{<=}@<5pt>}/<375,0>[{\phantom{O}}`{\phantom{O}};{\mu}]
}
\defdiag{rightsideoftheequationassociativityofmonad}{   
  \morphism(0,0)|a|/<-/<600,0>[{\BBB}`{\BBB};{T}]
  \morphism(0,-600)|l|/->/<0,600>[{\BBB}`{\BBB};{T}]
  \morphism(600,-600)|b|/<-/<-600,0>[{\BBB}`{\BBB};{T}]
  \morphism(600,-600)|r|/->/<0,600>[{\BBB}`{\BBB};{T}]
  \morphism(600,-600)|m|/->/<-600,600>[{\BBB}`{\BBB};{T}]
  \morphism(0,-450)|a|/{@{<=}@<5pt>}/<375,0>[{\phantom{O}}`{\phantom{O}};{\mu}]
  \morphism(225,-150)|a|/{@{<=}@<-5pt>}/<375,0>[{\phantom{O}}`{\phantom{O}};{\mu}]
}
\defdiag{firstsideoftheequationidenityofamonad}{   
  \morphism(0,0)|m|/->/<600,0>[{\BBB}`{\BBB};{T}]
  \morphism(600,0)|r|/->/<0,-600>[{\BBB}`{\BBB};{T}]
  \morphism(0,0)|b|/->/<600,-600>[{\BBB}`{\BBB};{T}]
  \morphism(0,0)/{@{=}@/^40pt/}/<600,0>[{\BBB}`{\BBB};]
  \morphism(300,322)|r|/=>/<0,-300>[{\phantom{O}}`{\phantom{O}};{\eta}]
  \morphism(300,-38)|r|/{@{=>}@<10pt>}/<0,-330>[{\phantom{O}}`{\phantom{O}};{\mu}]
}
\defdiag{secondsideoftheequationidenityofamonad}{   
  \morphism(600,0)|m|/->/<-600,0>[{\BBB}`{\BBB};{T}]
  \morphism(600,-600)|r|/->/<0,600>[{\BBB}`{\BBB};{T}]
  \morphism(600,-600)|b|/->/<-600,600>[{\BBB}`{\BBB};{T}]
  \morphism(0,0)/{@{=}@/^40pt/}/<600,0>[{\BBB}`{\BBB};]
  \morphism(300,322)|r|/=>/<0,-300>[{\phantom{O}}`{\phantom{O}};{\eta}]
  \morphism(300,-38)|r|/{@{=>}@<10pt>}/<0,-330>[{\phantom{O}}`{\phantom{O}};{\mu}]
}
\defdiag{compositionofadjunctionsrarilali}{   
  \morphism(0,0)|b|/{@{->}@/_16pt/}/<900,0>[{w}`{x};{g}]
  \morphism(900,0)|a|/{@{->}@/_16pt/}/<-900,0>[{x}`{w};{f}]
  \morphism(900,0)|b|/{@{->}@/_16pt/}/<900,0>[{x}`{y};{g'}]
  \morphism(1800,0)|a|/{@{->}@/_16pt/}/<-900,0>[{y}`{x};{f'}]
  \morphism(450,90)|r|/{@{-|}@<-7pt>}/<0,-180>[{\phantom{O}}`{\phantom{O}};{\left({v}{,}n\right)}]
  \morphism(1350,90)|r|/{@{-|}@<-7pt>}/<0,-180>[{\phantom{O}}`{\phantom{O}};{\left(v'{,}n'\right)}]
}
\defdiag{proofodthelalicancellation}{   
  \morphism(750,-675)|b|/->/<750,-225>[{x}`{y};{g'}]
  \morphism(0,-450)|b|/->/<750,-225>[{w}`{x};{g}]
  \morphism(750,-225)|a|/->/<-750,-225>[{x}`{w};{f}]
  \morphism(1500,0)|a|/->/<-750,-225>[{y}`{x};{f'}]
  \morphism(1500,0)/=/<0,-900>[{y}`{y};]
  \morphism(750,-225)/=/<0,-450>[{x}`{x};]
  \morphism(938,-450)|a|/<=/<375,0>[{\phantom{O}}`{\phantom{O}};{n'}]
  \morphism(262,-450)|a|/<=/<375,0>[{\phantom{O}}`{\phantom{O}};{n}]
}
\defdiag{compositionofadjunctionsrarilalicorollaryforisomorphisms}{   
  \morphism(0,0)|b|/{@{->}@/_16pt/}/<900,0>[{w}`{x};{g'}]
  \morphism(900,0)|a|/{@{->}@/_16pt/}/<-900,0>[{x}`{w};{f'}]
  \morphism(900,0)|b|/{@{->}@/_16pt/}/<900,0>[{x}`{y};{g}]
  \morphism(1800,0)|a|/{@{->}@/_16pt/}/<-900,0>[{y}`{x};{f}]
  \morphism(1800,0)|b|/{@{->}@/_16pt/}/<900,0>[{y}`{z};{g''}]
  \morphism(2700,0)|a|/{@{->}@/_16pt/}/<-900,0>[{z}`{y};{f''}]
}
\defdiag{existingadjunction1}{   
  \morphism(0,0)|b|/{@{->}@/_16pt/}/<900,0>[{w}`{x};{g}]
  \morphism(900,0)|a|/{@{->}@/_16pt/}/<-900,0>[{x}`{w};{f}]
  \morphism(450,90)|r|/{@{-|}@<-7pt>}/<0,-180>[{\phantom{O}}`{\phantom{O}};{\left({v}{,}n\right)}]
}
\defdiag{existingadjunction2}{   
  \morphism(0,0)|b|/{@{->}@/_16pt/}/<900,0>[{w}`{y};{\hat{g}}]
  \morphism(900,0)|a|/{@{->}@/_16pt/}/<-900,0>[{y}`{w};{ff'}]
  \morphism(450,90)|r|/{@{-|}@<-7pt>}/<0,-180>[{\phantom{O}}`{\phantom{O}};{\left({\hat{v}}{,}\hat{n}\right)}]
}
\defdiag{counitofthecancellationtheoremlaris}{   
  \morphism(0,0)|a|/->/<600,0>[{x}`{w};{f}]
  \morphism(600,0)|a|/->/<600,0>[{w}`{y};{\hat{g}}]
  \morphism(1200,0)|r|/->/<0,-600>[{y}`{x};{f'}]
  \morphism(1200,-600)|m|/->/<-600,0>[{x}`{w};{f}]
  \morphism(600,-600)|m|/->/<-600,0>[{w}`{x};{g}]
  \morphism(1200,-600)/{@{=}@/^60pt/}/<-1200,0>[{x}`{x};]
  \morphism(0,0)/=/<0,-600>[{x}`{x};]
  \morphism(600,0)/=/<0,-600>[{w}`{w};]
  \morphism(712,-300)|a|/<=/<375,0>[{\phantom{O}}`{\phantom{O}};{\hat{v}}]
  \morphism(210,-300)/=/<180,0>[{\phantom{O}}`{\phantom{O}};]
  \morphism(600,-772)/=/<0,-180>[{\phantom{O}}`{\phantom{O}};]
}
\defdiag{exampleofnocancellationpropertyoflalis}{   
  \morphism(600,0)|b|/->/<-600,0>[{\mathsf{2}}`{\mathsf{1}};{s^0}]
  \morphism(1200,0)|a|/{@{->}@/_20pt/}/<-1200,0>[{\mathsf{1}}`{\mathsf{1}};{s^0d^0}]
}
\defdiag{pullbackdiagramdefinition}{   
  \morphism(0,0)|a|/->/<0,-600>[{x\times_{(a,b)}w}`{x};{a^\ast{(b)}}]
  \morphism(0,0)|a|/->/<600,0>[{x\times_{(a,b)}w}`{w};{b^\ast{(a)}}]
  \morphism(0,-600)|l|/->/<600,0>[{x}`{y};{a}]
  \morphism(600,0)|r|/->/<0,-600>[{w}`{y};{b}]
}
\defdiag{twocellofpullbackdefinitionleftside}{   
  \morphism(0,0)|l|/{@{->}@/_28pt/}/<0,-750>[{z}`{x};{h_0}]
  \morphism(0,0)|r|/{@{->}@/^28pt/}/<0,-750>[{z}`{x};{h_0'}]
  \morphism(900,0)|r|/->/<0,-750>[{w}`{y};{b}]
  \morphism(0,0)|a|/->/<900,0>[{z}`{w};{h_1'}]
  \morphism(0,-750)|r|/->/<900,0>[{x}`{y};{a}]
  \morphism(510,-375)/=/<180,0>[{\phantom{O}}`{\phantom{O}};]
  \morphism(-150,-375)|a|/=>/<300,0>[{\phantom{O}}`{\phantom{O}};{\xi_0}]
}
\defdiag{twocellpullbackdefinitionrightside}{   
  \morphism(0,0)|a|/{@{->}@/^28pt/}/<900,0>[{z}`{w};{h_1'}]
  \morphism(0,-750)|b|/->/<900,0>[{x}`{y};{a}]
  \morphism(0,0)|b|/{@{->}@/_28pt/}/<900,0>[{z}`{w};{h_1}]
  \morphism(0,0)|l|/->/<0,-750>[{z}`{x};{h_0}]
  \morphism(900,0)|r|/->/<0,-750>[{w}`{y};{b}]
  \morphism(450,-435)/=/<0,-180>[{\phantom{O}}`{\phantom{O}};]
  \morphism(450,150)|l|/<=/<0,-300>[{\phantom{O}}`{\phantom{O}};{\xi_1}]
}
\defdiag{commadiagramdefinition}{   
  \morphism(0,0)|a|/->/<600,0>[{a\downarrow{b}}`{x};{a^\Rightarrow{(b)}}]
  \morphism(0,0)|a|/->/<0,-600>[{a\downarrow{b}}`{w};{b^\Leftarrow{(a)}}]
  \morphism(600,0)|r|/->/<0,-600>[{x}`{y};{a}]
  \morphism(0,-600)|r|/->/<600,0>[{w}`{y};{b}]
  \morphism(112,-300)|a|/<=/<375,0>[{\phantom{O}}`{\phantom{O}};{\liva^{a\downarrow{b}}}]
}
\defdiag{morphismcommadiagramdefinition}{   
  \morphism(0,0)|r|/->/<375,-375>[{z}`{a\downarrow{b}};{h}]
  \morphism(375,-375)|a|/->/<600,0>[{a\downarrow{b}}`{x};{a^\Rightarrow{(b)}}]
  \morphism(375,-375)|a|/->/<0,-600>[{a\downarrow{b}}`{w};{b^\Leftarrow{(a)}}]
  \morphism(375,-975)|r|/->/<600,0>[{w}`{y};{b}]
  \morphism(975,-375)|r|/->/<0,-600>[{x}`{y};{a}]
  \morphism(1575,-375)|a|/->/<600,0>[{z}`{x};{h_0}]
  \morphism(1575,-375)|a|/->/<0,-600>[{z}`{w};{h_1}]
  \morphism(1575,-975)|r|/->/<600,0>[{w}`{y};{b}]
  \morphism(2175,-375)|r|/->/<0,-600>[{x}`{y};{a}]
  \morphism(488,-675)|a|/<=/<375,0>[{\phantom{O}}`{\phantom{O}};{\liva^{a\downarrow{b}}}]
  \morphism(1688,-675)|a|/<=/<375,0>[{\phantom{O}}`{\phantom{O}};{\livb}]
  \morphism(1185,-675)/=/<180,0>[{\phantom{O}}`{\phantom{O}};]
}
\defdiag{commatwocellofpullbackdefinitionleftside}{   
  \morphism(900,0)|r|/{@{->}@/^28pt/}/<0,-750>[{z}`{x};{h_0}]
  \morphism(900,0)|l|/{@{->}@/_28pt/}/<0,-750>[{z}`{x};{h_0'}]
  \morphism(0,0)|l|/->/<0,-750>[{w}`{y};{b}]
  \morphism(900,0)|a|/->/<-900,0>[{z}`{w};{h_1'}]
  \morphism(900,-750)|b|/->/<-900,0>[{x}`{y};{a}]
  \morphism(38,-375)|a|/<=/<525,0>[{\phantom{O}}`{\phantom{O}};{\liva^{a\downarrow{b}}\ast\id{_{h'}}}]
  \morphism(750,-375)|a|/<=/<300,0>[{\phantom{O}}`{\phantom{O}};{\xi_0}]
}
\defdiag{commatwocellpullbackdefinitionrightside}{   
  \morphism(900,0)|b|/{@{->}@/^28pt/}/<-900,0>[{z}`{w};{h_1}]
  \morphism(900,-750)|b|/->/<-900,0>[{x}`{y};{a}]
  \morphism(900,0)|a|/{@{->}@/_28pt/}/<-900,0>[{z}`{w};{h_1'}]
  \morphism(900,0)|r|/->/<0,-750>[{z}`{x};{h_0}]
  \morphism(0,0)|l|/->/<0,-750>[{w}`{y};{b}]
  \morphism(450,-375)|r|/<=/<0,-300>[{\phantom{O}}`{\phantom{O}};{\liva^{a\downarrow{b}}\ast\id{_{h}}}]
  \morphism(450,150)|r|/<=/<0,-300>[{\phantom{O}}`{\phantom{O}};{\xi_1}]
}
\defdiag{modificationofthedefinitionoflaxidempotenttwomonad}{   
  \morphism(0,0)|a|/=/<1050,0>[{T^2}`{T^2};{\id_{T^2}}]
  \morphism(0,0)|l|/=>/<525,-525>[{T^2}`{T};{\mu}]
  \morphism(525,-525)|r|/=>/<525,525>[{T}`{T^2};{\eta{T}}]
  \morphism(525,-38)|r|/=>/<0,-375>[{\phantom{O}}`{\phantom{O}};{\Gamma}]
}
\defdiag{modificationofthedefinitionoflaxidempotenttwomonadtriangleidentityoone}{   
  \morphism(0,0)/=/<750,0>[{T^2(z)}`{T^2(z)};]
  \morphism(0,0)|l|/->/<375,-375>[{T^2(z)}`{T(z)};{\mu_z}]
  \morphism(375,-375)|r|/->/<375,375>[{T(z)}`{T^2(z)};{\eta_{T(z)}}]
  \morphism(750,0)|a|/->/<450,0>[{T^2(z)}`{T(z)};{\mu_z}]
  \morphism(375,-15)|r|/=>/<0,-300>[{\phantom{O}}`{\phantom{O}};{\Gamma_z}]
}
\defdiag{modificationofthedefinitionoflaxidempotenttwomonadtriangleidentityotwo}{   
  \morphism(450,0)/=/<750,0>[{T^2(z)}`{T^2(z)};]
  \morphism(450,0)|l|/->/<375,-375>[{T^2(z)}`{T(z)};{\mu_z}]
  \morphism(825,-375)|r|/->/<375,375>[{T(z)}`{T^2(z)};{\eta_{T(z)}}]
  \morphism(0,0)|a|/->/<450,0>[{T(z)}`{T^2(z)};{\eta_{T(z)}}]
  \morphism(825,-15)|r|/=>/<0,-300>[{\phantom{O}}`{\phantom{O}};{\Gamma_z}]
}
\defdiag{equationforidempotentmonadinterchangelaw}{   
  \morphism(0,0)|l|/->/<0,-375>[{\AAA}`{\AAA};{T}]
  \morphism(0,-375)|r|/->/<600,-300>[{\AAA}`{\AAA};{T}]
  \morphism(600,-675)|r|/->/<-600,-300>[{\AAA}`{\AAA};{T}]
  \morphism(0,-375)|l|/->/<0,-600>[{\AAA}`{\AAA};{T}]
  \morphism(0,0)/{@{=}@/^45pt/}/<0,-375>[{\AAA}`{\AAA};]
  \morphism(75,-675)|a|/<=/<300,0>[{\phantom{O}}`{\phantom{O}};{\mu}]
  \morphism(-22,-188)|a|/<=/<300,0>[{\phantom{O}}`{\phantom{O}};{\eta}]
}
\defdiag{equationforidempotentmonadinterchangelawtwo}{   
  \morphism(0,0)|l|/->/<0,-375>[{\AAA}`{\AAA};{T}]
  \morphism(0,-375)|m|/->/<600,-300>[{\AAA}`{\AAA};{T}]
  \morphism(600,-675)|r|/->/<-600,-300>[{\AAA}`{\AAA};{T}]
  \morphism(0,-375)|l|/->/<0,-600>[{\AAA}`{\AAA};{T}]
  \morphism(0,-375)/{@{=}@/^35pt/}/<600,-300>[{\AAA}`{\AAA};]
  \morphism(75,-675)|a|/<=/<300,0>[{\phantom{O}}`{\phantom{O}};{\mu}]
  \morphism(300,-248)|l|/{@{=>}@<7pt>}/<0,-300>[{\phantom{O}}`{\phantom{O}};{\eta}]
}
\defdiag{pastinginordertogetmorphismofalgebrasvertical}{   
  \morphism(0,0)/=/<0,-510>[{T(x)}`{T(x)};]
  \morphism(0,-510)|a|/->/<0,-510>[{T(x)}`{T(y)};{T(f)}]
  \morphism(0,-1020)|r|/->/<1200,0>[{T(y)}`{y};{b}]
  \morphism(0,0)|a|/->/<1200,0>[{T(x)}`{x};{a}]
  \morphism(1200,0)|m|/->/<-1200,-510>[{x}`{T(x)};{\eta_x}]
  \morphism(1200,-510)|m|/->/<-1200,-510>[{y}`{T(y)};{\eta_y}]
  \morphism(1200,-510)/=/<0,-510>[{y}`{y};]
  \morphism(1200,0)|r|/->/<0,-510>[{x}`{y};{f}]
  \morphism(202,-128)/=/<195,0>[{\phantom{O}}`{\phantom{O}};]
  \morphism(802,-892)/=/<195,0>[{\phantom{O}}`{\phantom{O}};]
  \morphism(502,-510)/=/<195,0>[{\phantom{O}}`{\phantom{O}};]
}
\defdiag{equationforidempotentmonadinterchangelawtwoadjunctiontwo}{   
  \morphism(0,-900)/=/<0,-600>[{\AAA}`{\AAA};]
  \morphism(0,-900)|m|/->/<450,-300>[{\AAA}`{\BBB};{G}]
  \morphism(450,-1200)|b|/->/<-450,-300>[{\BBB}`{\AAA};{F}]
  \morphism(0,0)|l|/->/<0,-300>[{\BBB}`{\AAA};{F}]
  \morphism(0,-300)|l|/->/<0,-300>[{\AAA}`{\BBB};{G}]
  \morphism(0,-600)|l|/->/<0,-300>[{\BBB}`{\AAA};{F}]
  \morphism(0,-600)/{@{=}@/^20pt/}/<450,-600>[{\BBB}`{\BBB};]
  \morphism(52,-900)|a|/<=/<345,0>[{\phantom{O}}`{\phantom{O}};{\eta}]
  \morphism(0,-1200)|a|/<=/<330,0>[{\phantom{O}}`{\phantom{O}};{\varepsilon}]
}
\defdiag{equationforidempotentmonadinterchangelawtwoadjunction}{   
  \morphism(450,0)/=/<0,-600>[{\BBB}`{\BBB};]
  \morphism(450,-900)/=/<0,-600>[{\AAA}`{\AAA};]
  \morphism(450,-900)|a|/->/<450,-300>[{\AAA}`{\BBB};{G}]
  \morphism(900,-1200)|b|/->/<-450,-300>[{\BBB}`{\AAA};{F}]
  \morphism(450,0)|a|/->/<-450,-300>[{\BBB}`{\AAA};{F}]
  \morphism(0,-300)|b|/->/<450,-300>[{\AAA}`{\BBB};{G}]
  \morphism(450,-600)|r|/->/<0,-300>[{\BBB}`{\AAA};{F}]
  \morphism(120,-300)|a|/<=/<330,0>[{\phantom{O}}`{\phantom{O}};{\eta}]
  \morphism(450,-1200)|a|/<=/<330,0>[{\phantom{O}}`{\phantom{O}};{\varepsilon}]
}
\defdiag{equivalencesplitepifullyfaithful}{   
  \morphism(0,0)|a|/->/<900,0>[{\AAA{\left({w},x\right)}}`{\AAA{\left({FG(w)},x\right)}};{\AAA({\varepsilon_w},x)}]
  \morphism(900,0)|a|/->/<900,0>[{\AAA{\left({FG(w)},x\right)}}`{\BBB{\left({G(w)},G(x)\right)}};{\cong}]
  \morphism(0,0)|b|/{@{->}@/_20pt/}/<1800,0>[{\AAA{\left({w},x\right)}}`{\BBB{\left({G(w)},G(x)\right)}};{G}]
}
\defdiag{coequalizerofthecounitforBecktheorem}{   
  \morphism(0,0)|a|/{@{->}@/^20pt/}/<900,0>[{FGFG(x)}`{FG(x)};{\varepsilon_{FG(x)}}]
  \morphism(0,0)|b|/{@{->}@/_20pt/}/<900,0>[{FGFG(x)}`{FG(x)};{FG\left(\varepsilon_x\right)}]
  \morphism(900,0)|a|/->/<600,0>[{FG(x)}`{x};{\varepsilon_x}]
}
\defdiag{compositionof2adjunctionsfirstdiagram}{   
  \morphism(0,0)|b|/{@{->}@/_16pt/}/<900,0>[{\AAA}`{\BBB};{G}]
  \morphism(900,0)|a|/{@{->}@/_16pt/}/<-900,0>[{\BBB}`{\AAA};{F}]
  \morphism(900,0)|b|/{@{->}@/_16pt/}/<900,0>[{\BBB}`{\CCC};{J}]
  \morphism(1800,0)|a|/{@{->}@/_16pt/}/<-900,0>[{\CCC}`{\BBB};{H}]
  \Loop(900,0){{\BBB}}(ur,ul)_{\TTTT}
  \morphism(450,90)|r|/{@{-|}@<-7pt>}/<0,-180>[{\phantom{O}}`{\phantom{O}};{\left(\varepsilon{,}\eta\right)}]
  \morphism(1350,90)|r|/{@{-|}@<-7pt>}/<0,-180>[{\phantom{O}}`{\phantom{O}};{\left(\delta{,}\rho\right)}]
}
\defdiag{compositionof2adjunctionsseconddiagram}{   
  \morphism(0,0)|b|/{@{->}@/_20pt/}/<1800,0>[{\AAA}`{\CCC};{J\circ{G}}]
  \morphism(1800,0)|a|/{@{->}@/_20pt/}/<-1800,0>[{\CCC}`{\AAA};{F\circ{H}}]
  \Loop(1800,0){{\CCC}}(rd,ru)_\RRRR
  \morphism(900,112)|r|/{@{-|}@<-7pt>}/<0,-225>[{\phantom{O}}`{\phantom{O}};{\left(\varepsilon\cdot\left(F\delta{G}\right){,\,\,}\left(J\eta{H}\right)\cdot\rho\right)}]
}
\defdiag{compositionoftwoadjunctionscompactHausdorfftopset}{   
  \morphism(0,0)/{@{->}@/_16pt/}/<900,0>[{\CmpHaus}`{\Top};]
  \morphism(900,0)/{@{->}@/_16pt/}/<-900,0>[{\Top}`{\CmpHaus};]
  \morphism(900,0)/{@{->}@/_16pt/}/<900,0>[{\Top}`{\Set};]
  \morphism(1800,0)/{@{->}@/_16pt/}/<-900,0>[{\Set}`{\Top};]
  \morphism(450,90)/-|/<0,-180>[{\phantom{O}}`{\phantom{O}};]
  \morphism(1350,90)/-|/<0,-180>[{\phantom{O}}`{\phantom{O}};]
}
\defdiag{twocellidentityfortheproofofthecriterionforsimplicity}{   
  \morphism(0,0)|a|/->/<750,0>[{\AAA}`{\BBB};{G}]
  \morphism(750,0)|r|/->/<0,-450>[{\BBB}`{\CCC};{J}]
  \morphism(750,-450)/=/<0,-1200>[{\CCC}`{\CCC};]
  \morphism(750,-450)|m|/->/<-375,-300>[{\CCC}`{\BBB};{H}]
  \morphism(375,-750)|a|/->/<-375,-300>[{\BBB}`{\AAA};{F}]
  \morphism(0,-1050)|l|/->/<375,-300>[{\AAA}`{\BBB};{G}]
  \morphism(375,-1350)|b|/->/<375,-300>[{\BBB}`{\CCC};{J}]
  \morphism(750,0)/{@{=}@/_30pt/}/<-375,-750>[{\BBB}`{\BBB};]
  \morphism(375,-750)/=/<0,-600>[{\BBB}`{\BBB};]
  \morphism(75,-1050)|a|/<=/<300,0>[{\phantom{O}}`{\phantom{O}};{\eta}]
  \morphism(412,-1050)|a|/<=/<300,0>[{\phantom{O}}`{\phantom{O}};{\rho}]
  \morphism(398,-375)|a|/<=/<330,0>[{\phantom{O}}`{\phantom{O}};{\delta}]
}
\defdiag{compositionofadjunctionsepsilondelta}{   
  \morphism(0,0)|b|/{@{=>}@/_15pt/}/<1350,0>[{JG}`{JGFG};{J\eta{G}}]
  \morphism(1350,0)|a|/{@{=>}@/_15pt/}/<-1350,0>[{JGFG}`{JG};{JG\varepsilon}]
  \morphism(1350,0)|b|/{@{=>}@/_15pt/}/<1350,0>[{JGFG}`{JGFHJG};{\vartheta}]
  \morphism(2700,0)|a|/{@{=>}@/_15pt/}/<-1350,0>[{JGFHJG}`{JGFG};{JGF\delta{G}}]
  \morphism(0,0)|b|/{@{=>}@<-5pt>@/_40pt/}/<2700,0>[{JG}`{JGFHJG};{\alpha_{JG}}]
  \morphism(2700,0)|a|/{@{=>}@<-5pt>@/_40pt/}/<-2700,0>[{JGFHJG}`{JG};{JG\left(\varepsilon\left({F}\delta{G}\right)\right)}]
}
\defdiag{laxmorphismofcoalgebras1cell}{   
  \morphism(0,0)|a|/->/<300,0>[{w}`{x};{f}]
}
\defdiag{twocelloflaxcommamorphism}{   
  \morphism(0,0)|a|/->/<600,0>[{w}`{x};{f}]
  \morphism(0,0)|l|/->/<300,-420>[{w}`{y};{a}]
  \morphism(600,0)|r|/->/<-300,-420>[{x}`{y};{b}]
  \morphism(128,-0)|a|/{@{<=}@<-20pt>}/<345,0>[{\phantom{O}}`{\phantom{O}};{\livb}]
}
\defdiag{compositionofmorphismslaxcommatwocategories}{   
  \morphism(0,0)|a|/->/<675,0>[{w}`{x};{f}]
  \morphism(675,0)|a|/->/<675,0>[{x}`{z};{g}]
  \morphism(0,0)|l|/->/<675,-600>[{w}`{y};{a}]
  \morphism(1350,0)|r|/->/<-675,-600>[{z}`{y};{c}]
  \morphism(675,0)|m|/->/<0,-600>[{x}`{y};{b}]
  \morphism(240,-150)|a|/<=/<345,0>[{\phantom{O}}`{\phantom{O}};{\livb}]
  \morphism(765,-150)|a|/<=/<345,0>[{\phantom{O}}`{\phantom{O}};{\liva}]
}
\defdiag{leftsideequationtwocellforlaxcommacategor}{   
  \morphism(0,0)|a|/{@{->}@/^18pt/}/<750,0>[{w}`{x};{f}]
  \morphism(0,0)|b|/{@{->}@/_18pt/}/<750,0>[{w}`{x};{f'}]
  \morphism(0,0)|l|/->/<0,-675>[{w}`{y};{a}]
  \morphism(750,0)|r|/->/<-750,-675>[{x}`{y};{b}]
  \morphism(375,150)|l|/=>/<0,-300>[{\phantom{O}}`{\phantom{O}};{\livc}]
  \morphism(-8,-338)|a|/<=/<390,0>[{\phantom{O}}`{\phantom{O}};{\livb{'}}]
}
\defdiag{rightsideequationtwocellforlaxcommacategory}{   
  \morphism(0,0)|a|/->/<690,0>[{w}`{x};{f}]
  \morphism(0,0)|l|/->/<0,-645>[{w}`{y};{a}]
  \morphism(690,0)|r|/->/<-690,-645>[{x}`{y};{b}]
  \morphism(22,-322)|a|/{@{<=}@<15pt>}/<375,0>[{\phantom{O}}`{\phantom{O}};{\livb}]
}
\defdiag{changeofthebasecommacategories}{   
  \morphism(0,0)|b|/{@{->}@/_16pt/}/<825,0>[{\AAA{/z}}`{\AAA{/y}};{c^\ast}]
  \morphism(825,0)|a|/{@{->}@/_16pt/}/<-825,0>[{\AAA{/y}}`{\AAA{/z}};{c!}]
  \morphism(412,90)/-|/<0,-180>[{\phantom{O}}`{\phantom{O}};]
}
\defdiag{changeofthebasecommatwocategories}{   
  \morphism(0,0)|b|/{@{->}@/_16pt/}/<825,0>[{\AAA{/z}}`{\AAA{/y}};{c^\ast}]
  \morphism(825,0)|a|/{@{->}@/_16pt/}/<-825,0>[{\AAA{/y}}`{\AAA{/z}};{c!}]
  \morphism(412,90)/-|/<0,-180>[{\phantom{O}}`{\phantom{O}};]
}
\defdiag{cstaronmorphismm}{   
  \morphism(0,0)|l|/->/<0,-600>[{w\times_{(a,c)}y}`{x\times_{(b,c)}y};{c^\ast{\left(f\right)}}]
  \morphism(0,0)/->/<600,0>[{w\times_{(a,c)}y}`{w};]
  \morphism(0,-600)|m|/->/<600,0>[{x\times_{(b,c)}y}`{x};{b^\ast{(c)}}]
  \morphism(600,0)|r|/->/<0,-600>[{w}`{x};{f}]
  \morphism(600,-600)|r|/->/<0,-600>[{x}`{z};{b}]
  \morphism(0,-600)|l|/->/<0,-600>[{x\times_{(b,c)}y}`{y};{c^\ast{(b)}}]
  \morphism(0,-1200)|b|/->/<600,0>[{y}`{z};{c}]
  \morphism(600,0)|r|/{@{->}@/^60pt/}/<0,-1200>[{w}`{z};{a}]
  \morphism(0,0)|l|/{@{->}@/_60pt/}/<0,-1200>[{w\times_{(a,c)}y}`{y};{c^\ast(a)}]
}
\defdiag{rightrsidefirstequationimagetwocellofthepbchangeofbase}{   
  \morphism(0,0)|m|/->/<0,-600>[{w\times_{(a,c)}{y}}`{w};{a^\ast{(c)}}]
  \morphism(0,-600)|r|/{@{->}@/^22pt/}/<0,-600>[{w}`{x};{f}]
  \morphism(0,-600)|l|/{@{->}@/_22pt/}/<0,-600>[{w}`{x};{f'}]
  \morphism(-188,-900)|a|/<=/<375,0>[{\phantom{O}}`{\phantom{O}};{\livc}]
}
\defdiag{leftsidefirstequationimagetwocellofthepbchangeofbase}{   
  \morphism(0,-600)|m|/->/<0,-600>[{x\times_{(b,c)}{y}}`{x};{b^\ast{(c)}}]
  \morphism(0,0)|r|/{@{->}@/^22pt/}/<0,-600>[{w\times_{(a,c)}{y}}`{x\times_{(b,c)}{y}};{c^\ast{\left(f\right)}}]
  \morphism(0,0)|l|/{@{->}@/_22pt/}/<0,-600>[{w\times_{(a,c)}{y}}`{x\times_{(b,c)}{y}};{c^\ast{\left(f'\right)}}]
  \morphism(-188,-300)|a|/<=/<375,0>[{\phantom{O}}`{\phantom{O}};{c^\ast{\left(\livc\right)}}]
}
\defdiag{rightrsidesecondequationimagetwocellofthepbchangeofbase}{   
  \morphism(0,0)|r|/{@{->}@/^15pt/}/<0,-1200>[{w\times_{(a,c)}{y}}`{y};{c^\ast{(a)}}]
  \morphism(0,0)|l|/{@{->}@/_15pt/}/<0,-1200>[{w\times_{(a,c)}{y}}`{y};{c^\ast{(a)}}]
  \morphism(-90,-600)/=/<180,0>[{\phantom{O}}`{\phantom{O}};]
}
\defdiag{leftsidesecondequationimagetwocellofthepbchangeofbase}{   
  \morphism(0,-600)|m|/->/<0,-600>[{x\times_{(b,c)}{y}}`{y};{c^\ast{(b)}}]
  \morphism(0,0)|r|/{@{->}@/^22pt/}/<0,-600>[{w\times_{(a,c)}{y}}`{x\times_{(b,c)}{y}};{c^\ast{(f)}}]
  \morphism(0,0)|l|/{@{->}@/_22pt/}/<0,-600>[{w\times_{(a,c)}{y}}`{x\times_{(b,c)}{y}};{c^\ast{(f')}}]
  \morphism(-188,-300)|a|/<=/<375,0>[{\phantom{O}}`{\phantom{O}};{c^\ast{\left(\livc\right)}}]
}
\defdiag{directimagedefinitiondiagram}{   
  \morphism(1350,0)|a|/->/<-675,0>[{\AAA{/y}}`{\AAA{/z}};{c!}]
  \morphism(675,0)/->/<-675,0>[{\AAA{/z}}`{\AAA{//z}};]
  \morphism(1350,0)|m|/{@{->}@/_35pt/}/<-1350,0>[{\AAA{/y}}`{\AAA{//z}};{c\!}]
}
\defdiag{commadiagramdefinitioncommaadjunction}{   
  \morphism(0,0)|a|/->/<600,0>[{b\downarrow{c}}`{x};{b^\Rightarrow{(c)}}]
  \morphism(0,0)|l|/->/<0,-600>[{b\downarrow{c}}`{y};{c^\Leftarrow{(b)}}]
  \morphism(600,0)|r|/->/<0,-600>[{x}`{z};{b}]
  \morphism(0,-600)|b|/->/<600,0>[{y}`{z};{c}]
  \morphism(112,-300)|a|/<=/<375,0>[{\phantom{O}}`{\phantom{O}};{\liva^{b\downarrow{c}}}]
}
\defdiag{twocelloflaxcommamorphismmmm}{   
  \morphism(0,0)|a|/->/<540,0>[{w}`{x};{f}]
  \morphism(0,0)|l|/->/<270,-420>[{w}`{z};{a}]
  \morphism(540,0)|r|/->/<-270,-420>[{x}`{z};{b}]
  \morphism(98,-0)|a|/{@{<=}@<-20pt>}/<345,0>[{\phantom{O}}`{\phantom{O}};{\livb}]
}
\defdiag{morphismfforthedefinitioncLeftarrow}{   
  \morphism(0,0)|a|/->/<255,0>[{w}`{x};{f}]
}
\defdiag{cleftarrowonmorphismm}{   
  \morphism(0,0)|m|/->/<600,-600>[{a\downarrow{c}}`{b\downarrow{c}};{c^\Leftarrow{\left(f,\livb\right)}}]
  \morphism(600,-600)|m|/->/<600,0>[{b\downarrow{c}}`{x};{b^\Rightarrow{(c)}}]
  \morphism(1200,-600)|r|/->/<0,-600>[{x}`{z};{b}]
  \morphism(600,-600)|m|/->/<0,-600>[{b\downarrow{c}}`{y};{c^\Leftarrow{(b)}}]
  \morphism(600,-1200)|b|/->/<600,0>[{y}`{z};{c}]
  \morphism(0,0)|l|/{@{->}@/_30pt/}/<600,-1200>[{a\downarrow{c}}`{y};{c^\Leftarrow{(a)}}]
  \morphism(0,0)|a|/{@{->}@/^30pt/}/<1200,-600>[{a\downarrow{c}}`{x};{f\,\cdot\,{a^\Rightarrow{(c)}}}]
  \morphism(210,-600)/=/<180,0>[{\phantom{O}}`{\phantom{O}};]
  \morphism(712,-900)|a|/<=/<375,0>[{\phantom{O}}`{\phantom{O}};{\liva^{b\downarrow{c}}}]
  \morphism(510,-300)/=/<180,0>[{\phantom{O}}`{\phantom{O}};]
}
\defdiag{cleftarrowonmorphismmrightside}{   
  \morphism(1200,-600)|r|/->/<-600,-600>[{x}`{z};{b}]
  \morphism(0,-1200)|b|/->/<600,0>[{y}`{z};{c}]
  \morphism(0,0)|l|/->/<0,-1200>[{a\downarrow{c}}`{y};{c^\Leftarrow{(a)}}]
  \morphism(600,0)|r|/->/<600,-600>[{w}`{x};{f}]
  \morphism(600,0)|m|/->/<0,-1200>[{w}`{z};{a}]
  \morphism(0,0)|a|/->/<600,0>[{a\downarrow{c}}`{w};{a^\Rightarrow{(c)}}]
  \morphism(712,-600)|a|/<=/<375,0>[{\phantom{O}}`{\phantom{O}};{\livb}]
  \morphism(112,-600)|a|/<=/<375,0>[{\phantom{O}}`{\phantom{O}};{\liva^{a\downarrow{c}}}]
}
\defdiag{rightrsidefirstequationimagetwocellofthecommachangeofbase}{   
  \morphism(0,0)|m|/->/<0,-600>[{a\downarrow{c}}`{w};{a^\Rightarrow{(c)}}]
  \morphism(0,-600)|r|/{@{->}@/^22pt/}/<0,-600>[{w}`{x};{f}]
  \morphism(0,-600)|l|/{@{->}@/_22pt/}/<0,-600>[{w}`{x};{f'}]
  \morphism(-188,-900)|a|/<=/<375,0>[{\phantom{O}}`{\phantom{O}};{\livc}]
}
\defdiag{leftsidefirstequationimagetwocellofthecommachangeofbase}{   
  \morphism(0,-600)|m|/->/<0,-600>[{b\downarrow{c}}`{x};{b^\Rightarrow{(c)}}]
  \morphism(0,0)|r|/{@{->}@/^22pt/}/<0,-600>[{a\downarrow{c}}`{b\downarrow{c}};{c^\Leftarrow{\left(f,\livb\right)}}]
  \morphism(0,0)|l|/{@{->}@/_22pt/}/<0,-600>[{a\downarrow{c}}`{b\downarrow{c}};{c^\Leftarrow{\left(f',\livb{'}\right)}}]
  \morphism(-188,-300)|a|/<=/<375,0>[{\phantom{O}}`{\phantom{O}};{c^\Leftarrow{\left(\livc\right)}}]
}
\defdiag{rightrsidesecondequationimagetwocellofthecommachangeofbase}{   
  \morphism(0,0)|r|/{@{->}@/^15pt/}/<0,-1200>[{a\downarrow{c}}`{y};{c^\Leftarrow{(a)}}]
  \morphism(0,0)|l|/{@{->}@/_15pt/}/<0,-1200>[{a\downarrow{c}}`{y};{c^\Leftarrow{(a)}}]
  \morphism(-90,-600)/=/<180,0>[{\phantom{O}}`{\phantom{O}};]
}
\defdiag{leftsidesecondequationimagetwocellofthecommachangeofbase}{   
  \morphism(0,-600)|m|/->/<0,-600>[{b\downarrow{c}}`{y};{c^\Leftarrow{(b)}}]
  \morphism(0,0)|r|/{@{->}@/^22pt/}/<0,-600>[{a\downarrow{c}}`{b\downarrow{c}};{c^\Leftarrow{\left(f,\livb\right)}}]
  \morphism(0,0)|l|/{@{->}@/_22pt/}/<0,-600>[{a\downarrow{c}}`{b\downarrow{c}};{c^\Leftarrow{\left(f',\livb{'}\right)}}]
  \morphism(-188,-300)|a|/<=/<375,0>[{\phantom{O}}`{\phantom{O}};{c^\Leftarrow{\left(\livc\right)}}]
}
\defdiag{changeofthebaselaxcommacategories}{   
  \morphism(0,0)|b|/{@{->}@/_18pt/}/<825,0>[{\AAA{//z}}`{\AAA{/y}};{c^\Leftarrow}]
  \morphism(825,0)|a|/{@{->}@/_18pt/}/<-825,0>[{\AAA{/y}}`{\AAA{//z}};{c\!}]
  \morphism(412,90)/-|/<0,-180>[{\phantom{O}}`{\phantom{O}};]
}
\defdiag{commadiagramdefinitionproofcommaobject}{   
  \morphism(0,0)|a|/->/<600,0>[{b\downarrow{c}}`{x};{b^\Rightarrow{(c)}}]
  \morphism(0,0)|l|/->/<0,-600>[{b\downarrow{c}}`{y};{c^\Leftarrow{(b)}}]
  \morphism(600,0)|r|/->/<0,-600>[{x}`{z};{b}]
  \morphism(0,-600)|b|/->/<600,0>[{y}`{z};{c}]
  \morphism(112,-300)|a|/<=/<375,0>[{\phantom{O}}`{\phantom{O}};{\liva^{b\downarrow{c}}}]
}
\defdiag{commaalongca}{   
  \morphism(0,0)|a|/->/<600,0>[{ca\downarrow{c}}`{w};{\left({ca}\right)^\Rightarrow{(c)}}]
  \morphism(0,0)|l|/->/<0,-600>[{ca\downarrow{c}}`{y};{c^\Leftarrow{c\!}{(a)}}]
  \morphism(600,0)|r|/->/<0,-600>[{w}`{z};{ca}]
  \morphism(0,-600)|b|/->/<600,0>[{y}`{z};{c}]
  \morphism(112,-300)|a|/<=/<375,0>[{\phantom{O}}`{\phantom{O}};{\liva^{ca\downarrow{c}}}]
}
\defdiag{morphismcommadiagramdefinitionrholinha}{   
  \morphism(0,0)|r|/->/<375,-375>[{w}`{ca\downarrow{c}};{\rho_{(w,a)}'}]
  \morphism(375,-375)|a|/->/<600,0>[{ca\downarrow{c}}`{w};{(ca)^\Rightarrow{(c)}}]
  \morphism(375,-375)|l|/->/<0,-600>[{ca\downarrow{c}}`{y};{c^\Leftarrow{(ca)}}]
  \morphism(375,-975)|b|/->/<600,0>[{y}`{z};{c}]
  \morphism(975,-375)|r|/->/<0,-600>[{w}`{z};{ca}]
  \morphism(1575,-375)|a|/->/<600,0>[{w}`{w};{\id_w}]
  \morphism(1575,-375)|l|/->/<0,-600>[{w}`{y};{a}]
  \morphism(1575,-975)|b|/->/<600,0>[{y}`{z};{c}]
  \morphism(2175,-375)|r|/->/<0,-600>[{w}`{z};{ca}]
  \morphism(488,-675)|a|/<=/<375,0>[{\phantom{O}}`{\phantom{O}};{\liva^{ca\downarrow{c}}}]
  \morphism(1785,-675)/=/<180,0>[{\phantom{O}}`{\phantom{O}};]
  \morphism(1185,-675)/=/<180,0>[{\phantom{O}}`{\phantom{O}};]
}
\defdiag{definitionofdeltalinhaleft}{   
  \morphism(0,0)|m|/->/<600,-600>[{c\cdot{c^\Leftarrow{(b)}}\downarrow{c}}`{b\downarrow{c}};{\updelta{'}}]
  \morphism(600,-600)|m|/->/<600,0>[{b\downarrow{c}}`{x};{b^\Rightarrow{(c)}}]
  \morphism(1200,-600)|r|/->/<0,-600>[{x}`{z};{b}]
  \morphism(600,-600)|m|/->/<0,-600>[{b\downarrow{c}}`{y};{c^\Leftarrow{(b)}}]
  \morphism(600,-1200)|b|/->/<600,0>[{y}`{z};{c}]
  \morphism(0,0)|l|/{@{->}@/_30pt/}/<600,-1200>[{c\cdot{c^\Leftarrow{(b)}}\downarrow{c}}`{y};{c^\Leftarrow{(a)}}]
  \morphism(0,0)|a|/{@{->}@/^30pt/}/<1200,-600>[{c\cdot{c^\Leftarrow{(b)}}\downarrow{c}}`{x};{f\,\cdot\,{a^\Rightarrow{(c)}}}]
  \morphism(210,-600)/=/<180,0>[{\phantom{O}}`{\phantom{O}};]
  \morphism(712,-900)|a|/<=/<375,0>[{\phantom{O}}`{\phantom{O}};{\liva^{b\downarrow{c}}}]
  \morphism(510,-300)/=/<180,0>[{\phantom{O}}`{\phantom{O}};]
}
\defdiag{definitionofdeltalinharight}{   
  \morphism(1425,-600)|r|/->/<-600,-600>[{x}`{z};{b}]
  \morphism(0,-1200)|b|/->/<825,0>[{y}`{z};{c}]
  \morphism(0,0)|l|/->/<0,-1200>[{c\cdot{c^\Leftarrow{(b)}}\downarrow{c}}`{y};{c^\Leftarrow{c\!}{c^\Leftarrow}(b)}]
  \morphism(825,0)|r|/->/<600,-600>[{b\downarrow{c}}`{x};{b^\Rightarrow{(c)}}]
  \morphism(825,0)|m|/->/<0,-600>[{b\downarrow{c}}`{y};{c^\Leftarrow{(b)}}]
  \morphism(825,-600)|m|/->/<0,-600>[{y}`{z};{c}]
  \morphism(0,0)|a|/->/<825,0>[{c\cdot{c^\Leftarrow{(b)}}\downarrow{c}}`{b\downarrow{c}};{\left({c}\cdot{c^\Leftarrow}(b)\right)^\Rightarrow{(c)}}]
  \morphism(938,-600)|a|/<=/<375,0>[{\phantom{O}}`{\phantom{O}};{\liva^{b\downarrow{c}}}]
  \morphism(150,-600)|a|/<=/<525,0>[{\phantom{O}}`{\phantom{O}};{\liva^{c\cdot{c^\Leftarrow}(b)\downarrow{c}}}]
}
\defdiag{secontriangleidentitycoomaobjecttwoadjunction}{   
  \morphism(0,0)|r|/->/<300,-300>[{b\downarrow{c}}`{c\cdot{c}^\Leftarrow{(b)}\downarrow{c}};{\rho_{c^\Leftarrow{(x,b)}}'}]
  \morphism(300,-300)|r|/->/<300,-300>[{c\cdot{c}^\Leftarrow{(b)}\downarrow{c}}`{b\downarrow{c}};{\updelta{'}}]
  \morphism(600,-600)|a|/->/<600,0>[{b\downarrow{c}}`{x};{b^\Rightarrow{(c)}}]
  \morphism(600,-600)|l|/->/<0,-600>[{b\downarrow{c}}`{y};{c^\Leftarrow{(b)}}]
  \morphism(600,-1200)|r|/->/<600,0>[{y}`{z};{c}]
  \morphism(1200,-600)|r|/->/<0,-600>[{x}`{z};{b}]
  \morphism(2100,-600)|a|/->/<600,0>[{b\downarrow{c}}`{x};{b^\Rightarrow{(c)}}]
  \morphism(2100,-600)|l|/->/<0,-600>[{b\downarrow{c}}`{y};{c^\Leftarrow{(b)}}]
  \morphism(2100,-1200)|r|/->/<600,0>[{y}`{z};{c}]
  \morphism(2700,-600)|r|/->/<0,-600>[{x}`{z};{b}]
  \morphism(712,-900)|a|/<=/<375,0>[{\phantom{O}}`{\phantom{O}};{\liva^{b\downarrow{c}}}]
  \morphism(2212,-900)|a|/<=/<375,0>[{\phantom{O}}`{\phantom{O}};{\liva^{b\downarrow{c}}}]
  \morphism(1560,-900)/=/<180,0>[{\phantom{O}}`{\phantom{O}};]
}
\defdiag{basicdiagramofcompositionwithcandinclusion}{   
  \morphism(1050,0)|b|/->/<-525,0>[{\AAA{/y}}`{\AAA{/z}};{c!}]
  \morphism(525,0)|b|/->/<-525,0>[{\AAA{/z}}`{\AAA{//z}};{\id_z\!}]
  \morphism(1050,0)|a|/{@{->}@/_30pt/}/<-1050,0>[{\AAA{/y}}`{\AAA{//z}};{c\!}]
}
\defdiag{compositionoftwoadjunctionschangeofbasecomma}{   
  \morphism(0,0)|b|/{@{->}@/_16pt/}/<675,0>[{\AAA{//z}}`{\AAA{/z}};{\id_z^\Leftarrow}]
  \morphism(675,0)|a|/{@{->}@/_16pt/}/<-675,0>[{\AAA{/z}}`{\AAA{//z}};{\id_z\!}]
  \morphism(675,0)|b|/{@{->}@/_16pt/}/<675,0>[{\AAA{/z}}`{\AAA{/y}};{c^\ast}]
  \morphism(1350,0)|a|/{@{->}@/_16pt/}/<-675,0>[{\AAA{/y}}`{\AAA{/z}};{c!}]
  \morphism(0,0)|b|/{@{->}@/_47pt/}/<1350,0>[{\AAA{//z}}`{\AAA{/y}};{c^\Leftarrow}]
  \morphism(1350,0)|a|/{@{->}@/_47pt/}/<-1350,0>[{\AAA{/y}}`{\AAA{//z}};{c\!}]
  \morphism(338,90)/-|/<0,-180>[{\phantom{O}}`{\phantom{O}};]
  \morphism(1012,90)/-|/<0,-180>[{\phantom{O}}`{\phantom{O}};]
}
\defdiag{changeofthebaseidentitycoherencelaxidempotent}{   
  \morphism(0,0)|b|/{@{->}@/_18pt/}/<825,0>[{\AAA{//y}}`{\AAA{/y}};{\id_y^\Leftarrow}]
  \morphism(825,0)|a|/{@{->}@/_18pt/}/<-825,0>[{\AAA{/y}}`{\AAA{//y}};{\id_y\!}]
  \morphism(412,90)|r|/-|/<0,-180>[{\phantom{O}}`{\phantom{O}};{\left(\delta{,}\rho\right)}]
}
\defdiag{counittwocelloftheidentitytwoadjunctionlaxcoomacomma}{   
  \morphism(0,0)|a|/->/<675,0>[{b\downarrow{\id_y}}`{x};{b^\Rightarrow{(\id_y)}}]
  \morphism(675,0)|r|/->/<0,-675>[{x}`{y};{b}]
  \morphism(0,0)|l|/->/<0,-675>[{b\downarrow{\id_y}}`{y};{\id_y^\Leftarrow{(b)}}]
  \morphism(0,-675)|b|/->/<675,0>[{y}`{y};{\id_y}]
  \morphism(128,-338)|a|/<=/<420,0>[{\phantom{O}}`{\phantom{O}};{\liva^{b\downarrow{\id_y}}}]
}
\defdiag{rholparaademonstracaodelaxidempotency}{   
  \morphism(0,0)|m|/->/<900,-450>[{x}`{b\downarrow{\id_y}};{\overline{\underline{\rho}}_{(x,b)}}]
  \morphism(900,-450)|m|/->/<900,0>[{b\downarrow{\id_y}}`{x};{b^\Rightarrow{(\id_y)}\,{=}\,\overline{\underline{\delta}}_{(x,b)}}]
  \morphism(1800,-450)|r|/->/<0,-900>[{x}`{y};{b}]
  \morphism(900,-450)|m|/->/<0,-900>[{b\downarrow{\id_y}}`{y};{\id_y^\Leftarrow{(b)}}]
  \morphism(900,-1350)|r|/->/<900,0>[{y}`{y};{\id_y}]
  \morphism(0,0)|l|/{@{->}@/_30pt/}/<900,-1350>[{x}`{y};{b}]
  \morphism(0,0)|r|/{@{->}@/^30pt/}/<1800,-450>[{x}`{x};{\id_x}]
  \morphism(360,-675)/=/<180,0>[{\phantom{O}}`{\phantom{O}};]
  \morphism(1125,-900)|a|/<=/<450,0>[{\phantom{O}}`{\phantom{O}};{\liva^{b\downarrow{\id_y}}}]
  \morphism(810,-225)/=/<180,0>[{\phantom{O}}`{\phantom{O}};]
}
\defdiag{identityonbtodefinerho}{   
  \morphism(0,0)|r|/{@{->}@/^28pt/}/<0,-1350>[{x}`{y};{b}]
  \morphism(0,0)|l|/{@{->}@/_28pt/}/<0,-1350>[{x}`{y};{b}]
  \morphism(-90,-675)/=/<180,0>[{\phantom{O}}`{\phantom{O}};]
}
\defdiag{firsttwocellinordertodefinecounitoftheadjunctionthatgivesarariii}{   
  \morphism(0,0)|a|/->/<300,-450>[{b\downarrow{\id_y}}`{x};{\overline{\underline{\delta}}_{(x,b)}}]
  \morphism(0,0)|b|/{@{->}@/_60pt/}/<600,-1500>[{b\downarrow{\id_y}}`{y};{\id_y^\Leftarrow{(b)}}]
  \morphism(300,-450)|a|/->/<300,-450>[{x}`{b\downarrow{\id_y}};{\overline{\underline{\rho}}_{(x,b)}}]
  \morphism(600,-900)|m|/->/<0,-600>[{b\downarrow{\id_y}}`{y};{\id_y^\Leftarrow{(b)}}]
  \morphism(300,-450)|l|/{@{->}@/_30pt/}/<300,-1050>[{x}`{y};{b}]
  \morphism(285,-975)/=/<180,0>[{\phantom{O}}`{\phantom{O}};]
  \morphism(-150,-750)|a|/<=/<300,0>[{\phantom{O}}`{\phantom{O}};{\liva^{b\downarrow{\id_y}}}]
}
\defdiag{secondtwocellinordertodefinecounitoftheadjunctionthatgivesarariiii}{   
  \morphism(0,0)|l|/->/<600,-900>[{b\downarrow{\id_y}}`{b\downarrow{\id_y}};{\id_{b\downarrow{\id_y}}}]
  \morphism(0,0)|a|/{@{->}@/^40pt/}/<1200,-900>[{b\downarrow{\id_y}}`{x};{\overline{\underline{\delta}}_{(x,b)}}]
  \morphism(600,-900)|a|/->/<600,0>[{b\downarrow{\id_y}}`{x};{\overline{\underline{\delta}}_{(x,b)}}]
  \morphism(510,-450)/=/<180,0>[{\phantom{O}}`{\phantom{O}};]
}
\defdiag{firsttwocellinordertodefinecounitoftheadjunctionthatgivesarari}{   
  \morphism(0,0)|a|/->/<300,-450>[{b\downarrow{\id_y}}`{x};{\overline{\underline{\delta}}_{(x,b)}}]
  \morphism(0,0)|b|/{@{->}@/_60pt/}/<600,-1500>[{b\downarrow{\id_y}}`{y};{\id_y^\Leftarrow{(b)}}]
  \morphism(300,-450)|a|/->/<300,-450>[{x}`{b\downarrow{\id_y}};{\overline{\underline{\rho}}_{(x,b)}}]
  \morphism(600,-900)|a|/->/<600,0>[{b\downarrow{\id_y}}`{x};{\overline{\underline{\delta}}_{(x,b)}}]
  \morphism(1200,-900)|r|/->/<0,-600>[{x}`{y};{b}]
  \morphism(600,-900)|m|/->/<0,-600>[{b\downarrow{\id_y}}`{y};{\id_y^\Leftarrow{(b)}}]
  \morphism(600,-1500)|b|/->/<600,0>[{y}`{y};{\id_y}]
  \morphism(300,-450)|l|/{@{->}@/_30pt/}/<300,-1050>[{x}`{y};{b}]
  \morphism(285,-975)/=/<180,0>[{\phantom{O}}`{\phantom{O}};]
  \morphism(750,-1200)|a|/<=/<300,0>[{\phantom{O}}`{\phantom{O}};{\liva^{b\downarrow{\id_y}}}]
  \morphism(-150,-750)|a|/<=/<300,0>[{\phantom{O}}`{\phantom{O}};{\liva^{b\downarrow{\id_y}}}]
}
\defdiag{secondtwocellinordertodefinecounitoftheadjunctionthatgivesarari}{   
  \morphism(0,0)|l|/->/<600,-900>[{b\downarrow{\id_y}}`{b\downarrow{\id_y}};{\id_{b\downarrow{\id_y}}}]
  \morphism(0,0)|a|/{@{->}@/^40pt/}/<1200,-900>[{b\downarrow{\id_y}}`{x};{\overline{\underline{\delta}}_{(x,b)}}]
  \morphism(600,-900)|a|/->/<600,0>[{b\downarrow{\id_y}}`{x};{\overline{\underline{\delta}}_{(x,b)}}]
  \morphism(1200,-900)|r|/->/<0,-600>[{x}`{y};{b}]
  \morphism(600,-900)|m|/->/<0,-600>[{b\downarrow{\id_y}}`{y};{\id_y^\Leftarrow{(b)}}]
  \morphism(600,-1500)|b|/->/<600,0>[{y}`{y};{\id_y}]
  \morphism(750,-1200)|a|/<=/<300,0>[{\phantom{O}}`{\phantom{O}};{\liva^{b\downarrow{\id_y}}}]
  \morphism(510,-450)/=/<180,0>[{\phantom{O}}`{\phantom{O}};]
}
\defdiag{basicliftingtwoadjunction}{   
  \morphism(0,0)|b|/{@{->}@/_16pt/}/<825,0>[{\AAA{/y}}`{\BBB{/G(y)}};{\check{G}}]
  \morphism(825,0)|a|/{@{->}@/_16pt/}/<-825,0>[{\BBB{/G(y)}}`{\AAA{/y}};{\varepsilon_{y}!\circ\check{F}}]
  \morphism(412,90)/-|/<0,-180>[{\phantom{O}}`{\phantom{O}};]
}
\defdiag{basicliftingtwoadjunctionlaxcommatwocategory}{   
  \morphism(0,0)|b|/{@{->}@/_16pt/}/<825,0>[{\AAA{//y}}`{\BBB{//G(y)}};{\check{G}}]
  \morphism(825,0)|a|/{@{->}@/_16pt/}/<-825,0>[{\BBB{//G(y)}}`{\AAA{//y}};{\varepsilon_{y}\l!\circ\check{F}}]
  \morphism(412,90)/-|/<0,-180>[{\phantom{O}}`{\phantom{O}};]
}
\defdiag{basicliftingtwoadjunctionn}{   
  \morphism(0,0)|b|/{@{->}@/_16pt/}/<825,0>[{\AAA{/y}}`{\BBB{/G(y)}};{\check{G}}]
  \morphism(825,0)|a|/{@{->}@/_16pt/}/<-825,0>[{\BBB{/G(y)}}`{\AAA{/y}};{\varepsilon_{y}!\circ\check{F}}]
  \morphism(412,90)/-|/<0,-180>[{\phantom{O}}`{\phantom{O}};]
}
\defdiag{basicliftingtwoadjunctionlaxcommatwocategoryy}{   
  \morphism(0,0)|b|/{@{->}@/_16pt/}/<825,0>[{\AAA{//y}}`{\BBB{//G(y)}};{\check{G}}]
  \morphism(825,0)|a|/{@{->}@/_16pt/}/<-825,0>[{\BBB{//G(y)}}`{\AAA{//y}};{\varepsilon_{y}\l!\circ\check{F}}]
  \morphism(412,90)/-|/<0,-180>[{\phantom{O}}`{\phantom{O}};]
}
\defdiag{basicliftingtwoadjunctionlaxcommatwocategoryagain}{   
  \morphism(0,0)|b|/{@{->}@/_16pt/}/<825,0>[{\AAA{//y}}`{\BBB{//G(y)}};{\check{G}}]
  \morphism(825,0)|a|/{@{->}@/_16pt/}/<-825,0>[{\BBB{//G(y)}}`{\AAA{//y}};{\varepsilon_{y}\!\circ\check{F}}]
  \morphism(412,90)/-|/<0,-180>[{\phantom{O}}`{\phantom{O}};]
}
\defdiag{basicadjunctionliftingagain}{   
  \morphism(0,0)|b|/{@{->}@/_16pt/}/<825,0>[{\AAA{/y}}`{\BBB{/G(y)}};{\check{G}}]
  \morphism(825,0)|a|/{@{->}@/_16pt/}/<-825,0>[{\BBB{/G(y)}}`{\AAA{/y}};{\varepsilon_{y}!\circ\check{F}}]
  \morphism(412,90)/-|/<0,-180>[{\phantom{O}}`{\phantom{O}};]
}
\defdiag{compositionof2adjunctionsfirstdiagramlaxcommatwocategorycomma}{   
  \morphism(0,0)|b|/{@{->}@/_20pt/}/<1200,0>[{\AAA{//F(y)}}`{\BBB{//GF(y)}};{\check{G}}]
  \morphism(1200,0)|a|/{@{->}@/_20pt/}/<-1200,0>[{\BBB{//GF(y)}}`{\AAA{//F(y)}};{\varepsilon_{F(y)}\l!\circ\check{F}}]
  \morphism(1200,0)|b|/{@{->}@/_20pt/}/<1200,0>[{\BBB{//GF(y)}}`{\BBB{/y}};{\eta_y^\Leftarrow}]
  \morphism(2400,0)|a|/{@{->}@/_20pt/}/<-1200,0>[{\BBB{/y}}`{\BBB{//GF(y)}};{\eta_y\!}]
  \Loop(1200,0){{\BBB{//GF(y)}}}(ur,ul)_{\TTTT}
  \morphism(600,90)|r|/{@{-|}@<-7pt>}/<0,-180>[{\phantom{O}}`{\phantom{O}};{\left(\varepsilon{,}\eta\right)}]
  \morphism(1800,90)|r|/{@{-|}@<-7pt>}/<0,-180>[{\phantom{O}}`{\phantom{O}};{\left(\delta{,}\rho\right)}]
}
\defdiag{compositionof2adjunctionsseconddiagrametacommalaxcomma}{   
  \morphism(0,0)|b|/{@{->}@/_23pt/}/<2400,0>[{\AAA{//F(y)}}`{\BBB{/y}};{\eta_y^\Leftarrow\circ\check{G}}]
  \morphism(2400,0)|a|/{@{->}@/_23pt/}/<-2400,0>[{\BBB{/y}}`{\AAA{//F(y)}};{\check{F}}]
  \Loop(2400,0){{\BBB{/y}}}(rd,ru)_\RRRR
  \morphism(1200,112)|r|/{@{-|}@<-15pt>}/<0,-225>[{\phantom{O}}`{\phantom{O}};{\left(\varepsilon\cdot\left(\id_{\check{F}}\ast\delta\ast\id_{\check{G}}\right){,}\,\,\,\alpha\right)}]
}
\defdiag{definitionofalphalinha}{   
  \morphism(0,0)|r|/->/<450,-450>[{x}`{GF(b)\downarrow{\eta_y}};{\alpha_{(x,b)}}]
  \morphism(450,-450)|a|/->/<900,0>[{GF(b)\downarrow{\eta_y}}`{GF(x)};{\left({G}F(b)\right)^\Rightarrow\left(\eta_y\right)}]
  \morphism(450,-450)|a|/->/<0,-750>[{GF(b)\downarrow{\eta_y}}`{y};{\eta_y^\Leftarrow\left({G}F(b)\right)}]
  \morphism(450,-1200)|b|/->/<900,0>[{y}`{GF(y)};{\eta_y}]
  \morphism(1350,-450)|r|/->/<0,-750>[{GF(x)}`{GF(y)};{GF(b)}]
  \morphism(2100,-450)|a|/->/<900,0>[{x}`{GF(x)};{\eta_x}]
  \morphism(2100,-450)|l|/->/<0,-750>[{x}`{y};{b}]
  \morphism(2100,-1200)|b|/->/<900,0>[{y}`{GF(y)};{\eta_y}]
  \morphism(3000,-450)|r|/->/<0,-750>[{GF(x)}`{GF(y)};{GF(b)}]
  \morphism(675,-825)|a|/<=/<450,0>[{\phantom{O}}`{\phantom{O}};{\liva^{ca\downarrow{c}}}]
  \morphism(2445,-825)/=/<210,0>[{\phantom{O}}`{\phantom{O}};]
  \morphism(1672,-825)/=/<210,0>[{\phantom{O}}`{\phantom{O}};]
}
\defdiag{basicliftingtwoadjunctionlaxcommatwocategoryylefttwoadjoint}{   
  \morphism(675,0)|b|/->/<-675,0>[{\BBB{//GF(y)}}`{\AAA{//F(y)}};{\varepsilon_{F(y)}\l!\circ\check{F}}]
  \morphism(1350,0)|b|/->/<-675,0>[{\BBB{/y}}`{\BBB{//GF(y)}};{\eta_y\!}]
  \morphism(1350,0)|m|/{@{->}@/_28pt/}/<-1350,0>[{\BBB{/y}}`{\AAA{//F(y)}};{\check{F}}]
}
\defdiag{basicliftingtwoadjunctionpullbacktwocategoryylefttwoadjoint}{   
  \morphism(675,0)|b|/->/<-675,0>[{\BBB{/GF(y)}}`{\AAA{/F(y)}};{\varepsilon_{F(y)}!\circ\check{F}}]
  \morphism(1350,0)|b|/->/<-675,0>[{\BBB{/y}}`{\BBB{/GF(y)}};{\eta_y!}]
  \morphism(1350,0)|m|/{@{->}@/_28pt/}/<-1350,0>[{\BBB{/y}}`{\AAA{/F(y)}};{\check{F}}]
}
\defdiag{2adjunctionidentityforthesouthafricantheorem}{   
  \morphism(0,0)|b|/{@{->}@/_16pt/}/<825,0>[{\AAA{//F(y)}}`{\AAA{/F(y)}};{\id_{F(y)}^{\Leftarrow}}]
  \morphism(825,0)|a|/{@{->}@/_16pt/}/<-825,0>[{\AAA{/F(y)}}`{\AAA{//F(y)}};{\id_{F(y)}\!}]
  \morphism(412,90)/-|/<0,-180>[{\phantom{O}}`{\phantom{O}};]
}
\defdiag{compositionof2adjunctionsadmissiblewrtbasifibration}{   
  \morphism(0,0)|b|/{@{->}@/_20pt/}/<1200,0>[{\AAA{/F(y)}}`{\BBB{/GF(y)}};{\check{G}}]
  \morphism(1200,0)|a|/{@{->}@/_20pt/}/<-1200,0>[{\BBB{/GF(y)}}`{\AAA{/F(y)}};{\varepsilon_{F(y)}!\circ\check{F}}]
  \morphism(1200,0)|b|/{@{->}@/_20pt/}/<1200,0>[{\BBB{/GF(y)}}`{\BBB{/y}};{\eta_y^\ast}]
  \morphism(2400,0)|a|/{@{->}@/_20pt/}/<-1200,0>[{\BBB{/y}}`{\BBB{/GF(y)}};{\eta_y!}]
  \morphism(0,0)|b|/{@{->}@<-5pt>@/_45pt/}/<2400,0>[{\AAA{/F(y)}}`{\BBB{/y}};{\eta_y^{\ast}\,\circ\,\check{G}}]
  \morphism(2400,0)|a|/{@{->}@<-5pt>@/_45pt/}/<-2400,0>[{\BBB{/y}}`{\AAA{/F(y)}};{\check{F}}]
  \morphism(600,90)|r|/{@{-|}@<-7pt>}/<0,-180>[{\phantom{O}}`{\phantom{O}};{\left(\varepsilon{,}\eta\right)}]
  \morphism(1800,90)/-|/<0,-180>[{\phantom{O}}`{\phantom{O}};]
}
\defdiag{theoremthatusespreservationofcommafirstleft}{   
  \morphism(0,0)|r|/{@{->}@/^20pt/}/<0,-900>[{\AAA{//F(y)}}`{\BBB{//GF(y)}};{\check{G}}]
  \morphism(0,-900)|l|/{@{->}@/^20pt/}/<0,900>[{\BBB{//GF(y)}}`{\AAA{//F(y)}};{\varepsilon_{F(y)}\l!\circ\check{F}}]
  \morphism(0,-900)|r|/{@{->}@/^20pt/}/<0,-900>[{\BBB{//GF(y)}}`{\BBB{/y}};{\eta_y^\Leftarrow}]
  \morphism(0,-1800)|l|/{@{->}@/^20pt/}/<0,900>[{\BBB{/y}}`{\BBB{//GF(y)}};{\eta_y\!}]
  \morphism(-105,-450)|a|/{@{-|}@<-7pt>}/<210,0>[{\phantom{O}}`{\phantom{O}};{\left(\varepsilon{,}\eta\right)}]
  \morphism(-105,-1350)|a|/{@{-|}@<-7pt>}/<210,0>[{\phantom{O}}`{\phantom{O}};{\left(\delta{,}\rho\right)}]
}
\defdiag{theoremthatusespreservationofcommasecond}{   
  \morphism(0,0)|r|/{@{->}@/^20pt/}/<0,-750>[{\AAA{//F(y)}}`{\BBB{//GF(y)}};{\check{G}}]
  \morphism(0,-750)|l|/{@{->}@/^20pt/}/<0,750>[{\BBB{//GF(y)}}`{\AAA{//F(y)}};{\varepsilon_{F(y)}\l!\circ\check{F}}]
  \morphism(0,-750)|r|/{@{->}@/^20pt/}/<0,-750>[{\BBB{//GF(y)}}`{\BBB{/GF(y)}};{\id_{GF(y)}^\Leftarrow}]
  \morphism(0,-1500)|l|/{@{->}@/^20pt/}/<0,750>[{\BBB{/GF(y)}}`{\BBB{//GF(y)}};{\id_{GF(y)}\!}]
  \morphism(0,-1500)|r|/{@{->}@/^20pt/}/<0,-750>[{\BBB{/GF(y)}}`{\BBB{/y}};{\eta_y^\ast}]
  \morphism(0,-2250)|l|/{@{->}@/^20pt/}/<0,750>[{\BBB{/y}}`{\BBB{/GF(y)}};{\eta_y!}]
  \morphism(-105,-375)|a|/{@{-|}@<-7pt>}/<210,0>[{\phantom{O}}`{\phantom{O}};{\left(\varepsilon{,}\eta\right)}]
  \morphism(-105,-1125)/-|/<210,0>[{\phantom{O}}`{\phantom{O}};]
  \morphism(-105,-1875)/-|/<210,0>[{\phantom{O}}`{\phantom{O}};]
}
\defdiag{theoremthatusespreservationofcommathird}{   
  \morphism(0,0)|r|/{@{->}@/^20pt/}/<0,-750>[{\AAA{//F(y)}}`{\AAA{/F(y)}};{\id_{F(y)}^\Leftarrow}]
  \morphism(0,-750)|l|/{@{->}@/^20pt/}/<0,750>[{\AAA{/F(y)}}`{\AAA{//F(y)}};{\id_{F(y)}\!}]
  \morphism(0,-750)|r|/{@{->}@/^20pt/}/<0,-750>[{\AAA{/F(y)}}`{\BBB{/GF(y)}};{\check{G}}]
  \morphism(0,-1500)|l|/{@{->}@/^20pt/}/<0,750>[{\BBB{/GF(y)}}`{\AAA{/F(y)}};{\varepsilon_{F(y)}!\circ\check{F}}]
  \morphism(0,-1500)|r|/{@{->}@/^20pt/}/<0,-750>[{\BBB{/GF(y)}}`{\BBB{/y}};{\eta_y^\ast}]
  \morphism(0,-2250)|l|/{@{->}@/^20pt/}/<0,750>[{\BBB{/y}}`{\BBB{/GF(y)}};{\eta_y!}]
  \morphism(-105,-375)/-|/<210,0>[{\phantom{O}}`{\phantom{O}};]
  \morphism(-105,-1125)|a|/{@{-|}@<-7pt>}/<210,0>[{\phantom{O}}`{\phantom{O}};{\left(\varepsilon{,}\eta\right)}]
  \morphism(-105,-1875)/-|/<210,0>[{\phantom{O}}`{\phantom{O}};]
}
\defdiag{definingtheimageofthelefttwoadjointinfamsecond}{   
  \morphism(0,0)/{@{->}@/^30pt/}/<975,0>[{\displaystyle\coprod_{j=1}^{n}x_j}`{\displaystyle\coprod_{j=1}^{m}y_j};]
  \morphism(0,0)/{@{->}@/_30pt/}/<975,0>[{\displaystyle\coprod_{j=1}^{n}x_j}`{\displaystyle\coprod_{j=1}^{m}y_j};]
  \morphism(488,150)/=>/<0,-300>[{\phantom{O}}`{\phantom{O}};]
}
\defdiag{definingtheimageofthelefttwoadjointinfam}{   
  \morphism(0,0)|a|/{@{->}@/^30pt/}/<750,0>[{x_i}`{y_{t_0(i)}};{t_i}]
  \morphism(0,0)|b|/{@{->}@/_30pt/}/<750,0>[{x_i}`{y_{t_0(i)}};{t_i'}]
  \morphism(750,0)/->/<675,0>[{y_{t_0(i)}}`{\displaystyle\coprod_{j=1}^{m}y_j};]
  \morphism(375,150)/=>/<0,-300>[{\phantom{O}}`{\phantom{O}};]
}
\defdiag{basicadjunctionoffam}{   
  \morphism(0,0)|b|/{@{->}@/_25pt/}/<825,0>[{\AAA}`{\Famf\left({\AAA}\right)};{I}]
  \morphism(825,0)/{@{->}@/_25pt/}/<-825,0>[{\Famf\left({\AAA}\right)}`{\AAA};]
  \morphism(412,105)|r|/{@{-|}@<-7pt>}/<0,-210>[{\phantom{O}}`{\phantom{O}};{\left(\varepsilon{,}\eta\right)}]
}
\defdiag{basiccommutativitydiagramoftheextensivity}{   
  \morphism(0,-525)|b|/->/<1125,0>[{\displaystyle\AAA{/}{\coprod_{j=1}^{n}y_j}}`{\Famf\left({\AAA}\right)/\left({y}_j\right)_{j\in\left\{{1},\ldots{,}n\right\}}};{\eta_{Y}^\ast\circ\,{\check{I}_\AAA}}]
  \morphism(0,-525)|l|/<-/<0,525>[{\displaystyle\AAA{/}{\coprod_{j=1}^{n}y_j}}`{\displaystyle\prod_{j=1}^n\AAA/y_j};{\simeq}]
  \morphism(1125,-525)|r|/<-/<0,525>[{\Famf\left({\AAA}\right)/\left({y}_j\right)_{j\in\left\{{1},\ldots{,}n\right\}}}`{\displaystyle\prod_{j=1}^n\Famf\left(\AAA/y_j\right)};{\simeq}]
  \morphism(0,0)|a|/->/<1125,0>[{\displaystyle\prod_{j=1}^n\AAA/y_j}`{\displaystyle\prod_{j=1}^n\Famf\left(\AAA/y_j\right)};{\prod_{j=1}^n{I_{\AAA/y_j}}}]
  \morphism(382,-262)|m|//<210,0>[{\phantom{O}}`{\phantom{O}};{\cong}]
}

\def\pu{}
\fi	


\title{Lax comma $2$-categories and admissible $2$-functors}
\author{Maria Manuel Clementino and Fernando Lucatelli Nunes}
\dedication{In memory of Marta Bunge}
\address{(1,2): University of Coimbra, CMUC, Department of Mathematics, 3000-143 Coimbra, Portugal.\\ 
(2): Departement Informatica, Universiteit Utrecht, Nederland}
\eaddress{mmc@mat.uc.pt and f.lucatellinunes@uu.nl}
\amsclass{18N10, 18N15, 18A05, 18A22, 18A40}


\thanks{This work was supported through the programme ``Oberwolfach Leibniz Fellows'' by the Mathematisches Forschungsinstitut Oberwolfach in 2022. 
	This research was partially supported  by the CMUC, Centre for Mathematics of the University of Coimbra - UIDB/00324/2020, funded by the Portuguese Government through FCT/MCTES, and by the Institut de Recherche en Math\'{e}matique et Physique (IRMP, UCLouvain, Belgium).}

\maketitle

\begin{abstract}
This paper is a contribution towards a two dimensional extension of the basic
ideas and results of Janelidze's Galois theory. In the present paper, we give a suitable counterpart notion
to that of \textit{absolute admissible Galois structure} for the lax idempotent context, compatible with
the context of \textit{lax orthogonal factorization systems}. As part of this work, we study lax comma $2$-categories,
giving analogue results to the basic properties of the usual comma categories. We  show that each morphism of a $2$-category induces a $2$-adjunction  between lax comma $2$-categories and comma $2$-categories, playing the role of the usual \textit{change-of-base functors}. With these induced $2$-adjunctions, we are able to show that each $2$-adjunction induces $2$-adjunctions between lax comma $2$-categories and comma $2$-categories, which are our analogues of the usual lifting to the comma categories used in Janelidze's Galois theory. We give sufficient conditions under which these liftings are $2$-premonadic and induce a lax idempotent $2$-monad, which corresponds to our notion of $2$-admissible $2$-functor. In order to carry out this work, we analyse when a composition of
$2$-adjunctions is a lax idempotent $2$-monad, and when it is $2$-premonadic. We give then examples of our $2$-admissible $2$-functors
(and, in particular, simple $2$-functors), especially using a result that says that all admissible ($2$-)functors in the
classical sense are also $2$-admissible (and hence simple as well).
\end{abstract}

\tableofcontents
\setcounter{secnumdepth}{-1}

\section{Introduction}
Categorical Galois theory, originally developed by Janelidze~\cite{MR1061480, MR1822890}, gives a unifying setting for most of the formerly introduced
Galois type theorems, even generalizing most of them.
It neatly gives a common ground for Magid's Galois theory
of commutative rings, Grothendieck's theory of \'{e}tale covering of schemes, and central extension of groups.
Furthermore, since its genesis, Janelidze's Galois theory has found several developments, applications and examples in new settings (see, for instance, \cite{MR1397399}, \cite{MR3275274}, \cite{MR3207214}, \cite[Theorem~4.2]{MR1245796}, and \cite[Theorem~9.8]{2016arXiv160604999L}).

The most elementary observation on factorization systems and
 Janelidze's Galois theory is that, in the suitable setting
of finitely complete categories, the notion of
absolute admissible Galois structure coincides with
that of a semi-left-exact reflective functor/adjunction (see, for instance, \cite[Section~5.5]{MR1822890} or \cite{zbMATH01024334}).

Motivated by the fact above and the theory of \textit{lax orthogonal factorization
	systems}~\cite{zbMATH07249998, MR3708821, MR3545937}, we have started a project whose aim is to
 investigate
a two dimensional extension of the basic
ideas and results of (absolute) Janelidze's Galois theory. We deal herein with a
key step of this endeavor, that is to say, we develop the basics in order to give a suitable counterpart notion
to that of \textit{absolute admissible Galois structure}.

We adopt the \textit{usual} viewpoint that the
$2$-dimensional analogue of an idempotent monad (full reflective functor) is that of a lax idempotent
monad (pre-Kock-Z\"{o}berlein $2$-functor).
Therefore the concept of an admissible Galois structure
within our context should be a lax idempotent
counterpart to the notion of \textit{semi-left exact reflective functor}; namely, an appropriate
notion of semi-left exact functor for the
context of \cite{MR3545937}.

We  study the lifting of $2$-adjunctions to comma type $2$-categories. We find two possible liftings which deserve interest. The underlying adjunction of the first type of lifting is the usual $1$-dimensional case, while the other one,
more relevant to our context, is
a counterpart to the lifting of the $2$-monad given in
\cite{MR3545937} by comma objects. The last one requires us to study the lax analogue notion for comma categories, the notion of \textit{lax comma $2$-categories} of the title.

We study the basic aspects of lax comma $2$-categories. Among them,
the $2$-adjunction between the usual comma $2$-category and the lax comma $2$-category (for each object), and a counterpart for the
usual change-of-base $2$-functors, which comes into play as a fundamental aspect of our work and, specially, to introduce
the definition of \textit{$2$-admissible $2$-adjunction}.

With these analogues of the change-of-base $2$-functors, we are able to introduce the lifting of
each $2$-adjunction to a $2$-adjunction between the lax comma $2$-category and the comma $2$-category as a composition of $2$-adjunctions. Namely, the composition of a straightforward lifting to the lax comma $2$-categories with a change-of-base $2$-functor induced by the appropriate component of the unit.
Fully relying on the study of properties of compositions of 2-adjunctions,
we  investigate the properties of these liftings of the $2$-adjunctions.
Namely, we  show under
which conditions these liftings induce lax idempotent
$2$-monads (the simple $2$-adjunctions of \cite{MR3545937}), recovering
 one characterization given in \cite{MR3545937} of their \textit{simple
	$2$-adjunctions}. We give also a characterization
of the $2$-functors whose introduced lifting is lax idempotent and $2$-premonadic,
the
\textit{$2$-admissible $2$-functors} within our context.

In Section \ref{section1 Preliminaries} we recall basic aspects and terminology of $2$-categories, such as $2$-adjunctions and  $2$-monads, finishing the section giving aspects on \textit{raris}, right-adjoint right-inverses (see Definition \ref{ralirarilalilaridefinition}) within a $2$-category. Taking the opportunity to fix notation, we also recall the universal properties of the
main two dimensional limits used in our work in Section \ref{twodimensionallimits section2}, that is to say, the definitions
of conical $2$-limits and comma objects.

In Section \ref{sectionkockzoberlein} we recall and show aspects on idempotent and lax idempotent $2$-monads needed to our work
on admissible and $2$-admissible $2$-functors, also introducing a characterization of
the $2$-adjunctions that induce lax idempotent $2$-monads, called herein lax idempotent $2$-adjunctions (see, for instance, Theorem \ref{characterizationlaxidempotent2adjunction}).

In Section \ref{compositionof2adjunctionssection}  we introduce the main concepts and results on composition of $2$-adjunctions
in order to introduce the notions of simple, admissible and $2$-admissible $2$-adjunctions (see, for instance,
Definitions \ref{Admissibletwoadjunction}, \ref{SIMPLEtwoadjunction}, and \ref{twoadmissibledefinition}).
The results focus on characterizing and giving conditions under which the composition of $2$-adjunctions is an idempotent/lax idempotent (full reflective/pre-Kock-Z\"{o}berlein) $2$-adjunction ($2$-functor). Most of them are analogues for the simpler case of
idempotent $2$-adjunctions (see, for instance, Theorem \ref{admissibilitywrtJH} which characterizes when the composition
of right $2$-adjoints is pre-Kock-Z\"{o}berlein).

In Section \ref{basicdefinitionsarticle} we recall the notion of lax comma $2$-categories $\AAA //y $, for each $2$-category $\AAA $ and object $y\in\AAA $ (see Definition \ref{definitionoflaxcommacategory}). This notion
has already appeared in the literature (see, for instance, \cite[I,5]{MR0371990}, \cite[\S~6]{MR0249483}, 
\cite[Exercise~5, p.~115]{MR1712872} and \cite[p.~305]{MR558494}). We, then, introduce the change-of-base $2$-functors for lax comma $2$-categories.
More precisely, we show that, for each morphism $c: y\to z $ in a $2$-category $\AAA $ with comma objects, we have an induced
$2$-adjunction
%
\pu
\begin{equation*}\label{changeofthebaselaxcommacategoriesequationchangeofbase}
\diag{changeofthebaselaxcommacategorieschangeofbase}.
\end{equation*}
between the lax comma $2$-category $\AAA //z $ and the comma $2$-category $\AAA /y$. We give an explicit construction of this $2$-adjunction: see Theorem \ref{teoremamudancadebaselassa}.

Provided that $\AAA $ has pullbacks and comma objects, these induced $2$-adjunctions, together with the classical change-of-base $2$-functors,
give the $2$-adjunctions
%
\pu
\begin{equation*}
\diag{compositionoftwoadjunctionschangeofbasecommacompositionofpullbackcomma}
\end{equation*}		
in which the composition of
	$c!\dashv c^\ast : \AAA /z \to \AAA /y $ with
	$\id _z\!\dashv id _z^\Leftarrow : \AAA //z \to \AAA /z $
	is, up to $2$-natural isomorphism, the $2$-adjunction $c\!\dashv c^\Leftarrow : \AAA //z  \to \AAA /y $ (see Theorem \ref{relationofchangeofbasecomma}).
We finish Section \ref{change-of-base functor} showing that, whenever it is well defined,  $\id_ y^\Leftarrow $ is pre-Kock-Z\"{o}berlein (Theorem \ref{laxidempotentcoherencelaxcommacomma}).

The main point of Section \ref{secao de admissibilidade} is to introduce our notions of admissibility and $2$-admissibility
(Definition \ref{maindefinition}), relying on the definitions previously introduced in Section \ref{compositionof2adjunctionssection}. We also use the main results of Section \ref{compositionof2adjunctionssection}
to characterize and give conditions under which a $2$-functor is $2$-admissible (see, for instance, Corollaries \ref{oneofthemaincorollaries} and \ref{preKockZoberleinFdashvG}).

We finish Section \ref{secao de admissibilidade} with a fundamental observation on admissibility and $2$-admissibility, namely,  Theorem \ref{southafricantheorem}.  It says that, provided that
$\AAA $ has comma objects,
if $F\dashv G $ is admissible in the classical
sense (called herein \textit{admissible w.r.t. the basic fibration}), meaning that $G$ itself is full reflective and the compositions
$$ \eta_y^\ast\,\circ \check{G}: \AAA  /F(y) \to \BBB /y $$
are full reflective for all $y$, then $G$ is $2$-admissible, which means that the compositions
$$ \eta_y^\Leftarrow\,\circ \check{G}: \AAA  //F(y) \to \BBB /y $$
are pre-Kock-Z\"{o}berlein for all objects $y$.
We discuss examples of $2$-admissible $2$-functors (and hence also simple $2$-functors) in Section
\ref{sectionexamplesexamples}. Most examples
are about cocompletion of $2$-categories, making use of Theorem  \ref{southafricantheorem}.

\setcounter{secnumdepth}{5}

\section{Preliminaries}\label{section1 Preliminaries}
Let $\Cat $ be the cartesian closed category of categories  in some universe. We denote the internal hom  by
$$\Cat (-,-): \Cat ^\op\times \Cat\to \Cat .$$

A $2$-category $\AAA $ herein is the same as a $\Cat $-enriched category.
We denote the enriched hom
of a $2$-category $\AAA $ by
$$\AAA (-,-): \AAA ^\op \times \AAA \to \Cat $$
which, again, is of course a $2$-functor.
As usual,
the composition of $1$-cells (morphisms) are denoted
by $\circ $, $\cdot $, or omitted whenever it is clear from the context.
The vertical composition of $2$-cells is denoted by $\cdot $  or omitted when it is clear,  while the horizontal composition is denoted
by $\ast$. Recall that, from the vertical and horizontal compositions, we construct the fundamental operation of \textit{pasting}~\cite{MR0357542, MR1040947}.

Finally, if $f: w\to x $, $g:y\to z $ are $1$-cells of $\AAA $, given a $2$-cell $\xi : h\Rightarrow h' : x\to y $,
motivated by the case of $\AAA = \Cat $, we use interchangeably the notations
%
\pu
%
\pu
\begin{equation}
\id_g\ast \xi\ast  \id_ f\quad =\quad \diag{whiskering1} \quad =\quad  g\xi f
\end{equation}
to denote the whiskering of $\xi $ with $f$ and $g$.

Henceforth, we consider the $3$-category of $2$-categories, $2$-functors, $2$-natural transformations and modifications, denoted by $2$-$\Cat $. We refer to \cite{MR0357542, MR0299653} for the basics on $2$-dimensional category theory, and, more particularly, to
the definitions of adjunctions, monads and Kan extensions inside a $2$-category.
Moreover, we also extensively assume  aspects of $2$-monad theory. The pioneering reference is \cite{MR1007911},
while we mostly follow the terminology (and results) of \cite{2016arXiv160703087L}.

In this paper, we consider the \textit{strict} versions of $2$-dimensional adjunctions and monads: the concepts
coincide with the $\Cat $-enriched ones. A \textit{$2$-adjunction}, denoted by $$(F\dashv G, \varepsilon, \eta): \AAA\to\BBB ,$$
consists of $2$-functors
%
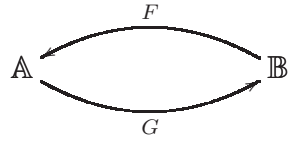
\pu
$$ \diag{twofoldtwofunctors} $$
with $2$-natural transformations $\varepsilon : FG\Longrightarrow \id _ \AAA $
and $\eta : \id _ \BBB \Longrightarrow GF $ playing the role of the \textit{counit} and the \textit{unit} respectively. More precisely, the equations of $2$-natural transformations
%

\pu

%

\pu

\begin{equation*}\tag{triangle identities}
\diag{triangleidentityadjunctiondiagram1}\, =\, \id _ G\quad\mbox{ and }\quad \diag{triangleidentityadjunctiondiagram2}\, =\, \id _F
\end{equation*}
hold.
We usually denote a $2$-adjunction 	$(F\dashv G, \varepsilon, \eta): \AAA\to\BBB $
by
%
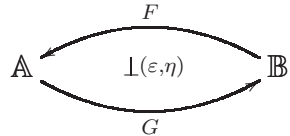
\pu
$$ \diag{basicadjunction} $$
or by $F\dashv G: \AAA\to \BBB $ for short, when the counit and unit are already given.

A \textit{ $2$-monad} on a $2$-category $\BBB $ is a triple $\mathcal{T} = (T, \mu , \eta) $ in which
$T: \BBB \to \BBB $ is an endo-$2$-functor and $\mu, \eta $ are $2$-natural transformations playing the role of the multiplication and the unit respectively.
That is to say, $\mu $ and $\eta $ are $2$-natural transformations such that the equations
%

\pu
%

\pu

%

\pu

%

\pu

\begin{equation*}\tag{associativity of a $2$-monad}
\diag{leftsideoftheequationassociativityofmonad} \quad =\quad \diag{rightsideoftheequationassociativityofmonad}
\end{equation*}
\begin{equation*}\tag{identity of a $2$-monad}
\diag{firstsideoftheequationidenityofamonad}\quad = \quad \diag{secondsideoftheequationidenityofamonad}\quad = \quad \id _T
\end{equation*}
hold.

Since the notions above coincide with the $\Cat $-enriched ones, it should be noted that the formal theory of monads applies to this case.
More precisely, every $2$-adjunction does induce a $2$-monad, and we have the usual Eilenberg-Moore and Kleisli factorizations
of a right $2$-adjoint functor (\textit{e.g}  \cite[Section~2]{MR0299653} or \cite[Section~3]{2019arXiv190201225L}),
which give rise respectively to the notions of $2$-monadic and Kleisli $2$-functors. Furthermore, we also have (the enriched version of) Beck's monadicity theorem~\cite[Theorem~II.2.1]{MR0280560}.

In this direction,
we use expressions like \textit{equivalence (or $2$-equivalence)},  and \textit{fully faithful $2$-functor} to mean
the (strict) $\Cat $-enriched notions: that is to say, respectively, \textit{equivalence} in the $2$-category of $2$-categories, and a $2$-functor that is \textit{locally an isomorphism}.

\subsection{Lalis and ralis}
To refer to adjunctions where the unit or counit is an identity, we adopt a terminology similar to the one introduced by Gray in \cite[0.3.B]{zbMATH03305157}. More precisely:
\begin{defi}\label{ralirarilalilaridefinition}
	Assume that $(f\dashv g, v, n )$ is an adjunction in a $2$-category $\AAA $.
	\begin{itemize}
		\renewcommand\labelitemi{--}
		\item If the counit $ v $ is the identity $2$-cell,  $(f\dashv g, v, n )$ is called a \textit{rari adjunction (or rari pair)}, or a \textit{lali adjunction}.
		
		If there is a rari adjunction $f\dashv g $,
		the morphism $f$ is called a \textit{lali (left-adjoint left-inverse)}, while the morphism
		$g$ is called a \textit{rari (right-adjoint right-inverse)}.

		\item If the unit  $ n $ is the identity $2$-cell,  $(f\dashv g, v, n )$ is called a \textit{rali adjunction}, or a \textit{lari adjunction}.
		
		If there is a rali adjunction $f\dashv g $,
		the morphism $f$ is called a \textit{lari}, while the morphism
		$g$ is called a \textit{rali}.
	\end{itemize}	
\end{defi}

Laris (ralis) are closed by composition, and have specific cancellation properties. We recall them below.

\begin{lem}\label{preparacaoparacancelamentolaris}
	Assume that
%
\pu
\begin{equation}
\diag{compositionofadjunctionsrarilali}
\end{equation}
are
adjunctions in $\AAA $.
\begin{enumerate}[a)]
	\item Assuming that $f\dashv g $ is a lari adjunction: we have that
	$ f  f'\dashv g'  g $ is a lari adjunction if, and only if, $f'\dashv g' $
	is a lari adjunction as well.\label{dualaralilaricancellation}
	\item Assuming that $f'\dashv g' $ is a lali adjunction: the adjunction
	$f f'\dashv g' g $ is a lali adjunction if, and only if, $f\dashv g $ is a lali adjunction as well. \label{aralilaricancellation}
\end{enumerate}
\end{lem}
\begin{proof}
Assuming that $ n  $ is an isomorphism, we have that the unit
%
\pu
\begin{equation}
\diag{proofodthelalicancellation}
\end{equation}
of the composition  $ff'\dashv g'g$ is invertible if, and only if, $ n' $ is invertible. This proves \ref{aralilaricancellation} and, dually, we get \ref{dualaralilaricancellation}.   	
\end{proof}

Of course, the situation is simpler when we consider isomorphisms. That is to say:

\begin{coro}\label{precisecancellationisomorphism}
Assume that
%
\pu
\begin{equation}
\diag{compositionofadjunctionsrarilalicorollaryforisomorphisms}
\end{equation}
are morphisms in $\AAA $ such that $ (f') ^{-1} = g' $ and $ (f'') ^{-1} = g'' $.
There is a lali (rali) adjunction $f'\cdot f\cdot f''\dashv g''\cdot g \cdot g' $
if and only if there is a lali (rali) adjunction $f\dashv g $.
\end{coro}
\begin{proof}
If $f\dashv g $ is a lali (rali) adjunction, since $f'\dashv g' $ and $f''\dashv g'' $ are of course lali and rali adjunctions, it follows that the composite
$$ f'\cdot f\cdot f''\dashv g''\cdot g \cdot g' $$ is a lali (rali) adjunction by Lemma \ref{preparacaoparacancelamentolaris}.

Conversely, if 	$f'\cdot f\cdot f''\dashv g''\cdot g \cdot g' $ is a lali (rali) adjunction,
since $g'\dashv f' $ and $g''\dashv f'' $ are lali and rali adjunctions, we get that the composite
$$  g'\cdot f'\cdot f\cdot f''\cdot g''\dashv f''\cdot g''\cdot g \cdot g'\cdot  f', $$
which is $f\dashv g $, is a lali adjunction.
\end{proof}

But we also have a stronger cancellation property:

\begin{theo}[Left cancellation property]\label{cancellationpropertyralis}
	Let $f: x\to w, f': y\to x $ be morphisms of a $2$-category $\AAA $.	
	\begin{enumerate}[a)]
		\item Assuming that $f: x\to w $ is a lari:
		the composite $ f  f': y\to w $
		is a lari if, and only if,
		$f': y\to x $  is a lari as well.\label{cancellationpropertyralisA}
		\item Assuming that $f $ is a rari: the composite
		$ff'$ is a rari if and only if $f' $ is a rari. \label{cancellationpropertyralisB}
	\end{enumerate}
\end{theo}
\begin{proof}
	By Lemma \ref{preparacaoparacancelamentolaris}, if $f $ and $f' $ are laris,  the composite $ff' $ is a lari as well.	
	
Conversely, assume that $f $ and $ff'$ are laris. This means that there are adjunctions
%
\pu
%
\pu

\begin{equation*}
\diag{existingadjunction1}\quad \diag{existingadjunction2}
\end{equation*}
in $\AAA $ such that $n=\id _{gf} $ and $\hat{n} = \id _{\hat{g}ff'}$.

We claim that
%
\pu
\begin{equation}
\left(f'\dashv \hat{g}f,\,  \diag{counitofthecancellationtheoremlaris} , \id_{\hat{g}ff'}\right)
\end{equation}
is a (lari) adjunction. In fact, the triangle identities follow from the facts that
the equations
$ \hat{v}ff' = \id_{ff'} $ and $\hat{g}\hat{v} = \id_{\hat{g}} $  hold.

Finally, the statement \ref{cancellationpropertyralisB} is the codual of \ref{cancellationpropertyralisA}.
\end{proof}

On the one hand, the \textit{left cancellation property} of Theorem \ref{cancellationpropertyralis} does not hold for lalis or ralis. For instance, in $\Cat $, we consider the terminal category $\mathsf{1} $ and the
category $\mathsf{2} $ with two objects and only one nontrivial morphism between them. The morphisms
%
\pu
\begin{equation}
\diag{exampleofnocancellationpropertyoflalis}
\end{equation}
are lalis. But the inclusion $d^0: \mathsf{1}\to \mathsf{2} $  of the terminal object of $\mathsf{2} $ is not a lali, since it does not have a right adjoint.
On the other hand, the dual of Theorem \ref{cancellationpropertyralis} gives a right cancellation property for ralis and lalis.

\begin{coro}[Right cancellation property]\label{leftcancellationpropertyralis}
	Let $f: x\to w, f': y\to x $ be morphisms of a $2$-category $\AAA $.	
	If $f': y\to x $ is a lali (rali): we have that
	$f: x\to w $  is a lali (rali) if, and only if,  the composite $ f  f': y\to w $
	is a lali (rali) as well.\label{leftcancellationpropertyralisA}
\end{coro}

\section{Two dimensional limits}\label{twodimensionallimits section2}
In this section, we recall basic universal constructions related to the results of this paper.
Two dimensional limits are the same as weighted limits in the $\Cat $-enriched context~\cite{MR0280560}. We refer, for instance, to  \cite{MR0401868} for the basics on $2$-dimensional limits.
We are particularly interested
in \textit{conical $2$-(co)limits} and \textit{comma objects}.

\subsection{Conical $2$-limits}
Two dimensional conical (co)limits are just weighted (co)limits with a weight constantly equal to the terminal category $\mathsf{1} $.
Henceforth, the words \textit{(co)product}, \textit{pullback/pushout} and \textit{(co)equalizer}
refer to the $2$-dimensional versions of each of those (co)limits. For instance, if $a :x\to y $, $b : w\to y $ are morphisms of a $2$-category $\AAA $, assuming its existence, the \textit{pullback} of
$b $ along $a $
is an object $\displaystyle x\times _{(a,b)} w $ together
with $1$-cells $ a^\ast (b):x\times _{(a,b)} w \to x  $ and $b^\ast (a)  :  x\times _{(a,b)} w\to w  $
making the diagram
%
\pu
\begin{equation}\label{definitionpullbackdiagram}
\diag{pullbackdiagramdefinition}
\end{equation}
commutative, and satisfying the following universal property. For every object $z $ and every pair of $2$-cells $$(\xi _ 0 : h_0\Rightarrow h_0' : z\to x ,\, \xi _ 1 : h_1\Rightarrow h_1' : z\to w ) $$
such that the equation
%
\pu
%
\pu
\begin{equation}
\diag{twocellpullbackdefinitionrightside} \quad =\quad \diag{twocellofpullbackdefinitionleftside}
\end{equation}
holds,
there is a unique $2$-cell $\xi : h\Rightarrow h' :  z\to x\times_{(a,b)}w $ satisfying the equations
\begin{center}
	$\id _ {a^\ast (b)}\ast \xi  = \xi _ 0 $ and $\id _ {b^\ast (a)  }\ast\xi = \xi _ 1 $.	
\end{center}

\begin{rem}
	It is clear that the concept of \textit{pullback} in locally discrete $2$-categories
	coincides with the concept of ($1$-dimensional) \textit{pullback} in the underlying
	categories.
	
	Moreover, when a \textit{pullback} exists in a $2$-category,
	it is isomorphic to the ($1$-dimensional) \textit{pullback} in the underlying category.
	
	Finally, both  the statements above are also true if
	\textit{pullback} is replaced by any type of conical $2$-limit
	with a locally discrete \textit{shape} (domain).
	
\end{rem}

\subsection{Comma objects}\label{definitionofcommaobjects}
If $a :x\to y $, $b : w\to y $ are morphisms of a $2$-category $\AAA $,
the comma object of $a$ along $b$, if it exists, is an object $a\downarrow b $  with the
following
universal property. There are $1$-cells
 $a^{\Rightarrow } (b) : a\downarrow b\to x  $ and $  b^{\Leftarrow } (a) :  a\downarrow b\to w $
and a $2$-cell
%
\pu
\begin{equation}\label{diagramdefinitioncommaobject}
\diag{commadiagramdefinition}
\end{equation}
such that:
\begin{enumerate}
	\item For every triple $(h_0: z\to x, h_1: z\to w,  \livb : a h_0 \Rightarrow b h_1  ) $ in which $h_0, h_1 $ are morphisms
	and $\livb $ is a $2$-cell of $\AAA $, there is a unique morphism $h:   z\to a\downarrow b $ such that the equations
	$h_0 =   a^\Rightarrow (b)  \cdot h $, $h_1 =  b^\Leftarrow (a)\cdot h   $ and
%
\pu
\begin{equation}
\diag{morphismcommadiagramdefinition}
\end{equation}
hold.
\item For every pair of $2$-cells $(\xi _ 0 : h_0\Rightarrow h_0' :z\to x,\, \xi _ 1 : h_1\Rightarrow h_1' : z\to w ) $ such that
%
\pu
%
\pu
\begin{equation}
\diag{commatwocellpullbackdefinitionrightside} \quad =\quad \diag{commatwocellofpullbackdefinitionleftside}
\end{equation}
holds,
	there is a unique $2$-cell $\xi : h\Rightarrow h' : z\to a\downarrow  b $ such that
	$\id _ {a^\Rightarrow (b)}\ast \xi   = \xi _ 0 $ and $\id _ {b^\Leftarrow (a)}\ast \xi = \xi _ 1 $.
\end{enumerate}

\begin{rem}
If $\AAA $ is a locally discrete $2$-category, the comma object of
a morphism
$a$ along $b$ has the same universal property of the pullback
of $a$ along $b$.

\end{rem}

\section{Lax idempotent $2$-adjunctions}\label{sectionkockzoberlein}

Herein, our standpoint is that the notion of \textit{pre-Kock-Z\"{o}berlein $2$-functor} is the $2$-dimensional counterpart
of the notion of \textit{full reflective functor}. In this section, we recall the basic definitions and give basic characterizations, but we refer to \cite{MR1359690, MR1432190} and \cite[Ch.~4]{zbMATH05036792} for
fundamental aspects and examples of lax idempotent $2$-monads.

\begin{defi}[Lax idempotent $2$-monad]
A \textit{lax
idempotent $2$-monad} is a $2$-monad $\mathcal{T} = (T, \mu , \eta) $
such that we have a rari adjunction $\mu \dashv \eta \ast \id_T   $.

An \textit{idempotent $2$-monad} is a $2$-monad $\mathcal{T} = (T, \mu , \eta) $  such that
$\mu $ is invertible or, in other words, it is a lax idempotent $2$-monad   such that  $\mu \dashv \eta \ast \id_T $ is a rali adjunction as well.
\end{defi}

More explicitly, a $2$-monad $\mathcal{T} = (T, \mu , \eta) $
on a $2$-category $\BBB $ is lax idempotent
if there is a modification
%
\pu
$$\diag{modificationofthedefinitionoflaxidempotenttwomonad} $$
such that, for each object $z\in \BBB $,
%
\pu
%
\pu
\begin{equation*}
\diag{modificationofthedefinitionoflaxidempotenttwomonadtriangleidentityotwo}\qquad\qquad \diag{modificationofthedefinitionoflaxidempotenttwomonadtriangleidentityoone}
\end{equation*}
are respectively the identity $2$-cells on $\eta_{T(z)} $ and on $\mu _z $.

\begin{rem}[Dualities and self-duality]\label{dualitylaxidempotentremark}
The concepts of lax idempotent and idempotent $2$-monads are actually notions that can
be defined inside any $3$-category
(or, more generally, tricategory~\cite{MR1261589}). Therefore they have eight dual notions each (counting the concept itself).

However, the notions of lax idempotent and idempotent $2$-monads are self-dual, that is to say,
the dual notion coincides with itself. More precisely, a triple $\mathcal{T} = (T, \mu , \eta) $ is a (lax) idempotent $2$-monad in the $3$-category $2\textrm{-}\Cat $ if and only if
the corresponding triple is also a (lax) idempotent $2$-monad in the $3$-category  $\left( 2\textrm{-}\Cat\right) ^\op $.

Furthermore, the notion of idempotent $2$-monad is self-$3$-dual, meaning that
the notion does not change when we invert the directions of the $3$-cells (which are, in our case, the modifications). However the $3$-dual of the notion
of lax idempotent $2$-monad is that of colax idempotent $2$-monad.

Finally, the notions obtained from the inversion of the directions of the $2$-cells, that is to say, the codual (or $2$-dual) concepts, are those of lax idempotent and idempotent $2$-comonads.
\end{rem}
\textit{Henceforth, throughout this section, we always assume that a $2$-adjunction
$$\diag{basicadjunction} $$
is given, and we denote by $\mathcal{T} = (T, \mu , \eta) $ the induced $2$-monad $(GF, G \varepsilon F, \eta ) $ on $\BBB $. }

\subsection{Idempotency}
There are several useful well-known characterizations of idempotent ($2$-)monads (see, for instance,
 \cite[p.~196]{MR1313497}).

\begin{lem}[Idempotent $2$-monad]\label{characterizationidempotentmonad}
The following statements are equivalent.
	\begin{enumerate}[i)]
		\item $\mathcal{T}$ is idempotent; \label{idempotentdefinitionbasictheorem}
		\item $T \eta $ (or $\eta T $) is an epimorphism;
		\item $\mu $ is a monomorphism; \label{trivialmonomorphismidempotent}
		\item $T\eta =\eta T  $; \label{idempotentequalunitsidempotent}
		\item $a: T(x)\to x $ is a $\mathcal{T}$-algebra structure if, and only if, $a\cdot \eta _ x = \id _x$;\label{novoestruturadealgebrapara2categoriadealgebrasestritas}
		\item $a: T(x)\to x $ is a $\mathcal{T}$-algebra structure if, and only if,
		$a $ is the inverse of $\eta _ x $;\label{2novoestruturadealgebrapara2categoriadealgebrasestritas}
		\item  the forgetful $2$-functor $ \mathcal{T}\textrm{-}\Alg _{\mathsf{s}}\to  \BBB $ between the $2$-category of strict $\mathcal{T}$-algebras and strict $\mathcal{T}$-morphisms (with modifications as $2$-cells) and the $2$-category $\BBB$  is fully faithful (that is to say, locally an isomorphism).\label{fullyfaithfulstrictalgebras}
	\end{enumerate}	
\end{lem}
\begin{proof}
Since
$\mu\cdot (\eta T) = \mu\cdot (T \eta ) = \id _ T $,
we have the following chain of equivalences:
$\mu $  is a monomorphism $\Leftrightarrow $ $ \mu $ is invertible $\Leftrightarrow $ $\eta T$ or $T\eta $ is invertible  $\Leftrightarrow $
$\eta T$ or $T\eta $ is an epimorphism. This proves the equivalence of the first three statements.

By the definition of monomorphism, \ref{trivialmonomorphismidempotent}  implies \ref{idempotentequalunitsidempotent}.
Conversely, assuming that $T\eta =\eta T  $, we have that $ T^2 \eta = T \eta T $ and, thus, we get that
%
\pu
%
\pu
$$ \left(T\eta\right)\cdot \mu \quad = \quad  \diag{equationforidempotentmonadinterchangelaw} \quad = \quad  \diag{equationforidempotentmonadinterchangelawtwo}\quad = \quad \id _{T^2}.  $$
Therefore $T \eta  $ is the inverse of $\mu $ and, hence, $\mu $ is a monomorphism.

Assuming one of the first four equivalent statements (and hence all of them), we have that,
 given a morphism $a: T(x)\to x $ such that
$a\cdot \eta _ x = \id _x $, the equation
\begin{equation}\label{rightinverseinversealgebrastructure}
\eta _ {x}\cdot  a = T(a) \cdot \eta _ {T(x)} = T(a\cdot \eta _x) = \id _ {T(x)}.
\end{equation}
holds. Thus, since $  \eta_ {T(x)}\cdot \eta_ x = T(\eta _x)\cdot \eta_ x  $ and $\mu = \left( T \eta \right)^{-1} $,
we conclude that
\begin{equation}\label{associativity trivially holds for algebras}
 a\cdot \mu _x =  \left(\eta_ {T(x)}\cdot \eta_ x\right)^{-1}
= \left(T\left(\eta_ {x}\right)\cdot \eta_ x\right)^{-1} = a\cdot T(a).
\end{equation}
This proves that \ref{novoestruturadealgebrapara2categoriadealgebrasestritas}
holds. Conversely, \ref{novoestruturadealgebrapara2categoriadealgebrasestritas} trivially implies \ref{trivialmonomorphismidempotent} (and, hence, all of the first four equivalent statements), since, for each $x\in\BBB $,  $\mu _x $ is a (free) $\mathcal{T}$-algebra structure for $x$.
Moreover, by Equations \eqref{rightinverseinversealgebrastructure} and \eqref{associativity trivially holds for algebras}, we conclude that the first four statements are also equivalent to
\ref{2novoestruturadealgebrapara2categoriadealgebrasestritas}.

Finally, recall that, for every $2$-monad $\mathcal{T}$ on a $2$-category $\BBB$, the forgetful functor $ \mathcal{T}\textrm{-}\Alg _{\mathsf{s}}\to  \BBB $ between the $2$-category of strict $\mathcal{T}$-algebras and strict $\mathcal{T}$-morphisms (with modifications as $2$-cells) and the $2$-category $\BBB$ 
is faithful. Assuming \ref{2novoestruturadealgebrapara2categoriadealgebrasestritas}, in order to verify that the forgetful functor is full, it is enough to see that, for any morphism
$f: x\to y $ of $\BBB$, if $a: T(x)\to x $, $b: T(y)\to y $ are  $\mathcal{T}$-algebra structures,
we have that the pasting
%
\pu
$$\diag{pastinginordertogetmorphismofalgebrasvertical} $$
is the identity $2$-cell and, hence, the morphism $f$ induces a morphism of algebras between $(x,a)$ and $(y,b)$.

Assuming
\ref{fullyfaithfulstrictalgebras},
we get that, for any object $x\in\BBB$,  $\eta_ {T(x)} $ induces
a morphism between the free $\mathcal{T}$-algebras $\left( T(x),\mu _x   \right)$ and $\left( T^2(x), \mu _{T(x)} \right) $.
That is to say,  $$\eta _ {T(x) }\cdot \mu _x = \mu _{T(x)}\cdot T(\eta _{T(x)} ) $$ and, since the right side of the equation above is equal to the identity on $T^2(x) $, we
conclude that $\mu _x $ is a split monomorphism. This proves that \ref{trivialmonomorphismidempotent} holds.
\end{proof}

A $2$-adjunction
induces an idempotent $2$-monad if, and only if,
 the induced $2$-comonad is also idempotent. More generally:

\begin{prop}\label{inducedidempotent}
The following statements are equivalent.
\begin{enumerate}[i)]
\item $\mathcal{T}$ is idempotent;\label{1oftheidempotentadjunction}
\item $F\eta $ (or $\eta G $) is an epimorphism; \label{equvalentdualidempotent2}
\item $\varepsilon F $ (or $G \varepsilon $) is a monomorphism;\label{equvalentdualidempotent}
\item The induced $2$-comonad is idempotent.\label{codualityequivalent4statement}
\end{enumerate}
\end{prop}
\begin{proof}
Since, by the triangle identities, we have that
\begin{center}
$\left( \varepsilon F\right) \cdot \left(F\eta\right) = \id _ F $ and
$\left( G\varepsilon \right) \cdot \left(\eta G\right) = \id _ G $,
\end{center}
we get that \ref{equvalentdualidempotent2} implies that $\varepsilon F $ or $G\varepsilon $ is invertible and, therefore, $G\varepsilon F = \mu $ is invertible.
Analogously, 	\ref{equvalentdualidempotent} implies \ref{1oftheidempotentadjunction}.

Moreover, if we assume that $\mathcal{T}$ is idempotent, by Lemma \ref{characterizationidempotentmonad}, we have that $$  GF \eta = \eta GF  $$ which, together with one of the triangle identities, implies that
%
\pu
%
\pu
\begin{equation*}
\left(F \eta\right)\cdot \left( \varepsilon F\right)\quad = \diag{equationforidempotentmonadinterchangelawtwoadjunction} =\quad
\diag{equationforidempotentmonadinterchangelawtwoadjunctiontwo} =\quad \id_{FGF}.
\end{equation*}
This proves that
\ref{1oftheidempotentadjunction} implies \ref{equvalentdualidempotent2} and
\ref{equvalentdualidempotent}.
Therefore we proved that \ref{1oftheidempotentadjunction}, \ref{equvalentdualidempotent2} and
\ref{equvalentdualidempotent} are equivalent statements.

Finally, since condition \ref{equvalentdualidempotent}
is codual and equivalent to condition \ref{equvalentdualidempotent2},
 we conclude that
\ref{1oftheidempotentadjunction} is equivalent to its codual -- that is to say,  to condition \ref{codualityequivalent4statement}.
\end{proof}

Motivated by the result above, we say that a $2$-adjunction
is \textit{idempotent} if it induces an idempotent $2$-(co)monad.

\begin{rem}\label{thinidempotency}
	If the $2$-adjunction $F\dashv G: \AAA\to\BBB $ is such that
	 the underlying category of $\AAA $ (or $\BBB $) is \textit{thin},
	then the induced $2$-monad is idempotent by Proposition \ref{inducedidempotent}. In particular, seeing categories
	as locally discrete $2$-categories and contravariant $2$-functors as covariant ones defined in the dual of the respective domains, any \textit{Galois connection}
	induces an idempotent ($2$-)(co)monad.
\end{rem}

If the $2$-adjunction $F\dashv G $ is idempotent and $G$ is $2$-monadic, $G$ is called a \textit{full reflective $2$-functor}. This terminology is justified by the well-known characterization below.

\begin{prop}[Full reflective $2$-functor]\label{fullreflective2functorcharacterization}
The following statements are equivalent.
\begin{enumerate}[i)]
		\renewcommand\labelitemi{--}
		\item $G$ is a full reflective $2$-functor;\label{1melhoradoidempotente}
		\item $F\dashv G $ is idempotent and $G$ is $2$-premonadic;\label{2melhoradoidempotente}
		\item $G$ is fully faithful;\label{3melhoradoidempotente}
		\item $\varepsilon $ is invertible.\label{4melhoradoidempotente}
\end{enumerate}
\end{prop}
\begin{proof}
Recall that a $2$-functor is $2$-premonadic if the (Eilenberg-Moore) comparison $2$-functor is fully faithful (that is to say, locally an isomorphism).

We have that \ref{1melhoradoidempotente} trivially implies \ref{2melhoradoidempotente}. Moreover, since the forgetful $2$-functor $ \mathcal{T}\textrm{-}\Alg _{\mathsf{s}} \to  \BBB $ is fully faithful whenever $ \mathcal{T}$ is idempotent,
we have that \ref{2melhoradoidempotente} implies \ref{3melhoradoidempotente}.

Since, for every pair of objects $w,x\in\AAA $,
the diagram
%
\pu
$$ \diag{equivalencesplitepifullyfaithful} $$
commutes, \ref{3melhoradoidempotente} and \ref{4melhoradoidempotente} are equivalent.

Assuming \ref{4melhoradoidempotente}, we have in particular that
$\varepsilon $ is a split epimorphism and $G$ reflects isomorphisms, hence, $G$ is $2$-monadic (see Proposition at  \cite[p.~236]{MR2056584}).
Furthermore, clearly, we also get that $G\varepsilon $ is a (split) monomorphism, which implies
that $F\dashv G $ is idempotent  by Proposition \ref{inducedidempotent}. Therefore \ref{4melhoradoidempotente} implies \ref{1melhoradoidempotente}.
\end{proof}

The dual notion of full reflective $2$-functor in $2$-$\Cat $ is called \textit{full co-reflective $2$-functor}. As a consequence of Proposition \ref{fullreflective2functorcharacterization}, we have:
\begin{coro}\label{coreflectiveplusreflectiveimpliesequivalence}
	If $F\dashv G $ is such that $F$ is full co-reflective and $G$ is full reflective, then
	$F\dashv G$ is a $2$-adjoint equivalence.
\end{coro}

\begin{rem}[Idempotent $2$-adjunction vs. full reflective $2$-functor]
It should be noted that there are non-$2$-monadic idempotent $2$-adjunctions.
Remark \ref{thinidempotency} gives a way of constructing easy examples.
For instance, given a $2$-category $\AAA $, the unique $2$-functor
$\AAA \to \mathsf{1} $
has a left $2$-adjoint if and only if $\AAA $ has an initial object. Assuming that $\AAA $ has an initial object and
 $\AAA $ is not ($2$-)equivalent to $\mathsf{1}$,  the $2$-functor
 $\AAA \to \mathsf{1} $ is not a reflective $2$-functor, although the $2$-adjunction is idempotent.

More generally, by Corollary \ref{coreflectiveplusreflectiveimpliesequivalence} any full reflective $2$-functor which is not an equivalence gives
 an example of an idempotent $2$-adjunction such that the left $2$-adjoint
 is not $2$-comonadic. Dually, any non-equivalence full co-reflective $2$-functor
gives an idempotent $2$-adjunction such that the right $2$-adjoint is not
a full reflective $2$-functor.
\end{rem}

\subsection{Kleisli vs. idempotent adjunctions}
Recall that a $2$-adjunction \textit{$F\dashv G $ is Kleisli if the Kleisli comparison $2$-functor} is an equivalence. This fact holds if, and only if, $F$ is essentially surjective on objects. Moreover,
a Kleisli $2$-adjunction is always premonadic, since the Kleisli $2$-category is equivalent to the full sub-$2$-category of free algebras of the $2$-category $\mathcal{T}\textrm{-}\Alg _{\mathsf{s}} $ of the strict algebras of the induced $2$-monad.

It should be noted that, by Proposition \ref{fullreflective2functorcharacterization},
we have that, whenever a $2$-adjunction $F\dashv G $ is idempotent, \textit{$G$ is $2$-premonadic if and only if $G$ is $2$-monadic}. Therefore by Lemma \ref{observationKleisliequivalenttoalgebras} below, this means that, whenever $\mathcal{T}$ is idempotent, the
Kleisli $2$-category is ($2$-)equivalent to the $2$-category  $\mathcal{T}\textrm{-}\Alg _{\mathsf{s}}$.
\begin{lem}\label{observationKleisliequivalenttoalgebras}
	The following statements are equivalent.
	\begin{enumerate}[i)]
		\item The Kleisli $2$-category w.r.t. $\mathcal{T}$ is $2$-equivalent to the $2$-category
		of (strict) $\mathcal{T} $-algebras.
		\item If $F'\dashv G' $ induces $\mathcal{T}$, then $G' $ is $2$-premonadic if, and only if,
		$G' $ is $2$-monadic.
	\end{enumerate}	
\end{lem}

By Proposition \ref{fullreflective2functorcharacterization}, we conclude the following well-known result:

\begin{coro}\label{KelislivsMonadicIdempotentcase}
	An idempotent $2$-adjunction $F\dashv G $ is $2$-monadic if, and only if, it is Kleisli.
\end{coro}

\subsection{Lax idempotency}
For this part, we assume the definition of strict algebras and lax $\mathcal{T}$-morphisms between them, which can be found, for instance,
in \cite[Definition~2.2]{arXiv:1711.02051}. 

Given a $2$-monad $ \mathcal{T} $, we denote by $ \mathcal{T}\textrm{-}\Alg  _{\ell }$
the $2$-category of strict algebras, lax $\mathcal{T}$-morphisms and modifications. In this case, 
$\mathcal{T}\textrm{-}\Alg _{\mathsf{s}}$ is the locally full sub-$2$-category of $ \mathcal{T}\textrm{-}\Alg  _{\ell }$ consisting of strict $\mathcal{T}$-algebras and strict $\mathcal{T}$-morphisms between them.

Theorem \ref{characterizationlaxidempotentmonad} is a well-known characterization of
lax idempotent $2$-monads~\cite{MR1359690}. We refer to \cite{MR1432190, MR1476422} for the proofs.

\begin{theo}[Lax idempotent $2$-monad]\label{characterizationlaxidempotentmonad}
	The following statements are equivalent.
\begin{enumerate}[i)]
		\item $\mathcal{T}$ is lax idempotent;	
		\item $\id _T\ast \eta \dashv \mu $ is a rali adjunction;
		\item $a: T(x)\to x $ is a $\mathcal{T}$-algebra structure if, and only if, there is a \textit{rari} adjunction $ a\dashv \eta_ x $;
		\item $a: T(x)\to x $ is a $\mathcal{T}$-pseudoalgebra structure if, and only if, there is an adjunction $ a\dashv \eta_ x $;
		\item the forgetful $2$-functor $ \mathcal{T}\textrm{-}\Alg  _{\ell }\to  \BBB $ between the $2$-category of strict $\mathcal{T}$-algebras and lax $\mathcal{T}$-morphisms and the $2$-category $\BBB$  is fully faithful.
\end{enumerate}	
\end{theo}

Similarly to the idempotent case, a $2$-adjunction
induces a lax idempotent $2$-monad
if and only if it induces a lax idempotent $2$-comonad.
Furthermore, we give below a lax idempotent analogue
of Proposition \ref{inducedidempotent}.

\begin{theo}[Lax idempotent $2$-adjunction]\label{characterizationlaxidempotent2adjunction}
The following statements are equivalent.
\begin{enumerate}[i)]
	\item $\mathcal{T}$ is lax idempotent;\label{laxidempotentadjunctioncharacterization1}		
	\item $G\varepsilon \dashv \eta G  $ is a lali adjunction;\label{laxidempotentadjunctioncharacterization2}	
	\item $F\eta \dashv \varepsilon F  $ is a rali adjunction;\label{laxidempotentadjunctioncharacterization3}	
	\item The induced $2$-comonad is lax idempotent.\label{laxidempotentadjunctioncharacterization4}
\end{enumerate}
\end{theo}
\begin{proof}
By Lemma \ref{characterizationlaxidempotentmonad}, it is clear that \ref{laxidempotentadjunctioncharacterization2} or \ref{laxidempotentadjunctioncharacterization3}
implies \ref{laxidempotentadjunctioncharacterization1}.
Conversely, assuming \ref{laxidempotentadjunctioncharacterization1},
we have by Lemma \ref{characterizationlaxidempotentmonad} that
$ \id _{GF}\ast \eta \dashv \id _G\ast\varepsilon \ast \id _ F $.
By \textit{doctrinal adjunction} (\textit{e.g.} \cite{MR0360749}), we conclude that
$F\left( \eta _ x\right)\dashv   \varepsilon _{F(x)} $ for every $x$
of $\BBB $. Finally, again, by doctrinal adjunction, we conclude that
$ \id _F\ast \eta\dashv \varepsilon \ast \id _ {F}  $. This proves that
\ref{laxidempotentadjunctioncharacterization1} implies \ref{laxidempotentadjunctioncharacterization3}.

Analogously, by doctrinal adjunction, we get that \ref{laxidempotentadjunctioncharacterization1} implies \ref{laxidempotentadjunctioncharacterization2}.
Hence we proved that the first three statements are equivalent.

Since
the condition \ref{laxidempotentadjunctioncharacterization2} is codual and equivalent to
\ref{laxidempotentadjunctioncharacterization3}, we get that \ref{laxidempotentadjunctioncharacterization1} is equivalent to its codual -- which means  \ref{laxidempotentadjunctioncharacterization4}.
\end{proof}

We follow Kelly's definition of regular epimorphism (also known as strict epimorphism) as outlined in \cite{zbMATH03271557}. We recall that, whenever a ($2$-)category $\AAA$ has kernel pairs, a morphism is a regular epimorphism if and only if its effective (that is, it is the coequalizer of its kernel pair).

\begin{defi}[pre-Kock-Z\"{o}berlein $2$-functor]
If the induced $2$-monad  $\mathcal{T}$ is
 lax idempotent, the $2$-adjunction $F\dashv G $ is
 \textit{lax idempotent}. In this case if, furthermore,
$G$ is $2$-premonadic,
$G$ is called a \textit{pre-Kock-Z\"{o}berlein $2$-functor}.
Finally, if it is also $2$-monadic, $G$ is a \textit{Kock-Z\"{o}berlein $2$-functor}.
\end{defi}

\begin{prop}\label{preKockZoberleintwofunctor}
	Assume that $ F\dashv G : \AAA\to\BBB $
	is lax idempotent. The following statements are equivalent.
	\begin{enumerate}[i)]
		\item $G$ is a pre-Kock-Z\"{o}berlein $2$-functor;	
		\item For each object $x\in \AAA $, $\varepsilon _x $ is a regular epimorphism;
		\item For each object $x\in \AAA $,
%
\pu
\begin{equation}
\diag{coequalizerofthecounitforBecktheorem}
\end{equation}
is a coequalizer.		
	\end{enumerate}
\end{prop}
\begin{proof}
	The result follows directly from the well-known characterization of ($2$-)premonadic ($2$-)functors due to Beck (see, for instance, \cite[p.~226]{MR2056584}).
\end{proof}

\begin{theo}
Assume that $ F\dashv G : \AAA\to\BBB $
is lax idempotent. The following statements are equivalent.
\begin{enumerate}[i)]
		\renewcommand\labelitemi{--}
		\item $G$ is a Kock-Z\"{o}berlein $2$-functor;	
		\item $G$ creates absolute coequalizers;
		\item $G$ is a pre-Kock-Z\"{o}berlein $2$-functor, and, whenever
		$\eta _ y $ is a \textit{rari}, there is $x\in \AAA $
		such that $G(x)\cong y $.
\end{enumerate}
\end{theo}
\begin{proof}
	The result follows from Proposition \ref{preKockZoberleintwofunctor}, and the
	characterization of algebra structures for lax idempotent $2$-monads
	recalled in Theorem \ref{characterizationlaxidempotentmonad}.
\end{proof}

\begin{rem}[Algebras and free algebras]
Corollary \ref{KelislivsMonadicIdempotentcase} says that a $2$-functor $G$
is Kleisli if and only if it is monadic, whenever $F\dashv G$ induces an idempotent $2$-monad.
This is not the case when $\mathcal{T}$ is only lax idempotent. The reference \cite{MR3673245} provides several counterexamples in this direction. Moreover, in our context, in Section \ref{change-of-base functor}, Theorem \ref{laxidempotentcoherencelaxcommacomma} also provides several examples: more precisely, given any $2$-category $\AAA $ and object $z\in\AAA $, the $2$-adjunction between the \textit{lax comma $2$-category $\AAA //z $ (see Definition \ref{definitionoflaxcommacategory})}
and the corresponding comma $2$-category $\AAA{/z} $ usually is a Kleisli $2$-adjunction which is not $2$-monadic.
\end{rem}

Finally, it should be noted that:
\begin{lem}\label{equivalenciadefullreflectivelaxidempotent}
	If $\AAA $ and $\BBB $ are locally discrete, we have that
	$F\dashv G $ is lax idempotent ($G$ is pre-Kock-Z\"{o}berlein) if and only if $F\dashv G $ is
	idempotent ($G$ is full reflective).
\end{lem}
\begin{proof}
It is enough to note that a $2$-monad defined on a locally discrete $2$-category is lax idempotent if and only if it is  idempotent. The rest follows from Proposition \ref{fullreflective2functorcharacterization}. More particularly, it follows from the fact that $2$-premonadicity and $2$-monadicity are equivalent properties for idempotent $2$-adjunctions.	
\end{proof}

\section{Composition of $2$-adjunctions}\label{compositionof2adjunctionssection}
Throughout this section,
%
\pu
\begin{equation}\label{compositionof2adjunctions}
\diag{compositionof2adjunctionsfirstdiagram}
\end{equation}
are given $2$-adjunctions, and $\TTTT = (T, \mu , \eta ) = (GF, G\varepsilon F , \eta ) $ is the $2$-monad induced by the $2$-adjunction $F\dashv G $.
Recall that we have the composition of $2$-adjunctions above given by
%
\pu
\begin{equation}
\diag{compositionof2adjunctionsseconddiagram}
\end{equation}
where   $\RRRR = ( R, v, \alpha ) $
denotes the $2$-monad induced by $FH\dashv JG $.

\subsection{Idempotent $2$-adjunctions}
If $J$ and $G$  are full reflective $2$-functors,
$JG$ is a full reflective $2$-functor and, in particular, $FH\dashv JG $
induces an idempotent $2$-monad. However, if $F\dashv G $ and $H\dashv J $
are only idempotent $2$-adjunctions, we cannot conclude that the composite is
idempotent. For instance,  consider the $2$-adjunctions
%
\pu
\begin{equation}\label{Exampleofcompositionoffullreflectivenotfullreflective}
\diag{compositionoftwoadjunctionscompactHausdorfftopset}
\end{equation}
in which $\Top $ is the locally discrete $2$-category of topological spaces
and continuous functions, $\CmpHaus $ is the full sub-$2$-category of compact Hausdorff spaces, and the right adjoints are the usual forgetful functors. Both
$2$-adjunctions are idempotent, but the composition induces the ultrafilter ($2$-)monad which is not idempotent.

Proposition \ref{1dimensionalsimple}
characterizes when the composition of the $2$-adjunctions is idempotent.
It corresponds to the characterization of the simple (reflective) functors in
the $1$-dimensional case.

\begin{prop}\label{1dimensionalsimple}
	Assume that $F\dashv G$ is idempotent.	
	The following statements are equivalent.
	\begin{enumerate}[i)]
		\item $FH\dashv JG $ is idempotent;\label{idmepotentcompositioncharacterization1}
		\item $JGF\delta G $ (or $F \delta GFH $) is a monomorphism;\label{idmepotentcompositioncharacterization2}
		\item $FH \alpha $ (or $\alpha JG $) is an epimorphism.\label{idmepotentcompositioncharacterization3}
	\end{enumerate}	
\end{prop}
\begin{proof}
	Since $F\dashv G$ is idempotent, $G\varepsilon $, $\varepsilon F$, $F\eta  $ and $\eta G $ are invertible.

	By Proposition         \ref{inducedidempotent},  the $2$-adjunction $FH\dashv JG $ is idempotent if, and only if,
	\begin{center}
	$JG\left(\varepsilon\cdot\left(F\delta{G}\right)\right) =  \left( JG \varepsilon\right)\cdot \left(JT\delta{G}\right) $,  or $\left(\varepsilon\cdot\left(F\delta{G}\right)\right) FH =\left(\varepsilon FH\right) \cdot \left(F\delta TH\right)  $,
	\end{center}
	is a monomorphism. Therefore, since $JG \varepsilon  $ and $\varepsilon FH $ are invertible, we get that $FH\dashv JG $ is idempotent if, and only if, $JT\delta G $, or $ F\delta TH$, is a monomorphism. This proves that
	\ref{idmepotentcompositioncharacterization1} is equivalent to \ref{idmepotentcompositioncharacterization2}.
	
	Finally, \ref{idmepotentcompositioncharacterization1} is equivalent to \ref{idmepotentcompositioncharacterization3} by
	Proposition         \ref{inducedidempotent}.
\end{proof}

\begin{coro}\label{criterion1dimensionalsimple}
If $J$ is full reflective and $F\dashv G $ is idempotent, then the composition
is idempotent.
\end{coro}
\begin{proof}
In this case, since $\delta $ is invertible, we have that $JGF\delta G $ is an isomorphism and, hence, a monomorphism.	
\end{proof}

\begin{defi}[Admissible $2$-functor]\label{Admissibletwoadjunction}
	The $2$-adjunction $F\dashv G $ is \textit{admissible} w.r.t. 	$H\dashv J$ if $JG$
	is a full reflective $2$-functor.
\end{defi}

If $G$ is full reflective, and the composition $JG $ is full reflective,
we generally cannot conclude that $J$ is full reflective. More precisely,
in this case, we have:

\begin{prop}\label{1dimensionalsemileftexact}
Assuming that $G$ is full reflective,
the horizontal composition $F\delta G $ is invertible
if and only if
 the $2$-adjunction $F\dashv G $ is admissible w.r.t. $H\dashv J $.
\end{prop}
\begin{proof}
Since $\varepsilon $ is invertible (by Proposition \ref{fullreflective2functorcharacterization}),
we get that $\left(F\delta{G}\right)$ is invertible if and only if the counit $\varepsilon\left(F\delta{G}\right)$ of
$FH\dashv JG $ is invertible.
By  Proposition \ref{fullreflective2functorcharacterization}, this fact completes the proof.
\end{proof}

\subsection{Lax idempotent $2$-adjunctions}
We turn our attention now to analogous results for the lax idempotent case. The main point is to investigate when the composition of the $2$-adjunctions is lax idempotent and premonadic.

\begin{defi}[Simplicity]\label{SIMPLEtwoadjunction}
	The $2$-adjunction $F\dashv G $ is \textit{simple} w.r.t. 	$H\dashv J$ if the
	composition $FH\dashv JG $ is lax idempotent.
\end{defi}

As a consequence of the characterization of lax idempotent $2$-adjunctions, we get:

\begin{theo}[Simplicity]\label{simplicityofFGwrtJH}
	Assume that $G$ is locally fully faithful. 	The $2$-adjunction
	$F\dashv G $ is simple w.r.t. 	$H\dashv J$ if and only if
	$$\left(\id _{TH}\ast \alpha\right) \dashv \left(\mu\ast \id _{H} \right)\cdot
	\left( \id _ T\ast \delta \ast \id _ {TH}\right) $$
	is a rali adjunction.
\end{theo}
\begin{proof}
By Theorem \ref{characterizationlaxidempotent2adjunction},
we conclude that
the $2$-adjunction $FH\dashv JG $ is lax idempotent
if and only if
$$\left( FH \alpha\right) \dashv \left(\varepsilon FH \right)\cdot
\left( F \delta \ast TH \right) $$
is a rali adjunction.
Since $G$ is locally fully faithful, we have the rali adjunction above if, and only if,
there is a rali adjunction
$ TH\alpha \dashv \left(\mu H \right)\cdot
\left( T \delta TH\right) $.
\end{proof}

The characterization of Theorem \ref{simplicityofFGwrtJH} turns out to be  difficult to apply for most of the examples, since  it involves several units
and counits of the given $2$-adjunctions. Therefore it seems useful to have suitable
sufficient conditions to get simplicity.

\begin{theo}\label{criterionsufficientforsimplicity}
\begin{enumerate}[a)]	
\item	Assume that $JGF\delta G $  is invertible:  $FH\dashv JG $ is lax idempotent if and only if 	there is a lali adjunction $JG\varepsilon\dashv   J\eta G $.
\item Assume that $ F \delta GFH $ is invertible: $FH\dashv JG $ is lax idempotent
if and only if there is a rali adjunction $F\eta H\dashv \varepsilon FH  $.
\end{enumerate}	
\end{theo}
\begin{proof}

We assume that $JGF\delta G $ is invertible. The other case is entirely analogous and,
in fact, dual ($3$-dimensional codual).

By hypothesis, there is a $2$-natural transformation $\vartheta : JGFG\Longrightarrow JGFHJG $
which is the inverse of $JGF\delta G $. Therefore, since
%
\pu
\begin{equation}
\left( JGF\delta G\right)\cdot \left( \alpha JG\right)\quad =\quad  \diag{twocellidentityfortheproofofthecriterionforsimplicity}
\quad  = \quad J\eta G,
\end{equation}
we conclude that
\begin{equation}
\vartheta\cdot \left( J\eta G\right) = \vartheta\cdot \left( JGF\delta G\right)\cdot \left( \alpha JG\right) = \alpha JG.
\end{equation}
Therefore we have the following situation
%
\pu
\begin{equation}
\diag{compositionofadjunctionsepsilondelta}
\end{equation}
in which $\vartheta ^{-1} = JGF\delta{G}  $. This is the hypothesis of Corollary
\ref{precisecancellationisomorphism}  and, thus, there is a lali adjunction
$$JG\left(\varepsilon\cdot\left({F}\delta{G}\right)\right)\dashv \alpha_{JG} $$
if, and only if, there is a lali adjunction $JG\varepsilon\dashv J\eta{G} $.
By Theorem \ref{characterizationlaxidempotent2adjunction}, this completes the proof.
\end{proof}
\begin{coro}\label{corollarytrivialforlaxidempotent}
Assume that $F\dashv G $ is lax idempotent.
\begin{enumerate}[a)]	
	\item	If $JGF\delta G $  is invertible, then  $FH\dashv JG $ is lax idempotent.
	\item If $ F \delta GFH $ is invertible, then $FH\dashv JG $ is lax idempotent.
\end{enumerate}		
\end{coro}
\begin{proof}
In fact, if $F\dashv G $ is lax idempotent, we have in particular that 	
there are a rali adjunction $F\eta H\dashv \varepsilon FH  $ and a lali adjunction $JG\varepsilon\dashv   J\eta G $. Therefore the result follows from Theorem \ref{criterionsufficientforsimplicity}.
\end{proof}

It should be noted that the $2$-adjunctions in \eqref{Exampleofcompositionoffullreflectivenotfullreflective}
show in particular that $FH\dashv JG $ might not be lax idempotent,
even if $F\dashv G$ and $H\dashv J $ are.  However, analogously to the idempotent case (see Corollary \ref{criterion1dimensionalsimple}), we have a nicer situation whenever
$J $ is full reflective.
\begin{coro}
If $J$ is full reflective, then $F\dashv G $ is lax idempotent if, and only if,
 $FH\dashv JG $ is lax idempotent.	
\end{coro}
\begin{proof}
Assuming that $J$ is full reflective, we get that $\delta $ is invertible and, thus, $JGF\delta G $ is invertible.

If $F\dashv G $ is lax idempotent, we get that the composite is lax idempotent
by Corollary \ref{corollarytrivialforlaxidempotent}.
Conversely, if $FH\dashv JG $ is lax idempotent,  by
Theorem \ref{criterionsufficientforsimplicity},
there is a lali adjunction
$$JG\varepsilon\dashv   J\eta G .$$
Since $J$ is locally an isomorphism, this implies that there is a lali adjunction $G\varepsilon\dashv   \eta G $
which proves that $F\dashv G $ is lax idempotent by Theorem \ref{characterizationlaxidempotent2adjunction}.
\end{proof}

\begin{defi}[$2$-admissibility]\label{twoadmissibledefinition}
	The $2$-adjunction $F\dashv G $ is \textit{$2$-admissible} w.r.t. 	$H\dashv J$ if the
	composition $FH\dashv JG $ is lax idempotent and premonadic (that is to say, $JG$  is pre-Kock-Z\"{o}berlein).
\end{defi}

As a consequence of Proposition \ref{preKockZoberleintwofunctor}
and Theorem \ref{simplicityofFGwrtJH}, we have:

\begin{theo}[$2$-admissibility]\label{admissibilitywrtJH}
	Assume that $G$ is pre-Kock-Z\"{o}berlein. The $2$-adjunction
	$F\dashv G $ is \textit{$2$-admissible} w.r.t.  $H\dashv J$
	if, and only if, the two conditions below hold.
	\begin{itemize}
		\renewcommand\labelitemi{--}
	    \item $TH \alpha \dashv \left(\mu\ast \id _{H} \right)\cdot
	    \left( \id _ T\ast \delta \ast \id _ {TH}\right)$ is a lari adjunction (or, equivalently, $F\dashv G $ is simple w.r.t. $H\dashv J $);
		\item For each object $z\in\CCC $, $\left(\varepsilon\cdot \left( F\delta G\right)\right)_z $
		is a regular epimorphism.
	\end{itemize}
\end{theo}
Recall that the composition of a regular epimorphism
with a split epimorphism is always a regular epimorphism~\textit{c.f.}~\cite{zbMATH03271557}. 
Therefore we also have that:

\begin{coro}
$F\dashv G $ is simple w.r.t. $H\dashv J $, and $F\delta G $ is a
split epimorphism, we conclude that $F\dashv G $ is $2$-admissible w.r.t. $H\dashv J $.	
\end{coro}
\begin{proof}
It follows directly from Theorem \ref{admissibilitywrtJH} and the observation above.	
\end{proof}

Since the composition of a regular epimorphism with an isomorphism is always a
regular epimorphism, we get:
\begin{coro}\label{themaincase}
If $F\delta G $ is an isomorphism and $G$ is pre-Kock-Z\"{o}berlein, then
$F\dashv G $ is $2$-admissible w.r.t. $H\dashv J $. In particular, if $J$ is full reflective and
$G$ is pre-Kock-Z\"{o}berlein, we conclude that $JG$ is pre-Kock-Z\"{o}berlein.
\end{coro}
\begin{proof}
Since $F\delta G $ is invertible, we get that $JG F\delta G $ is invertible. Therefore, by Corollary \ref{corollarytrivialforlaxidempotent}, we get the simplicity. Moreover
$\varepsilon \cdot\left( F\delta G\right) $
is a regular epimorphism since $\varepsilon $ is a regular epimorphism and
$\left( F\delta G\right) $ is invertible.	
\end{proof}

\section{Lax comma $2$-categories and change-of-base $2$-functors}\label{change-of-base functor}\label{basicdefinitionsarticle}
The notion of lax comma $2$-categories is well-known and has been considered in the literature in many contexts (see, for instance, \cite[I,5]{MR0371990}, \cite[\S~6]{MR0249483}, 
\cite[Exercise~5, p.~115]{MR1712872} or \cite[p.~305]{MR558494}). 
 We recall the definition in an elementary manner below, following the perspective of our setting.   
The main aim is to introduce the respective notions of change-of-base $2$-functors.

In a $2$-category 
$\AAA $ with products, given any object $y\in \AAA $, the endofunctor 
$$ \left( y\times - \right) : \AAA \to \AAA $$
has  a unique comonadic structure. 
In this context, the usual $2$-category of strict coalgebras
$\left({y}\times{-}\right)\textrm{-}\CoAlg_{\mathsf{s}}  $ is isomorphic to the
comma $2$-category $\AAA / y$. Moreover, the $2$-category $\left( y\times - \right)\textrm{-} \CoAlg _{\ell }$ of strict $ \left( y\times - \right)$-coalgebras and lax morphisms is what we call the \textit{lax comma $2$-category} $\AAA // y $.\footnote{See, for instance, \cite[Def.~4.1]{2016arXiv160703087L}.} More generally, we explicitly define the lax comma categories below.

\begin{defi}[Lax comma $2$-category]\label{definitionoflaxcommacategory}
	Given an object $y$ of a $2$-category $\AAA $, we denote by
	$\AAA // y  $ the $2$-category defined by the following.
	\begin{itemize}
		\renewcommand\labelitemi{--}	
		\item The objects are pairs $ (w, a) $ in which $w$ is an object
		of $\AAA $ and
		$$\xymatrix{ w\ar[r]|-{a} &y
		}$$
		is a morphism of $\AAA $.
		\item A morphism in $\AAA // y  $  between objects $(w,  a) $ and $ (x, b) $    is a pair
%
\pu	

%
\pu		
$$\left( \diag{laxmorphismofcoalgebras1cell}, \diag{twocelloflaxcommamorphism}        \right) $$
in which
$f: w\to x $ is a morphism of $\AAA $ and  $\livb $ is a $2$-cell of $\AAA $.
	
If $(f, \livb ): (w,  a)\to (x,  b)  $ and $ (g, \liva ): (x, b)\to (z, c) $ are morphisms  of $\AAA // y $, the composition is defined by
$(g\circ f, \livb\cdot \left(\liva \ast \id _f\right) ) $, that is to say,
the composition of the morphisms $g$ and $f$ with the pasting
%
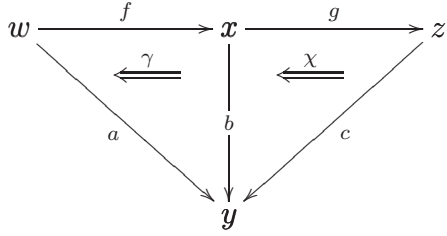
\pu		
$$\diag{compositionofmorphismslaxcommatwocategories} $$
of the $2$-cells $\liva $ and $\livb $. Finally, with the definitions above, the identity on the object $(w,a)$ is of course the morphism $(\id _w, \id _a ) $.
\item A $2$-cell between morphisms $(f, \livb  )  $ and $(f', \livb ' )  $
is given by a $2$-cell $\livc : f\Rightarrow f' $ such that the
equation
%
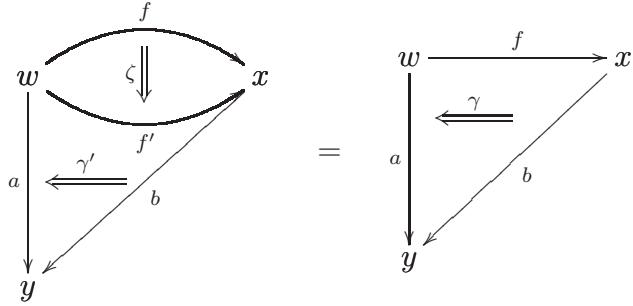
\pu	
%
\pu		
\begin{equation*}
\diag{leftsideequationtwocellforlaxcommacategor}\quad =\quad \diag{rightsideequationtwocellforlaxcommacategory}
\end{equation*}
holds.
	\end{itemize}
The $2$-category $\AAA // y$ is called the \textit{lax comma $2$-category} of $\AAA $ over $y $, while the
$2$-category $\AAA ^\co //y $ is called the 	\textit{colax comma $2$-category} of $\AAA $ over $y$.
\end{defi}

The concept of (co)lax comma $2$-category, possibly under other names, has already appeared in the literature. See, for instance, \cite[Exercise~5, p.~115]{MR1712872} or \cite[p.~305]{MR558494}.
As for our choice of the direction of the $2$-cells for the notion of lax comma $2$-categories, although we do not follow \cite[p.~305]{MR558494}, our choice is compatible with the usual definition of lax
natural transformation.

\begin{defi}[(Strict) comma $2$-category]
	Given an object $y$ of a $2$-category $\AAA $, we denote by
	$\AAA / y  $ the \textit{comma $2$-category} over $y$, defined
	to be the locally full \textit{wide} sub-$2$-category
	of 	$\AAA // y  $ in which
	a morphism from
	$\left( w,  a\right) $ to $ \left( x, b\right) $   is a morphism $$(f, \liva ): \left(w,  a\right) \to \left(x, b\right) $$  such that $\liva $ is the identity $2$-cell.
\end{defi}

\begin{rem}\label{Fcategoryforthefirsttime}
We have an inclusion $2$-functor $\AAA / y\to \AAA // y    $ obviously defined. The morphisms in the image of this inclusion
are called \textit{strict} (or \textit{tight}) morphisms of $ \AAA // y $. The $2$-category $ \AAA // y $ endowed with this inclusion
forms an enhanced $2$-category, or, more precisely, an $\mathfrak{F}$-category as defined in \cite{MR2854177}.
\end{rem}

\subsection{Classical (strict) change-of-base functor}
Assuming that $\AAA $ has pullbacks, given any morphism $c:y\to z $ of a $2$-category $\AAA $, it
is well known that it
induces a $2$-adjunction
%
\pu
\begin{equation}\label{equationtwoadjunctionchangeofthebasetwofunctorbasic}
\diag{changeofthebasecommacategories}
\end{equation}
between the (strict) comma $2$-categories
in which the right $2$-adjoint is called the \textit{change-of-base $2$-functor} induced by the morphism $c$ (see, for instance, \cite{MR1173011}).
Recall that $c^\ast $ is defined by the pullback along $c$,
and the left adjoint is defined by the composition with $c$, the so called direct-image $2$-functor $c! (w, a) = (w, ca)$.

In the present section, we give the analogue for lax comma $2$-categories, that is to say, the\textit{ change-of-base $2$-functors for the lax comma $2$-categories}, given in Proposition \ref{def:change-of-base-lax}. Firstly, we recall the classical case:

\begin{prop}[Change-of-base $2$-functor]\label{teoremamudancadebase}
	Let $\AAA $ be a $2$-category with pullbacks. If $c: y\to z $ is any morphism, we get a $2$-adjunction
%
\pu
\begin{equation}\label{equationtwoadjunctionchangeofthebasetwofunctorpullbackbasic}
	\diag{changeofthebasecommatwocategories}
\end{equation}
	in which $c^\ast $ is defined by the pullback along $c$.
	Explicitly, the assignment of objects of $c^\ast $ is given by $$(w, a) \mapsto  (w\times_{(a,c)}y, c^\ast (a) : w\times_{(a,c)}y \to y )$$
	while the action of $c^\ast $ on morphisms
	is given by
\begin{equation}\label{definitionofactiononmorphismstwofunctorcoalgebraslaxcommaa}
\left( w\xrightarrow{f} x,\, \id _a \right) : (w, a)\to (x, b) \,\,\mapsto\,
	\left( w\times_{(a,c)}y \xrightarrow{c^\ast\left(f, \id_a\right)} x\times_{(b,c)}y,\, \id _{c^\ast (a)}  \right) : c^\ast (a)\to c^\ast( b)
\end{equation}	
in which all the squares of
%
\pu

\begin{equation}\label{definitionpullbackonmorphismss}
\diag{cstaronmorphismm}
\end{equation}
are pullbacks. Finally, the image of a $2$-cell $\livc : f\Rightarrow f' : (w,a)\to (x,b) $ is defined by the unique $2$-cell $c^\ast \left(\livc \right)$ such that the equations
%
\pu
%
\pu

%
\pu
%
\pu
\begin{equation}
\diag{leftsidefirstequationimagetwocellofthepbchangeofbase} = \diag{rightrsidefirstequationimagetwocellofthepbchangeofbase}\quad\mbox{ and }\quad \diag{leftsidesecondequationimagetwocellofthepbchangeofbase} =\quad
\diag{rightrsidesecondequationimagetwocellofthepbchangeofbase}
\end{equation}
hold.
\end{prop}

By considering the comma object, every morphism $c: y\to z $ in a $2$-category $\AAA $ induces a $2$-functor $ c^\Leftarrow : \AAA //z \to \AAA / y  $. This gives what we call, herein, the (lax) \textit{change-of-base $2$-functor}. We start by defining the appropriate analogous of direct image $2$-functor in our setting, that is to say, the $2$-functor $c\! : \AAA / y \to \AAA / / z $
induced by each morphism $c: y\to z $ of $\AAA $. 
\begin{defi}[Direct image]\label{definitiondirectimage}
	If   $c: y\to z $ is any
	morphism of a $2$-category $\AAA $, we define the commutative diagram
%
\pu
\begin{equation}\label{directimagedefinitiondiagramequation}
	\diag{directimagedefinitiondiagram}
\end{equation}
	in which the unlabeled arrow is the obvious inclusion, and
	$$c! : \AAA /y \to \AAA / z $$
	is defined by
	$$(x, a)\mapsto (x, ca), \, (f, \id  )\mapsto (f, \id _ c \ast \id  ), \,
	\livc \mapsto \livc , $$
	that is to say, the usual
	direct image $2$-functor.
\end{defi}

Theorem \ref{teoremamudancadebaselassa} shows that, in the presence of suitable comma objects, for each morphism $c: y\to z $ in $2$-category, the direct image $2$-functor $c\! : \AAA / y \to \AAA // z $ has a right $2$-adjoint -- the (lax) change-of-base $2$-functor $c^\Leftarrow $ defined below. 

\begin{defi}[$c^\Leftarrow $]\label{def:change-of-base-lax}
	Let $\AAA $ be any $2$-category, and $c: y\to z $ a morphism of $\AAA $. Assume that $\AAA $ has comma objects along $c$. We denote by
	$$ c^\Leftarrow : \AAA //z \to \AAA / y                  $$
	the $2$-functor
	defined by the comma object along the morphism $c$. Explicitly, the action on objects of $c^\Leftarrow $
	is given by
	\begin{equation}
	(x, b) \mapsto   (b\downarrow c, c^\Leftarrow (b) : b\downarrow c\to y )
	\end{equation}
in which
%
\pu
\begin{equation}\label{definitionofthetwofunctorcommaobject}
\diag{commadiagramdefinitioncommaadjunction}
\end{equation}
is the comma object as in \ref{definitionofcommaobjects},
while the action on morphisms
is given by

%
\pu

%

\pu
\begin{equation}\label{definitiononmorphismsofcommaalongc}
\left( \diag{morphismfforthedefinitioncLeftarrow}, \diag{twocelloflaxcommamorphismmmm}\right)\quad\mapsto\quad
\left( a\downarrow{c}\xrightarrow{c^\Leftarrow\left(f, \livb\right)  } b\downarrow{c},\quad\id _{c^\Leftarrow (a) }  \right)
\end{equation}
in which $c^\Leftarrow (f, \livb )  $, sometimes only denoted by $c^\Leftarrow (f )$,   is the unique morphism of $\AAA $ such that the equations
$$ b^\Rightarrow (c)\,\cdot\,{c^\Leftarrow{\left( f, \livb \right)}} = f\,\cdot\,{a^\Rightarrow{(c)}}, \qquad c^\Leftarrow (b)\,\cdot\,{c^\Leftarrow{\left( f, \livb \right)}} = c^\Leftarrow (a), $$
%
\pu
%
\pu
\begin{equation}
\diag{cleftarrowonmorphismm}\quad =\quad \diag{cleftarrowonmorphismmrightside}
\end{equation}
hold. Finally, if $\livc : f\Rightarrow f' : (w,a)\to (x,b) $ is a $2$-cell between morphisms $(f, \livb) $ and $(f', \livb ' ) $ in $\AAA // z $, the $2$-cell
$c^\Leftarrow (\livc ) $ is the unique $2$-cell such that the equations
%
\pu
%
\pu

%
\pu
%
\pu
\begin{equation}
\diag{leftsidefirstequationimagetwocellofthecommachangeofbase} = \diag{rightrsidefirstequationimagetwocellofthecommachangeofbase}\quad\mbox{ and }\quad \diag{leftsidesecondequationimagetwocellofthecommachangeofbase} =\quad
\diag{rightrsidesecondequationimagetwocellofthecommachangeofbase}
\end{equation}
hold.
\end{defi}

\begin{theo}\label{teoremamudancadebaselassa}
	Let $\AAA $ be any $2$-category, and $c: y\to z $ a morphism in $\AAA $. If $\AAA $ has comma objects along  $c $, then we
	have a $2$-adjunction
%
\pu
\begin{equation}\label{changeofthebaselaxcommacategoriesequation}
\diag{changeofthebaselaxcommacategories}.
\end{equation}
\end{theo}

\begin{proof}
	We define below the counit, denoted by $\delta $,  and the unit, denoted by $\rho $,  of the $2$-adjunction $c\!\dashv c^\Leftarrow $.	
	
	For each object
	\begin{equation*}
		\left( x,   x\xrightarrow{b} z           \right)
	\end{equation*}
	of $\AAA // z  $,  we have the
	comma object
%
\pu	
\begin{equation}
\diag{commadiagramdefinitionproofcommaobject}
\end{equation}
as in \eqref{definitionofthetwofunctorcommaobject}. We define the counit on $(x,b) $, denoted by $\delta_ {(x,b)} $, to be
the morphism between $c\! c^\Leftarrow (x,b) $ and $(x, b) $ in $\AAA //z $ given by the pair $(b^\Rightarrow (c) , \liva ^{b\downarrow c} ) $.
	
Moreover, for each object
	\begin{equation*}
		\left( w,   w\xrightarrow{a} y        \right)
	\end{equation*}
	in $\AAA /y $, we have the comma object
%
\pu	
\begin{equation}
\diag{commaalongca}
\end{equation}	
	in $\AAA $. By the universal property of the comma object, there is a
	unique morphism $\rho _ {(w, a)}  '$ of $\AAA $
	such that the equations
%
\pu
\begin{equation}
\diag{morphismcommadiagramdefinitionrholinha}
\end{equation}
	\begin{equation*}
		\left(ca\right)^\Rightarrow (c)\, \cdot \, \rho _ {(w, a)}' = \id _w \qquad\mbox{and}\qquad
		c^\Leftarrow c\!  (a)\,\cdot\, \rho _ {(w, a)} '= a
	\end{equation*}
	hold.
	
	By the equation above, the pair $ (\rho _ {(w, a)}' , \id _a  ) $
	gives a morphism between $(w,a) $
	and $(ca\downarrow c, c^\Leftarrow c\!  (a)) $
	in $\AAA / y $.
	We claim that the component $\rho _ {(w, a)}  $ of the unit  of $c\! \dashv c^\Leftarrow $  on $(w, a) $ is  the morphism defined by the pair $ (\rho _ {(w, a)}' , \id _a  ) $.
	
	It is straightforward to see that the definitions above actually give $2$-natural transformations $\delta : c\! c^\Leftarrow\longrightarrow \id _{\AAA//z } $ and $\rho : \id _ {\AAA /y }\longrightarrow c^\Leftarrow c\!  $.
	We prove below that $\delta $ and $\rho $ satisfy the triangle identities.
	
Let $(w,a) $ be an object of $\AAA / y $.

The image of the morphism $\rho _{(w, a)} $
	by the $2$-functor $c\!  : \AAA /y\to \AAA //z  $ is the morphism
	$ (\rho _{(w,a)} ', \id _{ca} ) $ between $c\! (w, a ) =  (w, ca) $ and  $(ca\downarrow c, c\! c^\Leftarrow c\!  (a)) $  in $\AAA //z  $, while the component
	$\delta _ { c\! (w,a) } = \delta _ {(w, ca)}$ is the morphism
	$\left( (ca)^\Rightarrow (c), \liva ^{ca\downarrow c} \right) $.
	
	By the definition of $\rho ' _{(w,a)} $, we have that $(ca)^\Rightarrow (c) \cdot \rho _{(w,a)} ' = \id _w $
	and $        \liva ^{ca\downarrow c}   \ast \id _{\rho _{(w,a)} '} = \id _{ca} $. Therefore
	$\delta _ { c\! (w,a) }\, \cdot  \, c\!  \left( \rho _{(w,a)} \right)  $ is the identity on $c\!  (a) $.
	This proves the first triangle identity.
	
	Let $(x, b) $ be an object of $\AAA //z $. Denoting by  $(c\cdot c^\Rightarrow (b)\downarrow c , \liva ^{c\cdot c^\Leftarrow (b)\downarrow c } ) $
	the comma object of $c\cdot c^\Leftarrow (b) $ along $c $, we have that the morphism
	$$ c^\Leftarrow \left(\delta _{(x,b)}\right) : c^\Leftarrow c\! c^\Leftarrow (x,b) \to c^\Leftarrow (x,b) $$
	in $\AAA / y $
	is defined by the pair $(\updelta ' , \id _{c^\Leftarrow c\!  c^\Leftarrow (b) }  ) $ in which $\updelta '$
	is the unique morphism in $\AAA $
	making the diagrams
		\begin{equation*}
	\vcenter{
		\xymatrix@=4em{
			x
			&&
			b\downarrow c
			\ar[ll]_-{b^\Rightarrow (c) }
			\\
			&
			c\cdot c^\Leftarrow (b) \downarrow c
			\ar[ru]_-{\updelta ' }
			\ar[lu]^-{b^\Rightarrow (c) \cdot \left( c\cdot c^ \Leftarrow (b) \right) ^\Rightarrow (c) }
			&
		}
	}
	\qquad
	\vcenter{
		\xymatrix@=3em{
			b\downarrow c
			\ar[rd]^-{c^\Leftarrow (b) }&
			\\
			&
			y
			\\
			c\cdot c^\Leftarrow (b)\downarrow c
			\ar[ru]_-{c^\Leftarrow c\!  c^\Leftarrow (b) }
			\ar[uu]^-{\updelta ' }&
		}
	}
	\end{equation*}
commute, and the equation	
%
\pu
%
\pu
\begin{equation}
\diag{definitionofdeltalinhaleft}\quad =\quad \diag{definitionofdeltalinharight}
\end{equation}
holds.

	Since, by the definition of $\rho $, the underlying morphism $\rho _{c^\Leftarrow (x,b) }' $ of the component of $\rho $ on $ c^\Leftarrow (x, b) $ is such that
	the equations
	\begin{equation*}
		\liva ^ {c\cdot c^\Leftarrow (b)\downarrow c} \ast  \id _ {\rho ' _{c^\Leftarrow (x,b) } }
		= \id _{c\cdot c^\Leftarrow (b) },\,
		\left( c\cdot c^ \Leftarrow (b) \right)^\Rightarrow (c)\cdot \rho ' _{c^\Leftarrow (b) }
		= \id _{b\downarrow c },
		\, c^\Leftarrow c\!  c^\Leftarrow (b) \cdot \rho ' _{c^\Leftarrow (x, b) } = c\cdot c^\Leftarrow (a)
	\end{equation*}
	hold, we get that the equations
%
\pu
\begin{equation}
\diag{secontriangleidentitycoomaobjecttwoadjunction}
\end{equation}
	\begin{equation*}
		c^\Leftarrow (b)\cdot \delta ' \cdot \rho _{c^\Leftarrow (b) }' = c^\Leftarrow (b),
		\qquad
		b^\Rightarrow (c) \cdot \delta ' \cdot \rho _{c^\Leftarrow (b) }' =  b^ \Rightarrow (c)
	\end{equation*}
	hold. Since, by the universal property of the comma object of $b$ along $c$,
	the morphism satisfying the three equations above is unique,
	we conclude that $\delta ' \cdot \rho _{c^\Leftarrow (x,b) }'$ is the identity
	on $b\downarrow c $. This proves that
	$$c^\Leftarrow (\delta _{(x,b)} ) \cdot   \rho _{c^\Leftarrow (x,b) } = \id _{c^\Leftarrow (x,b)} $$
	which proves the second triangle identity.
\end{proof}

\begin{coro}\label{itisactuallythecoherencetwoadjunctioncorollary}
	If $\AAA $ has comma objects (along identities), then
	$\AAA /y \to \AAA // y  $ has a right $2$-adjoint
	which is defined by the comma object along the identity $\id _ y $.
\end{coro}
\begin{proof}
	It follows from Theorem \ref{teoremamudancadebaselassa} and the fact that the
	inclusion 	$\AAA /y \to \AAA // y  $  is actually given by the
	$2$-functor
	$\id _ y \!  : \AAA /y\to \AAA // y  $
	and, hence, it is left $2$-adjoint to the $2$-functor
	$$\id _y ^\Leftarrow : \AAA //y \to\AAA /y .$$
\end{proof}

By Theorem \ref{teoremamudancadebaselassa} and the fact that,
given a morphism $c: y\to z $ of a $2$-category $\AAA $,
%
\pu
\begin{equation}\label{basicdiagramofcompositionwithcandinclusionequation}
\diag{basicdiagramofcompositionwithcandinclusion}
\end{equation}
commutes, we get that:
\begin{theo}\label{relationofchangeofbasecomma}
	Let $\AAA$ be a $2$-category, and $c: y\to z $ a morphism of $\AAA $. If $\AAA $
		has comma objects and pullbacks along $c$, we have the following commutative diagram of $2$-adjunctions
%
\pu
\begin{equation}\label{compositionoftwoadjunctionschangeofbasecommaequation}
\diag{compositionoftwoadjunctionschangeofbasecomma}
\end{equation}		
	which means that the composition of the $2$-adjunction
	$c!\dashv c^\ast : \AAA /z \to \AAA /y $ with
	$\id _z\!\dashv id _z^\Leftarrow : \AAA //z \to \AAA /z $
	is, up to $2$-natural isomorphism, the $2$-adjunction $$c\!\dashv c^\Leftarrow : \AAA //z  \to \AAA /y .$$
\end{theo}

Given a $2$-category $\AAA $, it is clear that, for any object $y$ of $\AAA $, the $2$-adjunction
$ \id _y ! \dashv \id_y ^\ast : \AAA /y  \to \AAA /y $
is $2$-naturally isomorphic to the identity $2$-adjunction $ \id_{\AAA /y }\dashv \id_{\AAA /y  } $
and, in particular, is an idempotent $2$-adjunction.

In the setting of Theorem \ref{teoremamudancadebaselassa}, that is to say, the comma version of the change-of-base $2$-functor, the $2$-adjunction
$$ \id _y \! \dashv \id_y ^\Leftarrow : \AAA //y  \to \AAA /y ,$$
is far from being isomorphic to the identity $2$-adjunction. It is not even idempotent
in most of the cases. It is, however, always lax idempotent and a Kleisli $2$-adjunction. More precisely:

\begin{theo}\label{laxidempotentcoherencelaxcommacomma}
Let $\AAA $ be a $2$-category, and $y $ an object of $\AAA $. If $\AAA $ has comma objects along $\id _y $, then
the $2$-adjunction
%
\pu
\begin{equation}\label{changeofthebaseidentitycoherencelaxidempotentequationn}
\diag{changeofthebaseidentitycoherencelaxidempotent}
\end{equation}
is lax idempotent. Moreover, it is a Kleisli $2$-adjunction and, hence, $\id_y^\Leftarrow $ is a pre-Kock-Z\"{o}berlein
$2$-functor.
\end{theo}
\begin{proof}
In order to verify that \eqref{changeofthebaseidentitycoherencelaxidempotentequationn} is a Kleisli $2$-adjunction, it is enough to see that $\id_y\! $ is bijective on objects. In particular, we conclude that $\id_y^\Leftarrow $ is $2$-premonadic. Therefore, in order to prove that $\id_y^\Leftarrow $ is 	a pre-Kock-Z\"{o}berlein $2$-functor, it remains only to prove that the $2$-adjunction \eqref{changeofthebaseidentitycoherencelaxidempotentequationn} is lax idempotent.

We prove below that
\begin{equation}\label{proofofthelaxidempotencyofthetwoadjunctionoftheidentity}
\id_{\id_y\!} \ast \rho\dashv \delta \ast\id _{\id_y\! }
\end{equation}
is a rari adjunction and, hence, it satisfies the condition \ref{laxidempotentadjunctioncharacterization3} of Theorem \ref{characterizationlaxidempotent2adjunction}, which implies that the $2$-adjunction \eqref{changeofthebaseidentitycoherencelaxidempotentequationn} is lax idempotent.

For short, throughout this proof, we denote $\id_{\id_y\!} \ast \rho$ by $ \overline{\rho }  $, and $\delta \ast\id _{\id_y\! }$ by $\overline{\delta } $.

Recall that, given an object $(x, b)\in{\AAA /y} $, we have that $\overline{\delta }  _{(x,b) } $
is defined by the pair
%
\pu
\begin{equation*}
\left(\overline{\underline{\delta}}_{(x,b)},\,    \diag{counittwocelloftheidentitytwoadjunctionlaxcoomacomma}\right)
\end{equation*}
in which, as suggested by the notation, the $2$-cell is the comma object in $\AAA $, and  $$\overline{\underline{\delta}}_{(x,b)}:=  b^\Rightarrow{(\id_y)}.$$
Moreover, recall that, given an object $(x, b)\in{\AAA /y} $, we have that
 $\overline{\rho } _{(x,b) } = \left(\overline{\underline{\rho}}_ {(x,b)}, \id_b \right)  $ in which $\overline{\underline{\rho}}_ {(x,b)}$
 is the unique morphism of $\AAA $ such that the equations
%
\pu
%
\pu
\begin{center}
	$\overline{\underline{\delta}}_{(x,b)}\cdot \overline{\underline{\rho}}_{(x,b)}  = \id _x$,\quad   $ \id_y^\Leftarrow{(b)} \cdot \overline{\underline{\rho}}_{(x,b)}  =b$,\quad
	and
\end{center}	
\begin{equation}\label{rhodef}
\diag{identityonbtodefinerho} \quad = \quad  \diag{rholparaademonstracaodelaxidempotency}
\end{equation}
hold.

For each object $(x,b)\in\AAA /y $,   the pair of $2$-cells $\left( \liva^{b\downarrow{\id_y}}, \id _ {\overline{\underline{\delta}}_{(x,b)}} \right) $ satisfies the equation
%
\pu
%
\pu
\normalsize
%
\pu
%
\pu
\begin{equation}
\diag{firsttwocellinordertodefinecounitoftheadjunctionthatgivesarari}\, =\,\diag{secondtwocellinordertodefinecounitoftheadjunctionthatgivesarari}
\end{equation}
and, hence,  by the universal property of the comma object, there is a unique $2$-cell $\Gamma _ {(x,b)} $ such that the equations
$$\id_{\id_y^\Leftarrow{(b)} }\ast  \Gamma _ {(x,b)} = \liva^{b\downarrow{\id_y}} \quad\mbox{ and }\quad \id _ {\overline{\underline{\delta}}_{(x,b)}}\ast  \Gamma _ {(x,b)} =  \id _ {\overline{\underline{\delta}}_{(x,b)}} $$
hold. The $2$-cells $\Gamma _ {(x,b)} $ define a modification
$$\Gamma : \overline{\rho }\cdot \overline{\delta }\Longrightarrow \id _{\id_y\!\id_y^\Leftarrow\id_y\!} $$
which we claim to be the counit of the adjunction \eqref{proofofthelaxidempotencyofthetwoadjunctionoftheidentity}.

The first triangle identity holds, since, by the definition of $\Gamma $ above,
$$ \id _ {\overline{\underline{\delta}}_{(x,b)}}\ast  \Gamma _ {(x,b)} =  \id _ {\overline{\underline{\delta}}_{(x,b)}} $$
for every object $(x,b)\in\AAA /y $.

Finally, for each object $(x,b)\in\AAA /y $,
$\Gamma _ {(x,b)}\ast \id _{\overline{\underline{\rho}}_{(x,b)} } $ is such that
$$\id_{\id_y^\Leftarrow{(b)} }\ast  \Gamma _ {(x,b)}\ast \id _{\overline{\underline{\rho}}_{(x,b)} } = \liva^{b\downarrow{\id_y}}\ast \id _{\overline{\underline{\rho}}_{(x,b)} } = \id _b $$
by \eqref{rhodef}, and, of course,
$$\id _ {\overline{\underline{\delta}}_{(x,b)}}\ast  \Gamma _ {(x,b)}\ast \id _{\overline{\underline{\rho}}_{(x,b)} }  =  \id _ {\overline{\underline{\delta}}_{(x,b)}\cdot \overline{\underline{\rho}}_{(x,b)} } .$$
Therefore, by the universal property of the comma object $b\downarrow \id_y $, we get that $\Gamma _ {(x,b)}\ast \id _{\overline{\underline{\rho}}_{(x,b)} } = \id_{\id _{\overline{\underline{\rho}}_{(x,b)} }} $. This completes the proof that the second triangle identity holds.
\end{proof}

\section{Admissibility}\label{secao de admissibilidade}
Throughout this section,
$$\diag{basicadjunction} $$
is a given $2$-adjunction. By abuse of language, given any $2$-functor $H: \AAA\to \BBB $, for each object $x $ in $\AAA $, we denote by the same $\check {H} $ the $2$-functors
$$  \check{H} :  \AAA/x \to \BBB /H(x), \quad
\check{H} :  \AAA /x  \to \BBB // H(x), \quad  \check{H} :  \AAA // x  \to \BBB // H(x) $$
pointwise defined  by  $H$. Moreover, given a morphism $f:w\to x $ of $\AAA $, we denote by
$$ f\l! : \AAA //w\to \AAA //x $$
the $2$-functor defined by the \textit{direct image} between the lax comma $2$-categories, whose restriction to
$ \AAA /w $ is equal to $f\! $.

\begin{prop}
	If $G$ is a locally fully faithful $2$-functor then, for each object
	$x $ of $\AAA$, both
	$\check {G}: \AAA /x\to \BBB/ G(x)
	$ and $\check {G}: \AAA //x \to \BBB// G(x)
	$
	are locally fully faithful.
\end{prop}

\begin{theo}
	For any object $y\in \AAA $, we have two $2$-adjunctions
%
\pu
%
\pu
\begin{equation}
\diag{basicliftingtwoadjunction}\quad\mbox{ and }\quad\diag{basicliftingtwoadjunctionlaxcommatwocategory}
\end{equation}
where the counit and the unit of these $2$-adjunctions are defined
	pointwise by the counit and unit of $F\dashv G $.
\end{theo}

\begin{coro}\label{followsfromthepointwiseproperty}
	For each object $y\in \AAA $, the $2$-adjunctions
%
\pu
%
\pu 	
\begin{equation}
\diag{basicliftingtwoadjunctionn}\quad\mbox{ and }\quad\diag{basicliftingtwoadjunctionlaxcommatwocategoryy}
\end{equation}
are lax idempotent (premonadic) if, and only if, $F\dashv G $ is lax idempotent (premonadic).
\end{coro}

Henceforth, we further assume that $\BBB $ has
comma objects and pullbacks whenever necessary. Recall that, in this case,
by Section \ref{change-of-base functor}, for each object $y $ of $\BBB $,
we have $2$-adjunctions
\begin{center}
	$\eta_ y !\dashv  \eta _y ^\ast : \BBB / GF (y)\to \BBB/ y $ and $\eta_ y \!\dashv  \eta _y ^{\Leftarrow } : \BBB // GF (y) \to \BBB/ y   $
\end{center}
in which the right $2$-adjoints are given respectively by the pullback
and the comma object along $\eta_ y$.

\begin{defi}[Simple, admissible and $2$-admissible $2$-functors]\label{maindefinition}
The $2$-functor $G$ is called \textit{simple/$2$-admissible} if $F\dashv G $ is lax idempotent/pre-Kock-Z\"{o}berlein, and, for every $y\in \BBB$,
%
\pu
\begin{equation}
\diag{basicliftingtwoadjunctionlaxcommatwocategoryagain}
\end{equation}
is simple/$2$-admissible w.r.t. $\eta_ y !\dashv  \eta _y ^\Leftarrow $ (see Definitions \ref{SIMPLEtwoadjunction} and \ref{twoadmissibledefinition}).
	
	We say that $G$ is \textit{admissible
		w.r.t. the basic fibration}
	if $G$ is fully faithful, and, for every $y\in \BBB$,
%
\pu
\begin{equation}
\diag{basicadjunctionliftingagain}
\end{equation}	
is admissible w.r.t. $\eta_ y !\dashv  \eta _y ^\ast $.
\end{defi}

\begin{rem}
The notion of admissibility w.r.t. the basic fibration is
just the direct strict $2$-dimensional generalization of the classical notion of admissibility (also called semi-left-exact reflective functor)~\cite{MR779198, MR1822890}, while the notion of simplicity coincides with that introduced  in \cite{MR3545937}.	
\end{rem}

In order to establish the direct consequences of the results of Section \ref{compositionof2adjunctionssection} for
the case of $2$-admissibility and simplicity,
we set some notation below.
For each $y$ of $\BBB $, we consider the $2$-adjunctions
%
\pu
\begin{equation}\label{SIMPLECOMPOSITION}
	\diag{compositionof2adjunctionsfirstdiagramlaxcommatwocategorycomma}
\end{equation}
in which, by abuse of language, we denote respectively by $\varepsilon $ and $\eta $
the counit and unit defined pointwise, and $\TTTT = (T, \mu , \eta ) $ the
$2$-monad induced by $\varepsilon _ {F(y)}! \circ \check{F}\dashv \check{G} $.

In this case, the composition of $2$-adjunctions above is given by
%
\pu
\begin{equation}
	\diag{compositionof2adjunctionsseconddiagrametacommalaxcomma}
\end{equation}
where  $\alpha = \left(\id _ {\eta_ y ^\Leftarrow }\ast \eta \ast \id _{\eta _ y \! }\right)\cdot \rho $, and we denote by $\RRRR = ( R, v, \alpha ) $
the $2$-monad induced by $\check{F}\dashv \eta_y^\Leftarrow\circ\check{G}$.

\begin{rem}[$\alpha $]
Given an object $\left(x, b\right)\in \BBB / y   $,
$$\alpha _{\left(x, b\right) } : \left(x, b\right)\to \eta_ y ^\Leftarrow \check{G}\check{F} \left(x, b\right)    $$
is defined by the unique morphism $\alpha _b : w\to GF (b)\downarrow \eta _y  $ in $\BBB$ such that
the equations
%
\pu
\begin{equation}
	\diag{definitionofalphalinha}
\end{equation}
\begin{equation*}
	\left(GF(b)\right)^\Rightarrow  (\eta _ y )\, \cdot \, \alpha _b  = \eta _w \qquad\mbox{and}\qquad
	\eta _ y ^\Leftarrow  \left( GF(b) \right) \,\cdot\, \alpha _b  = b
\end{equation*}
hold.
\end{rem}
\begin{rem}
	The composition of $\varepsilon _ {F(y)}!\circ \check{F} $ with $\eta_ y !$ is
	given by $\check{F} $. More precisely, the diagrams
%
\pu
%
\pu
$$ \diag{basicliftingtwoadjunctionlaxcommatwocategoryylefttwoadjoint}\qquad \diag{basicliftingtwoadjunctionpullbacktwocategoryylefttwoadjoint} $$	
	commute.
\end{rem}

As direct consequences of the main results of Section \ref{compositionof2adjunctionssection}, we get the following corollaries.

\begin{coro}[Simplicity~\cite{MR3545937}]
	Let $G$ be pre-Kock-Z\"{o}berlein. The $2$-adjunction $$\left( F\dashv G, \varepsilon , \eta \right) : \AAA\to \BBB  $$ is simple if, and only if, for each $y\in \BBB $,
	$$\id _ T \ast\alpha  \dashv \mu \cdot \left(\id _ T \ast \delta \ast \id _ T\right)  $$
	in which $\left(\id _ T \ast\alpha \right)  $ is pointwise defined
	by $\left(\id _ T \ast\alpha \right) _b: = T( \alpha _ {(x,b)} ) $, and $\mu\cdot \left(\id _ T \ast \delta \ast \id _ T\right) $ is pointwise defined by
	\begin{equation*}
	\left(\mu\cdot \left(\id _ T \ast \delta \ast \id _ T\right)\right) _ b\quad :=\quad
	\vcenter{
		\xymatrix@=3.5em{
			T\left(T(b)\downarrow \eta _ y \right)\ar[r]^-{T\left( \delta _{T(b) } \right)  }
			\ar@{}[rd]|-{\xLeftarrow{T\left( \liva ^{ T(b)\downarrow \eta _ y }\right)   }}
			\ar[d]_-{T\left(\eta _ y ^\Leftarrow \left(T(b) \right)\right)  }
			& TT(x)
			\ar[d]|-{TT(b) } \ar[r]^-{\mu _ {x }}
			& T(x)\ar[d]^-{T(b)}
			\\
			T(y)
			\ar[r]_-{ T(\eta _y ) }
			&
			TT(y)
			\ar[r]_-{\mu_y}
			&
			T(y)
			\ar@{{}{ }{}}[lu]|-{=}
	} }
	\end{equation*}
\end{coro}
\begin{proof}
	The result follows from Corollary \ref{followsfromthepointwiseproperty} and Theorem \ref{simplicityofFGwrtJH}.
\end{proof}

\begin{coro}
Assume that $F\dashv G $ is lax idempotent. We have that $F\dashv G $
is simple provided that, for each $y\in \BBB $,  		
$\eta_y^{\Leftarrow }T\delta \check{G} $ or $ F \delta T $ is invertible.
\end{coro}
\begin{proof}
	It follows from Corollary \ref{followsfromthepointwiseproperty} and Corollary \ref{corollarytrivialforlaxidempotent}.
\end{proof}

\begin{coro}[$2$-admissibility]\label{oneofthemaincorollaries}
	Assume that $G$ is pre-Kock-Z\"{o}berlein. The $2$-adjunction $\left( F\dashv G, \varepsilon , \eta \right) : \AAA\to \BBB  $
	is $2$-admissible if and only if it is simple and, for every object $y\in \BBB$ and every object
	$a: w\to F(y)   $ of $\AAA // F(y)  $,
	the morphism defined by
	$$\xymatrix@=3.5em{
		F\left(G(a)\downarrow \eta _ y \right) \ar[d]_-{F\left(\eta _ y ^\Leftarrow \left(G(a) \right)\right)  }
		\ar[r]^-{F\left( \delta _{G(a) } \right)  }
		\ar@{}[rd]|-{\xLeftarrow{F\left( \liva ^{ G(a)\downarrow \eta _ y }\right)   }}
		&FG(w)\ar[d]|-{FG(a) } \ar[r]^-{\varepsilon _ {w }}  &
		w
		\ar[d]^-{a}
		\\
		F( y)
		\ar[r]_-{ F(\eta _ y) }	
		&
		FGF(y)\ar[r]_-{\varepsilon _ {F(y)} }
		&
		F(y)
		\ar@{{}{ }{}}[lu]|-{=}
	} $$
	in $\AAA //F(y)  $ is a regular epimorphism, \textit{i.e.} the morphism defined by
	$$\left(\varepsilon _ w\cdot F\left( \delta _{G(a) } \right) , \id _ {\varepsilon _{F(y)}}\ast F\left( \liva ^{ G(a)\downarrow \eta _ y }\right)    \right):\varepsilon _ {F(y)}!\,  \check{F}\, \eta _y ^\Leftarrow\, \check{G} (a)\to a $$
	is a regular epimorphism in $\AAA //F(y)  $.
\end{coro}
\begin{proof}
	The result follows from Corollary \ref{followsfromthepointwiseproperty} and Theorem
	\ref{admissibilitywrtJH}.
\end{proof}

\begin{coro}\label{preKockZoberleinFdashvG}
If $G$ is pre-Kock-Z\"{o}berlein then
$F\dashv G $ is $2$-admissible, provided that, for each $y\in\BBB $,  $\check{F}\delta \check{G} $ is invertible.
\end{coro}
\begin{proof}
	It follows from Corollary \ref{followsfromthepointwiseproperty} and Corollary \ref{themaincase}.
\end{proof}

It should be noted that
by Lemma \ref{equivalenciadefullreflectivelaxidempotent} we can conclude that the notion of simplicity w.r.t. the basic fibration (admissibility w.r.t. the basic fibration)
coincides with the notion of simplicity ($2$-admissibility) if $\AAA $ and $\BBB $ are locally discrete.
This shows that the notion of simplicity and $2$-admissibility can be seen as generalizations of
the classical notions of simplicity and admissibility/semi-left exact reflective functors~\cite{MR779198, MR1822890} when categories are seen as
locally discrete $2$-categories. Furthermore, Theorem \ref{southafricantheorem} shows that classical admissibility implies 2-admissibility in the presence of comma objects.
\begin{prop}\label{theoremthatfollowsfrompreservationandidentity}
	Assume that $F\dashv G $ is pre-Kock-Z\"{o}berlein, and $\AAA$ has comma objects. The $2$-adjunction $F\dashv G $ is simple ($2$-admissible) if, and only if,
	for each object $y\in\BBB $,
	the $2$-adjunction
%
\pu 	
\begin{equation}\label{idenittyonfycommaobject}
\diag{2adjunctionidentityforthesouthafricantheorem}
\end{equation}
	is simple ($2$-admissible) w.r.t.  the composite of the $2$-adjunctions
%
\pu
\begin{equation}
\diag{compositionof2adjunctionsadmissiblewrtbasifibration}
\end{equation}	
\end{prop}
\begin{proof}
By definition, $F\dashv G $ is simple ($2$-admissible)
if, and only if, for each object $y\in\BBB$, the
composition of the $2$-adjunctions of \eqref{SIMPLECOMPOSITION}
is lax idempotent (pre-Kock-Z\"{o}berlein). Since $G$ is right $2$-adjoint, it preserves comma
objects and, hence, we get that
%
\pu
%
\pu
%
\pu
\begin{equation}
\diag{theoremthatusespreservationofcommathird}
\qquad \cong \qquad \diag{theoremthatusespreservationofcommasecond}\qquad \cong\qquad
\diag{theoremthatusespreservationofcommafirstleft}
\end{equation}
in which the second $2$-natural isomorphism follows from Theorem \ref{relationofchangeofbasecomma}. By the definitions of simplicity and $2$-admissibility (see Definitions \ref{twoadmissibledefinition} and \ref{SIMPLEtwoadjunction}), the proof is complete.	
\end{proof}

\begin{theo}\label{southafricantheorem}
Provided that $\AAA $ has comma objects, if $\left(F\dashv G\right): \AAA\to \BBB $ is 	admissible w.r.t. the basic fibration, then it is $2$-admissible.
\end{theo}
\begin{proof}
By Theorem \ref{laxidempotentcoherencelaxcommacomma}, the $2$-functor
$\id_{F(y)}^{\Leftarrow } $ (the right $2$-adjoint of
\eqref{idenittyonfycommaobject})
is a pre-Kock-Z\"{o}berlein $2$-functor for every $y\in \BBB $.	
	
If $F\dashv G $ is admissible w.r.t. the basic fibration,
we get that, for every $y\in\BBB$,  $\eta _ y ^\ast \circ \check{G} $
is full reflective. Therefore
$\eta _ y ^\ast \circ \check{G} \circ \id_ {F(y)} ^{\Leftarrow } $
is a pre-Kock-Z\"{o}berlein $2$-functor by Corollary \ref{themaincase}. By Proposition \ref{theoremthatfollowsfrompreservationandidentity}, this means that
$F\dashv G $ is $2$-admissible.
\end{proof}

\section{Examples}\label{sectionexamplesexamples}
The references \cite{MR3545937, MR3708821} provide several examples of simple $2$-adjunctions/monads. In this section,
we give examples of $2$-admissible $2$-adjunctions which, in particular, are also examples of simple $2$-adjunctions.

Our first example of $2$-admissible $2$-adjunction is the identity. The result below follows directly from Theorem \ref{laxidempotentcoherencelaxcommacomma}.

\begin{lem}\label{firstexampleofadmissible}
	Let $\AAA $ be any $2$-category with comma objects. The $2$-adjunction $\id_\AAA \dashv \id _\AAA $ is $2$-admissible.
\end{lem}

Of course, the identity is also an example of admissible $2$-functor w.r.t. the basic fibration.
 Moreover, by Theorem \ref{southafricantheorem}, examples of admissible $2$-functors w.r.t. the basic fibrations
give us a wide class of examples of $2$-admissible $2$-functors.

\begin{theo}
Let $\ord $ be the $2$-category of preordered sets, and $\cat $ the $2$-category of small categories.
The inclusion $2$-functor $\ord\to\cat $ has a left $2$-adjoint and it is admissible w.r.t. the basic fibration (and, hence, also $2$-admissible).
\end{theo}

\begin{proof}
It is known that the underlying adjunction is admissible (w.r.t. the basic fibration)~\cite{MR2075604}. Since $\cat $ is a complete $2$-category, we get that the $2$-adjunction is admissible w.r.t. the basic fibration. 	
\end{proof}

Free cocompletions of $2$-categories also give us a good source for examples of admissibility w.r.t. the basic fibration. In particular, the most basic cocompletion is the free addition of the initial object.

\begin{theo}
Let $\AAA $ be a $2$-category with pullbacks and an initial object $0$. We denote by $\overline{\AAA } $ the free addition of an initial object. If $\AAA (-, 0) : \AAA ^\op\to\Cat $
is constantly equal to the empty category, the canonical $2$-functor
$$G:\AAA\to \overline{\AAA } $$
is admissible w.r.t. the basic fibration (and, hence, if $\AAA $ has comma objects, it is $2$-admissible as well).
\end{theo}
\begin{proof}
In fact $\AAA\to \overline{\AAA } $ has a left $2$-adjoint if and only if $\AAA $ has initial object.
Moreover, provided that $\AAA $ has initial object, we denote by $\eta $ the unit of this $2$-adjunction and by $\overline{0} $
the initial object freely added.

We have that
$\eta _x $ is invertible whenever $x\neq \overline{0} $. Therefore, in this case,
$$\eta _{x} ^\ast\circ \check{G} :  \AAA /x \to \overline{\AAA } /x $$
is fully faithful.

Moreover, $\eta _{\overline{0}} ^\ast\circ \check{G} :  \AAA /0 \to \overline{\AAA } /\overline{0} $
is clearly an isomorphism, since $\AAA /0$ and $\overline{\AAA } /\overline{0}$ are both empty.

This completes the proof
that $G$ is admissible w.r.t. the basic fibration and, hence, $2$-admissible provided that it has comma objects.
\end{proof}

Another example is the free cocompletion of a $2$-category under (finite) coproducts.

\begin{defi}\label{cocompletionunderfinitecoproducts}
Let $\AAA $ be a $2$-category. We define the $2$-category $\Famf\left(\AAA \right)  $ as follows.
The objects of $\Famf\left(\AAA \right)  $ are finite families of objects of $\AAA $, which
can be seen as (possibly empty) lists of objects
$$\left( x_1, \ldots , x_n\right). $$
In this case, a morphism $\left( x_1, \ldots , x_n \right)\to \left( y_1, \ldots , y_m\right)$
is a list $t = \left( t_0, \ldots , t_n \right) $
in which
$$t_0 :\left\{ 1, \ldots , n \right\}  \to \left\{ 1, \ldots , m \right\} $$
is a function,
and, for $j>0 $,
$$ t_j : x_j\to y _{t_0(j) } $$
is a morphism of $\AAA $. The composition and, hence, the identities are defined pointwise.
Finally, given morphisms
$$t = \left( t_0, \ldots , t_n \right),  t' = \left( t'_0, \ldots , t'_n \right): \left( x_1, \ldots , x_n \right)\to \left( y_1, \ldots , y_m\right)$$
of $\Famf\left(\AAA \right)  $, there is no $2$-cell $t\Rightarrow t' $, provided that $t_0\neq t_0' $. Otherwise, a
$2$-cell $\tau : t\Rightarrow t' $ is a finite family of $2$-cells
$$ \left( \tau _ j : t_j \Rightarrow t'_j:  x_j\to y_ {t_0 (j) }  \right)_{j\in\left\{ 1,\ldots , n\right\} } $$
of $\AAA $. The horizontal and vertical compositions are again defined pointwise.

\end{defi}
There is an obvious full faithful $2$-functor
$ I_\AAA: \AAA \to \Famf\left(\AAA \right)$
which takes each object $x$ to the family $(x) $. As observed above, the $2$-category $ \Famf\left(\AAA\right) $ is the \textit{free cocompletion }of $\AAA $
under finite coproducts. In particular, we have:

\begin{prop}
	The fully faithful $2$-functor
	$$I_\AAA : \AAA \to \Famf\left(\AAA \right)$$
has a left $2$-adjoint if and only if $\AAA $ has finite coproducts. In this case,
the left $2$-adjoint is given by the coproduct. More precisely, a $2$-cell
$$ \left(\tau _ 1, \ldots , \tau _n\right): (t_0, \ldots , t_n)\Longrightarrow (t_0', \ldots , t_n') : (x_1, \ldots , x_n)\to (y_1, \ldots , y_m) $$
in $ \Famf\left(\AAA \right) $ is taken to the unique $2$-cell
%
\pu 	
\begin{equation}\label{definingtheimageofthelefttwoadjointinfamlabelsecond}
\diag{definingtheimageofthelefttwoadjointinfamsecond}
\end{equation}
induced by the $2$-cells
%
\pu 	
\begin{equation}\label{definingtheimageofthelefttwoadjointinfamlabel}
\left(\diag{definingtheimageofthelefttwoadjointinfam}\right)_{i\in\left\{1,\ldots , n\right\} }
\end{equation}
in which the second arrows are the components of the universal cocone that gives the coproduct.
\end{prop}

\begin{rem}
	If we replace \textit{finite families} with \textit{arbitrary families} in
	 Definition \ref{cocompletionunderfinitecoproducts},
	we get the concept of  $\mathsf{Fam}\left(\AAA \right)  $ which corresponds to the free cocompletion of $\AAA $ under coproducts. 	
\end{rem}

We say that a $2$-category \textit{$\AAA $ has finite limits }if it has finite products, pullbacks and comma objects.
The well-known notion of extensive category has an obvious (strict) $2$-dimensional analogue.
In order to simplify the hypothesis on completion of the $2$-category $\AAA $, we are going to consider lextensive $2$-categories.

\begin{defi}[Lextensive $2$-category]\label{definitionofextensivetwocategory}
A $2$-category $\AAA $  is lextensive if it has finite limits and coproducts, and, for every finite
family of objects $\left( y_1, \ldots , y_n\right) $, the $2$-functor
\begin{eqnarray*}
	\displaystyle\prod _{j=1}^n \AAA/y_j &\to & \AAA / \coprod  _{j=1}^n  y_j \\
	\left(a_j : w_j\to y_j   \right)_{j\in\left\{1, \ldots n\right\}} &\mapsto &\coprod_{j=1}^n a_j
\end{eqnarray*}	
defined pointwise by the coproduct is a ($\Cat $-)equivalence.
\end{defi}

\begin{theo}
Let $\AAA $ be a lextensive $2$-category. We consider the $2$-adjunction
%
\pu
$$ \diag{basicadjunctionoffam} $$ 	
in which the right $2$-adjoint is the canonical inclusion. For each
finite family $Y=\left( y_j \right) _ {j\in\left\{ 1, \ldots , n\right\} } $
of objects in $\AAA $, there is a (canonical) $2$-natural isomorphism
%
\pu
\begin{equation}\label{diagramwhichprovestheadmissibilityofthecaseofFam}
\diag{basiccommutativitydiagramoftheextensivity}
\end{equation}
\end{theo}
\begin{proof}
The equivalence $2$-functor
$$\displaystyle\prod_{j=1}^n\Famf\left(\AAA/y_j\right)\to\Famf\left({\AAA}\right)/\left({y}_j\right)_{j\in\left\{{1},\ldots{,}n\right\}}
$$
is such that each object
$$A= \left( \left( a_{(1,1)}, \ldots, a_{(1,m_1)} \right), \ldots , \left( a_{(n,1)}, \ldots, a_{(n,m_n)}\right)\right)
$$
is taken to
$$ t^A = \left( t^A_ l \right)_ {l\in \left\{ 0, \left(1,1\right),\ldots , \left(1,m_1\right), \ldots , \left(n,m_n\right)\right\} }$$
in which $ t^A_0(j,k) := j $ and $t^A_{(j,k)}:= a_{(j,k)} $. The action on morphisms and $2$-cells is then
pointwise defined.
\end{proof}

\begin{coro}
Let $\AAA $ be a lextensive $2$-category. The $2$-functor $I_\AAA : \AAA \to \Famf\left( \AAA \right)$
is admissible w.r.t. the basic fibration and, hence, $2$-admissible.
\end{coro}
\begin{proof}
In fact, since products of fully faithful $2$-functors are fully faithful,  we get that
$\eta _ Y^\ast I_ \AAA $ is fully faithful by the $2$-natural isomorphism \eqref{diagramwhichprovestheadmissibilityofthecaseofFam}.	
\end{proof}

\begin{rem}
Definition \ref{definitionofextensivetwocategory} has an obvious infinite analogue, the definition of \textit{infinitary lextensive $2$-category}.
For an infinitary lextensive $2$-category $\AAA $, we have an analogous result w.r.t. $\mathsf{Fam}\left( \AAA \right) $. More precisely,
$$I_\AAA : \AAA  \to \mathsf{Fam}\left( \AAA \right) $$
is admissible w.r.t. the basic fibration (and, hence, $2$-admissible) whenever $\AAA $ is infinitary extensive.
\end{rem}

\section*{Acknowledgments}
We express our gratitude to George Janelidze for kindly hosting us at the University of Cape Town in Nov/2019 and for providing us with valuable suggestions on admissible functors. We also extend our thanks to Marino Gran and Tim Van der Linden for hosting us at the UCL, Louvain la Neuve, in May/2018, where we commenced this project. Lastly, we appreciate the anonymous referee's input, including suggestions and key bibliographic references.

\bibliographystyle{plain-abb}
\bibliography{references}

\begin{thebibliography}{10}

\bibitem{MR1007911}
R.~Blackwell, G.~M. Kelly, and A.~Power.
\newblock Two-dimensional monad theory.
\newblock {\em J. Pure Appl. Algebra}, 59(1):1--41, 1989.

\bibitem{MR1313497}
F.~Borceux.
\newblock {\em Handbook of categorical algebra. 2}, volume~51 of {\em
  Encyclopedia of Mathematics and its Applications}.
\newblock Cambridge University Press, Cambridge, 1994.

\bibitem{MR1822890}
F.~Borceux and G.~Janelidze.
\newblock {\em Galois theories}, volume~72 of {\em Cambridge Studies in
  Advanced Mathematics}.
\newblock Cambridge University Press, Cambridge, 2001.

\bibitem{zbMATH05036792}
M.~Bunge and J.~Funk.
\newblock {\em Singular coverings of toposes}, volume 1890 of {\em Lect. Notes
  Math.}
\newblock Berlin: Springer, 2006.

\bibitem{MR1397399}
A.~Carboni, G.~Janelidze, and A.~Magid.
\newblock A note on the {G}alois correspondence for commutative rings.
\newblock {\em J. Algebra}, 183(1):266--272, 1996.

\bibitem{MR779198}
C.~Cassidy, M.~H\'{e}bert, and G.~M. Kelly.
\newblock Reflective subcategories, localizations and factorization systems.
\newblock {\em J. Austral. Math. Soc. Ser. A}, 38(3):287--329, 1985.

\bibitem{MR3708821}
M.~M. Clementino and I.~L\'{o}pez-Franco.
\newblock Lax orthogonal factorisations in monad-quantale-enriched categories.
\newblock {\em Log. Methods Comput. Sci.}, 13(3):Paper No. 32, 16, 2017.

\bibitem{zbMATH07249998}
M.~M. Clementino and I.~L{\'o}pez~Franco.
\newblock Lax orthogonal factorisations in ordered structures.
\newblock {\em Theory Appl. Categ.}, 35:1379--1423, 2020.

\bibitem{MR3275274}
M.~Clementino, D.~Hofmann, and A.~Montoli.
\newblock Covering morphisms in categories of relational algebras.
\newblock {\em Appl. Categ. Structures}, 22(5-6):767--788, 2014.

\bibitem{MR3545937}
M.~Clementino and I.~L\'{o}pez~Franco.
\newblock Lax orthogonal factorisation systems.
\newblock {\em Adv. Math.}, 302:458--528, 2016.

\bibitem{MR0280560}
E.~Dubuc.
\newblock {\em Kan extensions in enriched category theory}.
\newblock Lecture Notes in Mathematics, Vol. 145. Springer-Verlag, Berlin-New
  York, 1970.

\bibitem{MR1261589}
R.~Gordon, A.~Power, and R.~Street.
\newblock Coherence for tricategories.
\newblock {\em Mem. Amer. Math. Soc.}, 117(558):vi+81, 1995.

\bibitem{MR0371990}
J.~Gray.
\newblock {\em Formal category theory: adjointness for {$2$}-categories}.
\newblock Lecture Notes in Mathematics, Vol. 391. Springer-Verlag, Berlin-New
  York, 1974.

\bibitem{zbMATH03305157}
J.~W. Gray.
\newblock Fibred and cofibred categories.
\newblock Proc. {Conf}. {Categor}. {Algebra}, {La} {Jolla} 1965, 21-83, 1966.

\bibitem{MR0249483}
J.~Gray.
\newblock The categorical comprehension scheme.
\newblock In {\em Category {T}heory, {H}omology {T}heory and their
  {A}pplications, {III} ({B}attelle {I}nstitute {C}onference, {S}eattle,
  {W}ash., 1968, {V}ol. {T}hree)}, pages 242--312. Springer, Berlin, 1969.

\bibitem{MR3673245}
D.~Hofmann and L.~Sousa.
\newblock Aspects of algebraic algebras.
\newblock {\em Log. Methods Comput. Sci.}, 13(3):Paper No. 4, 25, 2017.

\bibitem{MR1061480}
G.~Janelidze.
\newblock Pure {G}alois theory in categories.
\newblock {\em J. Algebra}, 132(2):270--286, 1990.

\bibitem{zbMATH01024334}
G.~Janelidze and G.~M. Kelly.
\newblock The reflectiveness of covering morphisms in algebra and geometry.
\newblock {\em Theory Appl. Categ.}, 3:132--159, 1997.

\bibitem{MR1245796}
G.~Janelidze, D.~Schumacher, and R.~Street.
\newblock Galois theory in variable categories.
\newblock {\em Appl. Categ. Structures}, 1(1):103--110, 1993.

\bibitem{MR1173011}
G.~Janelidze and W.~Tholen.
\newblock How algebraic is the change-of-base functor?
\newblock In {\em Category theory ({C}omo, 1990)}, volume 1488 of {\em Lecture
  Notes in Math.}, pages 174--186. Springer, Berlin, 1991.

\bibitem{zbMATH03271557}
G.~M. Kelly.
\newblock Monomorphisms, epimorphisms, and pull-backs.
\newblock {\em J. Aust. Math. Soc.}, 9:124--142, 1969.

\bibitem{MR0360749}
G.~M. Kelly.
\newblock Doctrinal adjunction.
\newblock In {\em Category {S}eminar ({P}roc. {S}em., {S}ydney, 1972/1973)},
  pages 257--280. Lecture Notes in Math., Vol. 420. Springer, Berlin, 1974.

\bibitem{MR1476422}
G.~M. Kelly and S.~Lack.
\newblock On property-like structures.
\newblock {\em Theory Appl. Categ.}, 3:No. 9, 213--250, 1997.

\bibitem{MR0357542}
G.~M. Kelly and R.~Street.
\newblock Review of the elements of {$2$}-categories.
\newblock {\em Category {S}eminar ({P}roc. {S}em., {S}ydney, 1972/1973)}, pages
  75--103. Lecture Notes in Math., Vol. 420, 1974.

\bibitem{MR1359690}
A.~Kock.
\newblock Monads for which structures are adjoint to units.
\newblock {\em J. Pure Appl. Algebra}, 104(1):41--59, 1995.

\bibitem{MR2854177}
S.~Lack and M.~Shulman.
\newblock Enhanced 2-categories and limits for lax morphisms.
\newblock {\em Adv. Math.}, 229(1):294--356, 2012.

\bibitem{2016arXiv160703087L}
F.~Lucatelli~Nunes.
\newblock On lifting of biadjoints and lax algebras.
\newblock {\em Categories and General Algebraic Structures with Applications},
  9(1):29--58, 2018.

\bibitem{2016arXiv160604999L}
F.~Lucatelli~Nunes.
\newblock Pseudo-{K}an extensions and descent theory.
\newblock {\em Theory Appl. Categ.}, 33:No. 15, 390--444, 2018.

\bibitem{arXiv:1711.02051}
F.~{Lucatelli Nunes}.
\newblock Pseudoalgebras and non-canonical isomorphisms.
\newblock {\em Appl. Categ. Struct.}, 27(1):55--63, 2019.

\bibitem{2019arXiv190201225L}
F.~{Lucatelli Nunes}.
\newblock Semantic factorization and descent.
\newblock {\em Appl. Categ. Struct.}, 30(6):1393--1433, 2022.

\bibitem{MR1712872}
S.~Mac~Lane.
\newblock {\em Categories for the working mathematician}, volume~5 of {\em
  Graduate Texts in Mathematics}.
\newblock Springer-Verlag, New York, second edition, 1998.

\bibitem{MR2056584}
J.~MacDonald and M.~Sobral.
\newblock Aspects of monads.
\newblock In {\em Categorical foundations}, volume~97 of {\em Encyclopedia
  Math. Appl.}, pages 213--268. Cambridge Univ. Press, Cambridge, 2004.

\bibitem{MR1432190}
F.~Marmolejo.
\newblock Doctrines whose structure forms a fully faithful adjoint string.
\newblock {\em Theory Appl. Categ.}, 3:No.\ 2, 24--44, 1997.

\bibitem{MR3207214}
A.~Montoli, D.~Rodelo, and T.~Van~der Linden.
\newblock A {G}alois theory for monoids.
\newblock {\em Theory Appl. Categ.}, 29:No. 7, 198--214, 2014.

\bibitem{MR1040947}
A.~Power.
\newblock A {$2$}-categorical pasting theorem.
\newblock {\em J. Algebra}, 129(2):439--445, 1990.

\bibitem{MR0299653}
R.~Street.
\newblock The formal theory of monads.
\newblock {\em J. Pure Appl. Algebra}, 2(2):149--168, 1972.

\bibitem{MR0401868}
R.~Street.
\newblock Limits indexed by category-valued {$2$}-functors.
\newblock {\em J. Pure Appl. Algebra}, 8(2):149--181, 1976.

\bibitem{MR558494}
R.~Street, W.~Tholen, M.~Wischnewsky, and H.~Wolff.
\newblock Semitopological functors. {III}. {L}ifting of monads and adjoint
  functors.
\newblock {\em J. Pure Appl. Algebra}, 16(3):299--314, 1980.

\bibitem{MR2075604}
J.~Xarez.
\newblock Separable morphisms of categories via preordered sets.
\newblock In {\em Galois theory, {H}opf algebras, and semiabelian categories},
  volume~43 of {\em Fields Inst. Commun.}, pages 543--549. Amer. Math. Soc.,
  Providence, RI, 2004.

\end{thebibliography}

\pu

\end{document}
